\documentclass[a4paper,10pt]{article}
\usepackage[utf8]{inputenc}     
\usepackage{lmodern}
\usepackage{graphicx}

\usepackage[english, french]{babel}
\usepackage[french]{babel}
\usepackage{verbatim}
\usepackage{graphicx}
\usepackage{algorithmicx}
\usepackage{xcolor, soul}
\sethlcolor{red}
\usepackage{tikz}

\usepackage{subfigure}
\usepackage{amsmath, amsfonts, dsfont, amssymb}
\usepackage{framed}
\usepackage{mdframed}
\usepackage[amsmath, standard, thmmarks]{ntheorem}
\usepackage[noend]{algpseudocode}
\usepackage{hyperref}
\usepackage{url}
\usepackage{enumerate}
\usepackage[margin=3cm]{geometry}
\usepackage{color}
\usepackage{mathrsfs} 
\usepackage[french]{minitoc} 

\usepackage{appendix}

\usepackage{fancyhdr}
\usepackage{authblk}

\newcommand{\calA}{\mathcal{A}}

\newcommand{\calC}{\mathcal{C}}

\newcommand{\calO}{\mathcal{O}}

\newcommand{\calT}{\mathcal{T}}

\newcommand{\fre}{\mathfrak{e}}

\newcommand{\frt}{\mathfrak{t}}

\newcommand{\frA}{\mathfrak{A}}

\newcommand{\frT}{\mathfrak{T}}

\newcommand{\bbE}{\mathbb{E}}

\newcommand{\bbN}{\mathbb{N}}

\newcommand{\bbP}{\mathbb{P}}

\newcommand{\bbR}{\mathbb{R}}

\newcommand{\bbZ}{\mathbb{Z}}

\newcommand{\bfh}{\mathbf{h}}
\newcommand{\bfe}{\mathbf{e}}
\newcommand{\vfre}{\vec{\fre}}
\newcommand{\bfd}{\mathbf{d}}
\newcommand{\bfT}{\mathbf{T}}
\newcommand{\bddelta}{\boldsymbol{\delta}}
\newcommand{\bdphi}{\boldsymbol{\phi}}
\newcommand{\bdcalT}{\boldsymbol{\mathcal{T}}}
\newcommand{\bdtheta}{\boldsymbol{\theta}}
\newcommand{\bdTheta}{\boldsymbol{\Theta}}

\newcommand{\bdmu}{\boldsymbol{\mu}}
\newcommand{\bdsigma}{\boldsymbol{\sigma}}
\newcommand{\bdSigma}{\boldsymbol{\Sigma}}
\newcommand{\bdgamma}{\boldsymbol{\gamma}}

\newcommand{\bdcalV}{\boldsymbol{\mathcal{V}}}
\newcommand{\bdcalH}{\boldsymbol{\mathcal{H}}}
\newcommand{\bdcalS}{\boldsymbol{\mathcal{S}}}

\newcommand{\bfv}{\mathbf{v}}
\newcommand{\bfw}{\mathbf{w}}

\newcommand{\bdcalC}{\boldsymbol{\mathcal{C}}}
\newcommand{\bdfrT}{\boldsymbol{\frT}}
\newcommand{\tildebdfrT}{\boldsymbol{\widetilde{ \frT}}}
\newcommand{\tildebdcalC}{\boldsymbol{\widetilde{ \calC}}}

\newcommand{\bdcalR}{\boldsymbol{\mathcal{R}}}
\newcommand{\bdomega}{\boldsymbol{\omega}}

\newcommand{\bfB}{\mathbf{B}}
\newcommand{\bdfrA}{\boldsymbol{\frA}}
\newcommand{\bfH}{\mathbf{H}}
\newcommand{\bfX}{\mathbf{X}}
\newcommand{\bfN}{\mathbf{N}}

\numberwithin{equation}{section}

\DeclareSymbolFont{symbolsC}{U}{txsyc}{m}{n}
\SetSymbolFont{symbolsC}{bold}{U}{txsyc}{bx}{n}
\DeclareFontSubstitution{U}{txsyc}{m}{n}
\DeclareMathSymbol{\multimapboth}{\mathrel}{symbolsC}{"13}

\newcommand{\theoremname}{Theorem}
\addto\captionsfrench{\renewcommand{\theoremname}{Th\'eor\`eme}}
\addto\captionsenglish{\renewcommand{\theoremname}{Theorem}}

\makeatletter
\newtheoremstyle{MyNonumberplain}%
 {\item[\theorem@headerfont\hskip\labelsep ##1\theorem@separator]}
  {\item[\theorem@headerfont\hskip\labelsep ##3\theorem@separator]}
\makeatother
\theoremstyle{MyNonumberplain}
\theorembodyfont{\upshape}
\newtheorem{prove}{Proof}

\def\myhfill{
  \parfillskip=0pt
  \widowpenalty=10000
  \displaywidowpenalty=10000
  \finalhyphendemerits=0
  \unskip\nobreak\null\hfil\penalty50
  \hskip2em\null\hfill
}
 
\def\qedsymb{\ensuremath\square}
\def\qedd{\myhfill\qedsymb\par}

\title{Percolation on supercritical causal triangulations}
\author{David Corlin Marchand\thanks{Laboratoire de Math\'ematiques Rapha\"{e}l Salem, UMR CNRS 6085, Universit\'e de Rouen Normandie, avenue de l’Universit\'e, Technop\^{o}le du Madrillet, 76801 Saint-\'{E}tienne-du-Rouvray,
France and Institut Mines Télécom Nord Europe, Cit\'e scientifique, rue Guglielmo Marconi,
59650 Villeneuve-d’Ascq, France. 

\emph{Email address:}
\href{mailto:david.corlin.marchand@gmail.com}{david.corlin.marchand@gmail.com}.}}

\begin{document}

\selectlanguage{english}
\maketitle

\begin{abstract}
We study oriented percolation on random causal triangulations, those are random planar graphs obtained roughly speaking by adding horizontal connections between vertices of an infinite tree. When the underlying tree is a geometric Galton--Watson tree with mean~$m>1$, we prove that the oriented percolation undergoes a phase transition at~$p_c(m)$, where
$$ p_c(m) = \frac{\eta}{1+\eta} \quad \mbox{ with } \quad \eta = \frac{1}{m+1} \sum_{n \geq 0} \frac{m-1}{m^{n+1}-1}.$$
We establish that strictly above the threshold $p_c(m)$, infinitely many infinite components coexist in the map. This is a typical percolation result for graphs with a hyperbolic flavour.
We also demonstrate that large critical oriented percolation clusters converge after rescaling towards the Brownian continuum random tree. The proof is based on a Markovian exploration method, similar in spirit to the peeling process of random planar maps.
\end{abstract}

\section{Introduction and main results}

\begin{figure}[h!]
    \centering 
    \includegraphics[scale=0.65]{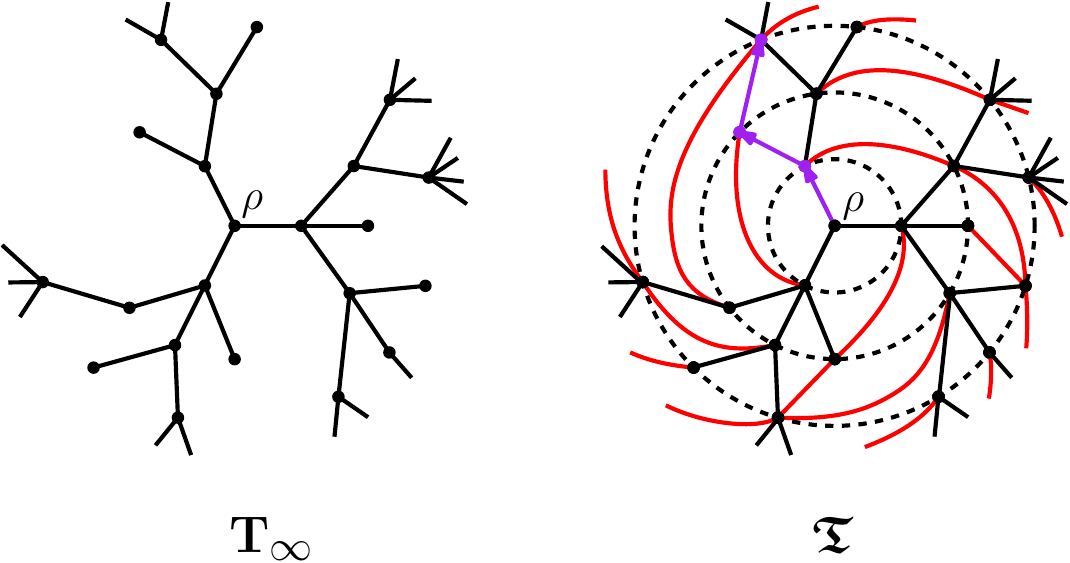}
    \caption{A piece of an infinite supercritical Galton--Watson tree $\bfT_\infty$, rooted in $\rho$, and its associated causal triangulation $\bdfrT$. In dashed lines are the horizontal edges connecting at each level consecutive vertices. They do not play any role in our work. Red edges are added to triangulate the map. Edges in purple illustrate the orientation---away from the root---chosen in our model.}
\label{fig:defcausaltriangulations}
\end{figure}

In these pages we study oriented percolation on a random "hyperbolic" type of planar graph which is constructed as follows. Let $t$ be any infinite rooted plane tree. For any $r \geq 1$, we add horizontal edges between consecutive vertices belonging to the boundary of the ball of radius $r$, so that it forms a cycle. Then, we triangulate each face by connecting vertices on its bottom side---except the rightmost one---to the top-right vertex. 
See Figure~\ref{fig:defcausaltriangulations} for an illustration.
The planar triangulation that we get is called a \textit{causal triangulation}.

We are interested in the case where $t$ is a random supercritical Galton--Watson tree conditioned to survive, denoted by $\bfT_\infty$. Specifically, we shall focus on the case where the offspring distribution~$\bdmu$ of the tree is geometric. 
That is for some $\alpha \in (1/2,1)$, we have for every~$k \geq 0$, $$\bdmu(k)=\bdmu_{\alpha}(k):=\alpha^k (1-\alpha).$$
We recall that the mean of $\bdmu_\alpha$ is $$m:=\alpha \cdot (1-\alpha)^{-1}>1.$$
We will denote by $\bdfrT$ the random causal triangulation built from $\bfT_\infty$, and $\bbP_{\bdfrT}$ the associated probability measure. We call the model \textbf{supercritical causal triangulation}, abbreviated to~SCT. In what follows, we will write~$\mathbf{GW}_\beta$ to designate the distribution of a Galton--Watson tree with offspring law~$\bdmu_\beta$---without conditioning on the survival event---and where~$\beta$ is more generally in $(0,1)$. \\

\noindent\fbox{\parbox[c]{\textwidth}{\textbf{\quad In the rest of the paper, we assume that $\alpha  \in (1/2,1)$ is fixed. We will omit $\alpha$ in our notation most of the time, although what will be stated will often depend on it. }}} \\ \\

Causal triangulations first appeared in the physics literature \cite{ambjorn1998non}. Random versions have been recently studied from a mathematical point of view \cite{curien2020geometric,budzinski2019supercritical}. In particular, the author of \cite{budzinski2019supercritical} focuses on the causal maps derived from a supercritical (non necessarily geometric) Galton--Watson tree conditioned to survive. Causal maps are constructed from a rooted tree in the same manner as causal triangulations, except that they are not triangulated. It is shown in the mentioned paper that they have "hyperbolic" properties (anchored Gromov hyperbolicity, positive speed of the associated simple random walk, a non-trivial Poisson boundary, etc).

\paragraph{Results.} Our contribution here is to grasp the SCT model through the lens of \emph{directed percolation}. Fix a (deterministic) causal triangulation $\theta$. The edges are directed away from the root, \textbf{the horizontal connections play no role}. A Bernoulli bond percolation is applied on them: for some fixed parameter $p \in [0,1]$, each (directed) edge is declared either \textit{open} or \textit{closed} with probabilities~$p \in [0,1]$ and~$1-p$ respectively, independently of the others. We denote by $\bbP_{p,\theta}$ the distribution of the percolation process on $\theta$. Now, we introduce $\bbP_p$ the overall distribution
$$\int \bbP_{\bdfrT}(\mathrm{d} \theta) \int \mathrm{d} \bbP_{p,\theta}$$
which averages the percolation process on supercritical causal triangulations. This is an \emph{annealed} percolation distribution, while $\bbP_{p,\theta}$ is a \emph{quenched} one. 
 
In our first main result, we show the existence of a non trivial \emph{annealed} phase transition: 

\begin{theorem}\label{theorem:phasetransitionSCQ}
\NoEndMark Let $\bdcalC$ be the (directed) percolation cluster of the root in $\bdfrT$. We define $\Theta(p)$ the \emph{annealed} probability that $\bdcalC$ is infinite and set
\begin{align}
\label{theorem:exactexpressioncriticalpointdirectedpercolationforSCQ} p_c=p_c(m) := \frac{\eta}{1+\eta} \quad \mbox{ with } \quad \eta = \frac{1}{m+1} \sum_{n \geq 0} \frac{m-1}{m^{n+1}-1}.
\end{align}
Then:
\begin{align}
\label{eq:probaclusterorigininfinite}
\Theta(p) = \left\{
    \begin{array}{ll}
        0 & \mbox{if } p \leq p_c \\
	 >0 & \mbox{if } p > p_c.
    \end{array}
\right.
\end{align} 
Furthermore, when $p_c < p <1$, there are $\bbP_p$-almost surely infinitely many disjoint infinite clusters in $\bdfrT$.
\end{theorem}

Since the environment is directed, the notion of disjoint clusters may seem ambiguous. Imagine indeed two clusters $\bdcalC_u$ and $\bdcalC_v$ associated to two distinct vertices $u$ and $v$, which do not share any vertex in common, but are both included in a bigger cluster emanating from a third vertex nearer to the root. We do not definitely want to consider them as disjoint. So, to be more precise, our statement is that as $p_c<p<1$, there are $\bbP_p$-almost surely infinitely many disjoint \emph{non-directed} clusters (i.e.~when edges are not directed anymore), each containing at least one infinite \emph{directed} cluster within it.

We easily derive from Theorem~\ref{theorem:exactexpressioncriticalpointdirectedpercolationforSCQ} a \emph{quenched} version of it. In a classic way, we define for any causal triangulation $\theta$:
$$p_c (\theta):=\inf\big\{p \in [0,1]: \ \bbP_{p,\theta}(\text{there exists an infinite (directed) cluster in $\theta$})=1\big\}.$$
The theorem ensures that $p_c (\theta)$ is $\bbP_{\bdfrT}$-almost surely equal to $p_c$. Also, we have $\bbP_{\bdfrT}$-almost surely  infinitely many disjoint infinite clusters in $\theta$ with probability one.

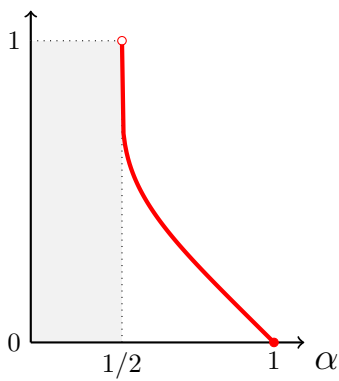
\begin{figure}[h!]
\centering
\begin{tikzpicture}[scale=4]
\fill[color=gray!10] (0.2,0) -- (0.2,1) -- (0.5,1) -- (0.5,0) -- cycle;
\draw[thick] (0.2,0) node[left] {$0$} -- (0.2,1) node[left] {$1$};
\draw[->, thick] (1,0) -- (1.1,0) node[below right] {\Large{$\alpha$}};
\draw[dotted] (0.5,0) node[below] {$1/2$} -- (0.5,1);
\draw[->, thick] (0.2,1) -- (0.2,1.1);
\draw[ultra thick, red] plot file {dataplotcriticalpointSCQ.txt};
\draw[dotted] (0.2,1) -- (0.5,1);
\filldraw[white] (0.5,1) circle (0.3pt);
\draw[red] (0.5,1) circle (0.4pt);
\draw[thick] (0.2,0) -- (1,0) node[below] {$1$};
\filldraw[red] (1,0) circle (0.4pt);
\end{tikzpicture}
\caption{Plot of the map $\alpha \mapsto p_c (\alpha)$ on the interval $(1/2,1)$.}
\label{fig:plotcriticalpointSCQ}
\end{figure}

The graph of $p_c$ as a function of $\alpha$ is plotted on Figure~\ref{fig:plotcriticalpointSCQ}. We remark that~$p_c(\alpha) \to 1$ as~$\alpha \to 1/2$. It is worth to investigate the case $\alpha=1/2$. 
The Galton--Watson tree $\bfT$ is then critical and almost surely finite. By conditioning its size to be increasingly large, an infinite random tree emerges as local limit in distribution, namely the Kesten's tree with offspring distribution~$\bdmu_{1/2}$~\cite{abraham2015introduction}. A causal triangulation can be defined from it. Furthermore, this infinite tree is the local limit as~$\alpha \to (1/2)^+$ of Kesten's trees with offspring distribution~$\bdmu_\alpha$.  The convergence  holds for the related causal triangulations as well. Now, we claim that the percolation threshold of the causal triangulation built from the Kesten's tree with reproduction law~$\bdmu_\alpha$ is equal to the~$p_c$ defined in Theorem~\ref{theorem:phasetransitionSCQ}, because the distribution of the tree is absolutely continuous with respect to that of~$\bfT_\infty$ \cite{lyons1995conceptual}. Since the critical percolation threshold is continuous with respect to the local topology, we finally deduce that the phase transition indeed degenerates, that is to say~$p_c=1$ when~$\alpha=1/2$. 
Although the percolation problem is trivial on the critical causal triangulation, its geometric and spectral properties are not and were studied in \cite{durhuus2010spectral} and \cite{curien2020geometric}. 

The annealed phase transition described in Theorem~\ref{theorem:phasetransitionSCQ} is of the same nature when bond percolation is performed on the underlying Galton--Watson tree alone. Only the value of the threshold~$p_c$ change into $m^{-1}$. It stems from the fact that in such context, the cluster of the root vertex is itself a Galton--Watson tree of mean $p  m$. The conclusion remains valid by taking for instance an infinite regular tree instead of a supercritical Galton--Watson tree. Such kind of phase transition---including a non degenerate phase $p_c < p <p_u$ where infinitely many infinite clusters coexist in the map---is typical of percolation on "hyperbolic" type of graphs---that is roughly those where the volume of balls of radius $r$ grows exponentially fast as $r \to +\infty$. See \cite{benjamini1996percolation} or \cite{hutchcroft2019percolation} for deterministic environments. 
Similar results were also established in random environment contexts, like for site percolation on hyperbolic half-planar triangulations in~\cite{ray2014geometry}, or bond percolation on hyperbolic bipartite half-planar maps in~\cite{curien2019peeling}. In the former reference, the author furthermore shows the existence of another non degenerate phase $p_u <p$, where there is exactly one (and only one) infinite cluster in the map. In our case, the phase is unfortunately trivial with~$p_u=1$, as it is also in the tree set-up mentioned above.

In a second result, we compute several annealed exponents related to percolation cluster:
~\begin{theorem} 
\label{theorem:criticalexponentstheoremSCQ}
\NoEndMark We have:
\begin{enumerate}[(i)]
\item $\bbP_{p_c}(\lvert \bdcalC \rvert  \geq n) = \Omega\big(n^{-1/2}\big)$  and $\bbP_{p_c}(\lvert \bdcalC \rvert  \geq n) =\calO\big(\log{(n)} \cdot n^{-1/2}\big)$;
\item $C^{-1} \cdot (p-p_c) \leq \Theta(p) \leq C \cdot (p-p_c)$ for any $p>p_c$ and some constant $C>0$.
\end{enumerate}
\end{theorem}

Here and later, for any two sequences $f,g : \bbN \mapsto (0,+\infty)$, we use the notation $\displaystyle{f(n)=\calO\big(g(n)\big)}$, resp.~$\displaystyle{f(n)=\Omega\big(g(n)\big)}$, if the ratio $f/g$ is bounded below, resp.~above, by some positive constant~$C>~0$. Sometimes, we also write $f(n)=o\big(g(n)\big)$ to mean that $f(n)/g(n) \xrightarrow[]{}~0$. 

We expect the second asymptotic of item (i) to hold with no logarithmic correction to the polynomial term. We insist on the fact that exponents of Theorem~\ref{theorem:criticalexponentstheoremSCQ} are annealed versions, i.e.~averaged on the random map and the percolation process. We do not know if quenched exponents exist and if their values coincide with the annealed ones. The same critical and off-critical exponents arise in the context of bond percolation on a Galton--Watson tree, or on deterministic hyperbolic graphs \cite{hutchcroft2019percolation}. We do not have knowledge of any equivalent result in other random hyperbolic models. See \cite{angel2015percolations, curien2015percolation, gorny2018geometry} or \cite{curien2019peeling} about non-hyperbolic environments, where the values of critical exponents significantly differ from those observed in our model.

In our last main result, we demonstrate that a large critical oriented percolation cluster admits the Brownian continuum random tree as scaling limit:

\begin{theorem}
\label{theorem:scalinglimitSCQ}
\NoEndMark At $p=p_c$, we have:
\begin{align}
\label{theorem:scalinglimitSCQconvergenceequality}
\big( n^{-1/2} \cdot \bdcalC \ \big\vert \ \vert \bdcalC \rvert \geq n \big) \xrightarrow[n \to +\infty]{d_{GH}} \kappa' \cdot \bdcalT_{\geq 1},
\end{align}
where $\kappa'=\kappa'(\alpha)>0$ is a positive number only depending on $\alpha$
and $\bdcalT_{\geq 1}$ is the Continuum Random Tree of mass greater than~1 \cite{aldous1991continuum, aldous1993continuum}. The convergence \eqref{theorem:scalinglimitSCQconvergenceequality} holds in distribution for the Gromov--Hausdorff distance.
\end{theorem}

Again, a comparable behaviour is witnessed in the Galton--Watson tree cleared of the additional edges going with the causal triangulation. Indeed, as~$p=p_c=m^{-1}$, the percolation cluster of the root is distributed as a critical Galton--Watson. When the latter is conditioned to have a large size, it converges to the CRT \cite{aldous1993continuum, le2005random}. The emergence of such scaling limit in a hyperbolic context is not new. See for instance~\cite{chen2017long} where the authors prove that the CRT is the limit of long Brownian bridges in hyperbolic space. For a precise definition of the Gromov--Hausdorff topology and of the metric space $\bdcalT_{\geq 1}$, we refer the reader to Section~\ref{sec:scalinglimitpercocluster}. 

\paragraph{Techniques.} The main tool to prove Theorems~\ref{theorem:phasetransitionSCQ}, \ref{theorem:criticalexponentstheoremSCQ} and \ref{theorem:scalinglimitSCQ} is a Markovian exploration of~$\bdfrT$ along the directed cluster of the origin. More precisely, we show that in a simplified "half-plane" model $ \tildebdfrT$ of~$\bdfrT$, we can explore the underlying map step-by-step using an algorithm tailored to the percolation.  This exploration yields a random walk with independent increments which roughly does the contour of our percolation cluster. This is reminiscent of the peeling process in the theory of random planar maps, and in particular of the fact that percolation on random half-planar maps can be studied  (using the peeling process) via a random walk with independent increments. See~\cite{angel2003growth,angel2005scaling} for the pioneer works of Angel on the subject and \cite{curien2019peeling} for a comprehensive treatment. With this tool at hands, a version of Theorems~\ref{theorem:phasetransitionSCQ}, \ref{theorem:criticalexponentstheoremSCQ} and \ref{theorem:scalinglimitSCQconvergenceequality} is easily established in our simplified model~$\tildebdfrT$. The results are then transferred back to $\bdfrT$.

\begin{paragraph}{Organisation of the paper}
In Section~\ref{sec:toymodelSCTHP}, we introduce the toy model of random triangulations of the half-plane, and outline its main properties. We describe a peeling procedure for these maps which exhibits a spatial Markov property. In Section~\ref{sec:directedpercoontoymodelSCTHP}, we study percolation on the toy model and prove an equivalent of our three main theorems. We deal with the SCT model in Section~\ref{sec:directedpercoonSCQ}. Potential extensions of our results are discussed in Section~\ref{sec:openproblems}.
\end{paragraph}

\begin{paragraph}{Acknowledgement}
This work was mostly carried out during my PhD at the Laboratoire de Math\'{e}matiques d'Orsay. I deeply thank Nicolas Curien for his careful reading of the earlier versions of the paper and all his valuable comments.
\end{paragraph}

\section{The supercritical causal triangulation of the half-plane}
\label{sec:toymodelSCTHP}

We present in a first section the definition of a random model of triangulations of the half-plane, which is supposed to approach SCT. The idea is that it will capture essential features of $\bdfrT$ at large distances from the root vertex.
In the next two sections, we outline some basic properties of this model, and introduce a Markovian exploration of the map.

\subsection{Definition of the model} \label{sec:deftoymodelSCTHP}

We start by introducing a family of random triangulation of the cylinder $\bbR \times [0,1]$, with vertex set~$\bbN \times \{0,1\}$. The maps are then \emph{labelled}. Any face of such a map is a triangle, which has to be either \emph{top-oriented}---if it shares an edge with the upper boundary of the cylinder---or \emph{bottom-oriented}---if it does with the lower boundary. We can rank them according to their order of apparition from left to right. In our random model, we assume that each triangle is top-oriented with probability~$\alpha$, or bottom-oriented with probability $1-\alpha$, independently of the others. Then, the triangulation can be summed up by a sequence of i.i.d.~biased Bernoulli trials of parameter~$\alpha$, where the $i$-th trial indicates whether the $i$-th leftmost triangle is top or bottom-oriented.

\begin{figure}[h!]
    \centering 
    \includegraphics[scale=0.45]{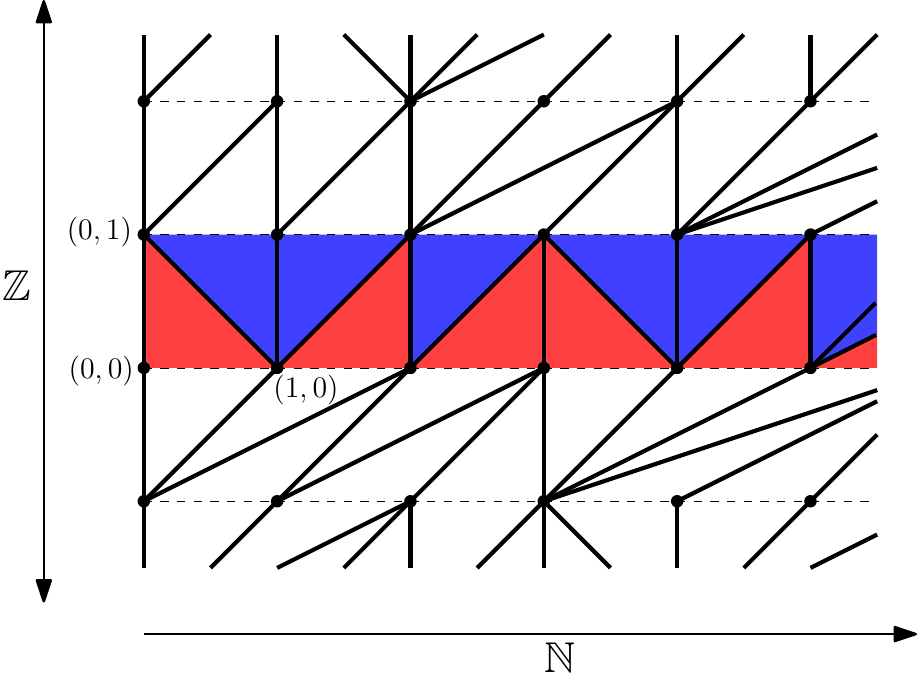}
    \caption{A portion of a causal triangulation of the half-plane $\tildebdfrT$. In the middle cylinder between heights~$0$ and $1$, top-oriented triangles are coloured in blue, bottom-oriented triangles in red.}
\label{fig:ExampleSCTHP}
\end{figure}

Consider now an i.i.d.~sequence (indexed by $\bbZ$) of such random maps. We glue them so that for any $i \in \bbZ$, the lower boundary of the $i$-th cylinder coincides with the upper boundary of the~$i-1$-th. After a simple relabelling operation of vertices, it yields a random triangulation of the half-plane~$\bbR_+ \times \bbR$, denoted by $\tildebdfrT$, with vertex set $\bbN \times \bbZ$. We call it \textbf{$\alpha$-supercritical causal triangulation of the half-plane}, abbreviated to $\alpha$-SCTHP or just SCTHP. See Figure~\ref{fig:ExampleSCTHP} for an illustration.

As in the SCT model, \textbf{horizontal edges here never play an effective role} since our aim is to study directed percolation with edges being oriented in the bottom-up direction. This is why they are drawn in dashed lines in our figures. Their presence is however helpful to understand the peeling algorithm that we will introduce later.

\subsection{Basic properties} \label{sec:basicpropertiesofthetoymodelSCTHP}

Some natural notion of direction---both vertically and horizontally---arises in SCTHP. For any vertex~$v=(i,j)$ of a SCTHP, we will designate by~$h(v)$---\textit{the height of $v$}---the second coordinate~$j$. Vertices $w=(i',j)$---with same height than $v$---shall be said \textit{on the left} of $v$ if $i' \leq i$,  or \textit{on the right} if $i'\geq i$. Given the structure of SCTHP, the vertex~$v$ is necessarily adjacent to vertices $w$ which are either such as $$h(w)=h(v)-1, \quad \text{or } \quad h(w)=h(v), \quad \text{or also} \quad h(w)=h(v)+1.$$Vertices of the first kind shall be called \textit{parents of~$v$}, those of the last kind \textit{offspring of $v$}. We remark that if $w$ is on the left, resp.~on the right, of~$v$, then its offspring are on the left, resp.~on the right, of offspring of $v$. The same holds for the parents.

Of course, the parent and offspring notions can be reversed, and it allows to exhibit some nice duality property of $\alpha$-SCTHP. Indeed, by flipping vertically the map, we transform a top-oriented triangle into a bottom-oriented, and vice versa. The $\alpha$-SCTHP becomes an $(1-\alpha)$-SCTHP.

\begin{figure}[h!]
    \centering 
    \includegraphics[scale=0.3]{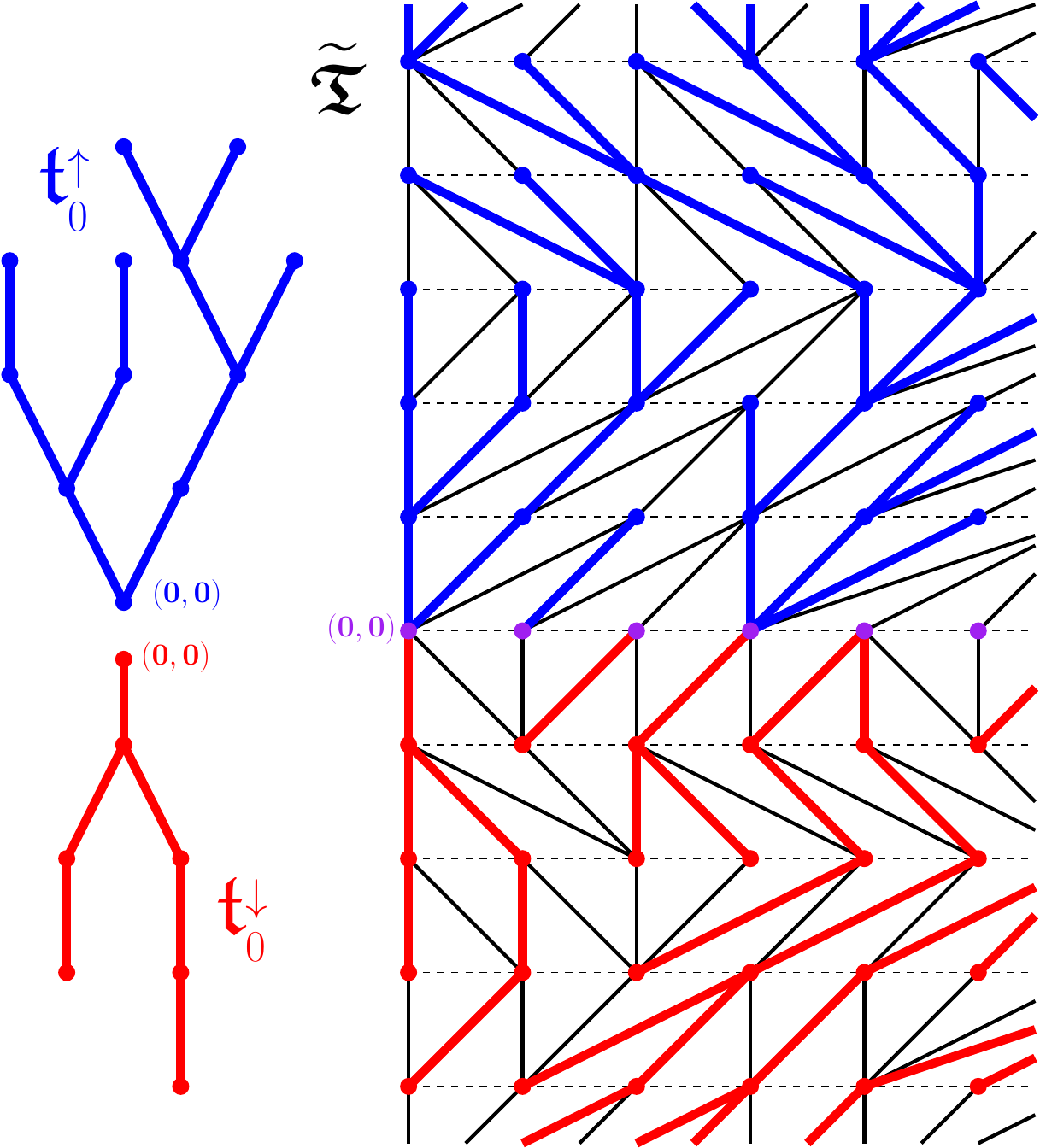}
    \caption{Ascending (in blue) and descending (in red) trees emerging from vertices of $\bbN \times \{ 0 \} \subset \tildebdfrT$. Note how large is $\frt_{(3,0)}^{\uparrow}$ while $\frt_{(2,0)}^{\uparrow}$ and $\frt_{(1,0)}^{\downarrow}$ are for instance made of a single vertex. On the left-hand side of the figure are depicted the ascending and tree associated to the origin vertex. The thin continuous or dashed black lines are the edges deterministically added to recover the whole map, in the case where the only provided data are the blue and red trees. The operation is the same than that to build $\bdfrT$ from $\bfT$ in the SCT model. See Figure~\ref{fig:defcausaltriangulations}.}
\label{fig:exampleascendingdescendingtrees}
\end{figure}

It is worthwhile to note that some tree structure underlies SCTHP:

\begin{definition}[Ascending and descending tree]
\label{def:ascendingdescendingtreeDHP}
\NoEndMark For any vertex $v$ of a SCTHP, we define its associated \emph{ascending tree} $\frt_{v}^{\uparrow}$ as the plane tree rooted in $v$, containing all its offspring except the rightmost one, and the ascending trees associated to these offspring. Conversely, we also define its associated \emph{descending tree} $\frt_{v}^{\downarrow}$ as the plane tree rooted in $v$, containing all its parents except the rightmost one, and the descending trees associated to these offspring. See Figure~\ref{fig:exampleascendingdescendingtrees} below for an illustration.
\end{definition}

At any given height $h \in \bbZ$, the collection of ascending and descending trees emerging from vertices with height $h$ is enough data to recover the whole map. Indeed, these ascending trees, resp.~descending trees, form a partition of the set of vertices with height greater than $h$, resp.~less than~$h$. Edges to add to get the associated SCTHP are deterministically deduced from the trees, in the same way that the SCT $\bdfrT$ was obtained from the tree $\bfT$. See Figure~\ref{fig:exampleascendingdescendingtrees} for an illustration. 

As with parent and offspring, the notions of ascending and descending trees switch by flipping vertically the SCTHP: the ascending tree of some vertex $v$ becomes its descending tree in the new map, and conversely.

The ascending and descending trees have actually a well-known distribution:
\begin{proposition}
\label{prop:ascendingdescendingtreesHDHPareGW}
\NoEndMark For any height $h \in \bbZ$, the ascending trees~$\big(\frt_{v}^{\uparrow}\big)_{v}$ and the descending trees~$\big(\frt_{v}^{\downarrow}\big)_{v }$ emerging from vertices with height $h$ are all mutually independant. Moreover, the ascending trees are distributed like $\mathbf{GW}_\alpha$, and the descending trees like $\mathbf{GW}_{1-\alpha}$.
\end{proposition}

\begin{prove}[Proof of Proposition~\ref{prop:ascendingdescendingtreesHDHPareGW}]
As explained in Section~\ref{sec:deftoymodelSCTHP}, triangles in a layer of the SCTHP are i.i.d.~and top-oriented with probability $\alpha$. This gives some trivial spatial Markov property: assume that we know the orientation of the $K$ first triangles (ordered from left to right), the next triangles remain i.i.d.~with the same distribution. In particular, the probability that the next $L$ first triangles are top-oriented is $\alpha^{L}=\bdmu_{\alpha}(\cdot \geq L)$. Since the number of offspring of some vertex~$v$ is equal to the number of top-oriented triangles which emerge from it, the above facts imply that it is distributed as $\bdmu_\alpha$ and independent of the corresponding quantity for the other vertices. Triangulations of the cylinders being independent, we deduce that the ascending trees~$\big(\frt_{v}^{\uparrow}\big)_{v}$ emerging from vertices with some height $h \in \bbZ$ are i.i.d.~with distribution $\mathbf{GW}_\alpha$.
Due to duality, a similar conclusion holds for the descending trees, but the offspring distribution is in this case~$\bdmu_{1-\alpha}$. Finally, the ascending and descending trees are mutually independent since the former are measurable with respect to cylinders above height $h$, while the latter are with respect to those below height~$h$. \qedd
\end{prove}

A straightforward consequence of Proposition~\ref{prop:ascendingdescendingtreesHDHPareGW} is that ascending trees in SCTHP are supercritical Galton--Watson trees---since $\sum_{k \geq 0} k \cdot \bdmu_{\alpha}(k) = \frac{\alpha}{1-\alpha}>1$ when $\alpha \in (\frac{1}{2},1)$---while descending trees are subcritical. Later on, it will be useful to have some more information on the latter. Below, we prove that the distribution of their height is exponentially tailed. We also express the first moment as an infinite series. 

\begin{proposition}
\label{prop:heightsubcriticalGWtrees}
Consider $\frt ^{\downarrow} \overset{(d)}{\sim} \mathbf{GW}_{1-\alpha}$ for $\alpha > 1/2$. We note~$\bfh(\cdot)$ the function returning the height of any rooted tree. Recall that $m:= \frac{\alpha}{1-\alpha}>1$. Then:
$$\mathbb{P}\Big(\bfh(\frt ^{\downarrow})\geq n\Big) = \frac{m-1}{m^{n+1}-1} \quad \text{and so} \quad \mathbb{E}[\bfh(\frt ^{\downarrow})] = \sum_{n \geq 1} \frac{m-1}{m^{n+1}-1}.$$
\end{proposition}

\begin{prove}[Proof of Proposition~\ref{prop:heightsubcriticalGWtrees}]
Let $(\mathbf{Z}_n)_{n \geq 0}$ be a Galton--Watson process started from one particle and with $\bdmu_{1-\alpha}$ as offspring distribution. We have:
\begin{align}
\label{relationtailheightgwprocessexpectation}
\bbE[\bfh(\frt^{\downarrow})]= \sum_{n \geq 0} \mathbb{P}\Big(\bfh(\frt^{\downarrow})\geq n\Big) \underset{\text{Proposition~\ref{prop:ascendingdescendingtreesHDHPareGW}}}{=} \sum_{n \geq 0} 1-\underset{:=u_n}{\underbrace{\bbP(\mathbf{Z}_n =0)}}
\end{align}
The stake is to derive an exact formula for the $u_{n}$. 
It is classic and easy to see that $u_n$ satisfies~$u_0=~0$, and for every $n\geq 0$, the recurrence relation
\begin{align}
\label{recurrencerelation}
u_{n+1} = \phi(u_{n}),
\end{align} 
where $\phi$ is the generative function $\phi: [0,\frac{1}{1-\alpha}) \to \mathbb{R}_{+}$ of the distribution $\bdmu_{1-\alpha}$, which is equal to
\begin{align}
\label{generativefunction}
\phi(z) := \sum_{n \geq 0} \bdmu_{1-\alpha}(n) \cdot z^n = \frac{\alpha}{1-(1-\alpha)\cdot z}.
\end{align}
The next step is then to find the fixed points of the homography:
\begin{align}
\phi(z)=z \Leftrightarrow \alpha \cdot z^2 - z + (1-\alpha) = 0  \Leftrightarrow z=1 \ \text{or} \ z=\frac{\alpha}{1-\alpha}=m>1. \nonumber
\end{align}
Let us introduce now the new sequence $v_{n}:=\frac{u_{n}-m}{u_{n}-1}$. We check that it is geometric of ratio~$m$
and immediately deduce that for every $n \geq 1$:
\begin{align}
v_{n}= m^{n-1} \cdot v_{1} = m^{n-1} \cdot \frac{\alpha^2}{(1-\alpha)^2} = m^{n+1}. \nonumber
\end{align}
Which yields the conclusion. \qedd
\end{prove}

\subsection{A peeling exploration and a spatial Markov property} \label{sec:peelingexplorationSCTHPsmp}

In this section, we describe a particular way to \textit{peel} a SCTHP, that exhibits a nice spatial Markov property. By peeling, we mean that the map is "revealed" step by step through some process: at each time, we only know and observe a portion of the map---a \textit{submap} of it. In the next definition---in which we introduce the fundamental pattern of our \textit{peeling process}, it will become clear for the reader what is a submap.

\begin{figure}[!h]
\centering 
\subfigure[An upward revelation after the peeling of the edge~$(0,0)\leftrightarrow (0,1)$.]{
\centering
\includegraphics[scale=0.345]{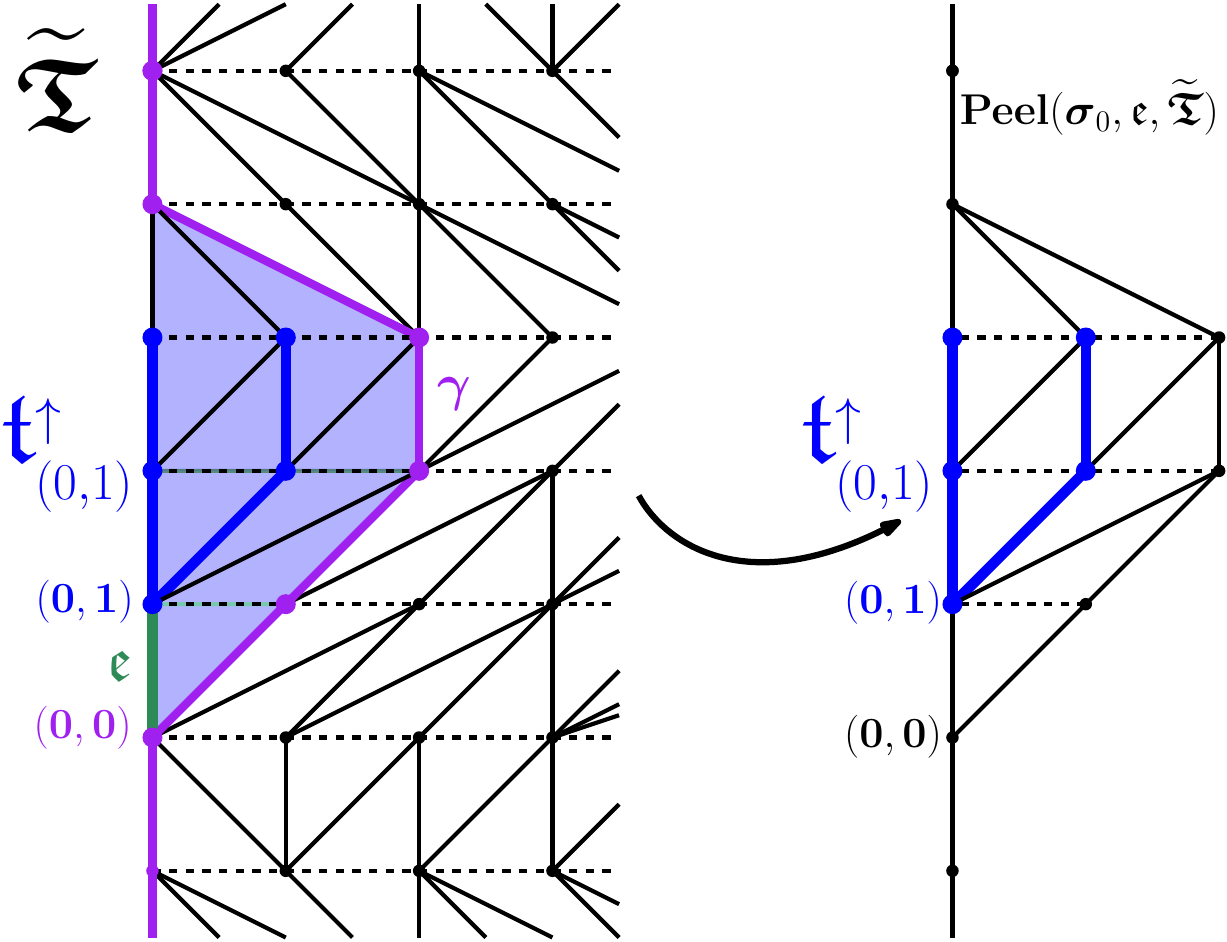}
}
\quad
\subfigure[A downward revelation after the peeling of the edge~$(0,0)\leftrightarrow (0,1)$.]{
\includegraphics[scale=0.345]{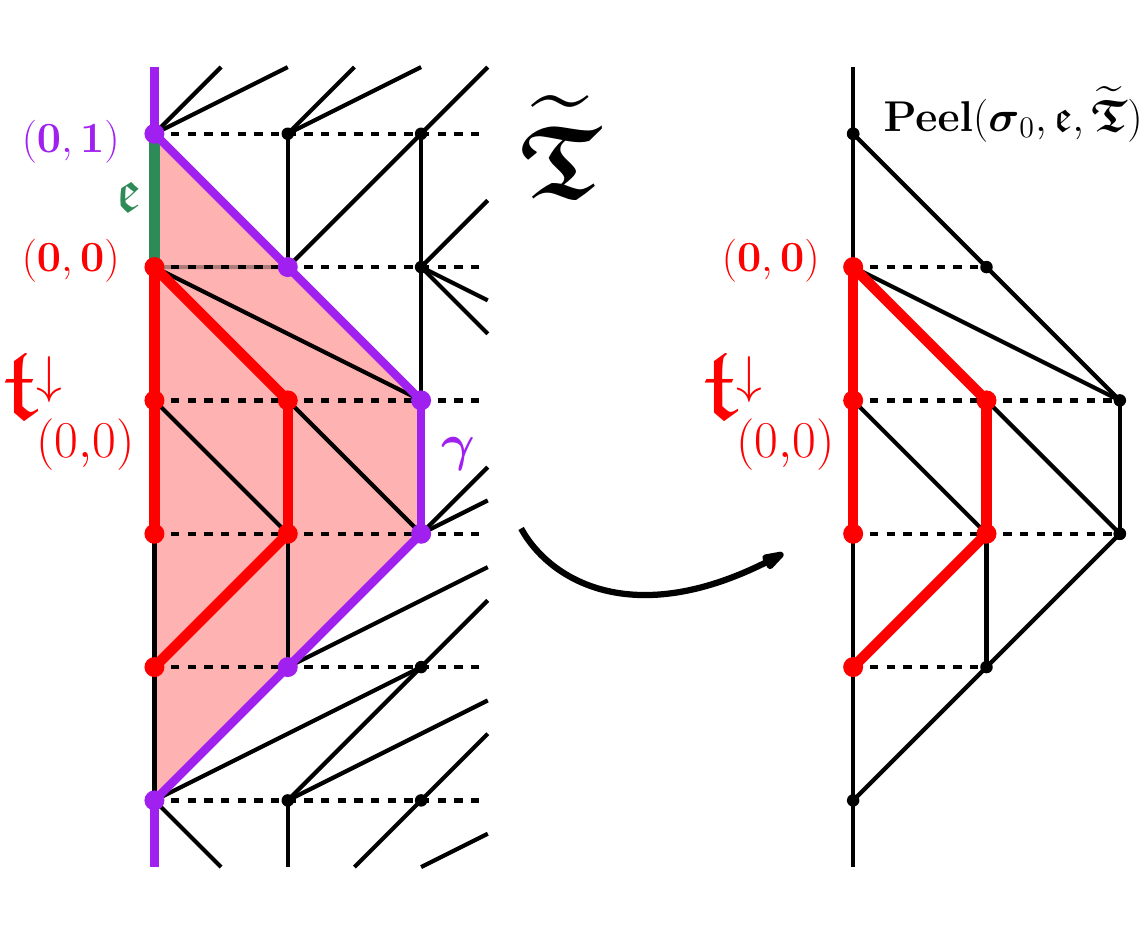}
}
\caption{Illustration of the peeling of the edge $\fre=(0,0)\leftrightarrow (0,1)$. Filled in light blue or red are the faces glued to form~$\text{Peel}(\bdsigma_0,\fre,\tildebdfrT)$. They are bordered on the right by the leftmost path $\bdgamma$ defined in Remark~\ref{remark:otherconstructionpeelingstep}.}
\label{figure:exampleupwardanddownwardrevelationfortheuniversalsubmap}
\end{figure}

\begin{definition}
\label{def:peelinguniversalsubmap}
\NoEndMark Let $\tildebdfrT$ be a SCTHP and $\fre$ be an edge of its \emph{left boundary} $\bdsigma_0:=\left\{ 0 \right\} \times \bbZ$. Note that~$\fre=~(x,y)$ is necessarily such that $x=(0,n)$ and $y=(0,n+1)$ for some $n \in \bbZ$. Then:
\begin{itemize}
\item \textbf{(upward revelation)} if the triangle bordering $\fre$ on its right is \emph{bottom-oriented}, we reveal all the vertices of the ascending tree $\frt_{y}^{\uparrow}$, and all the edges emerging from them;
\item \textbf{(downward revelation)} otherwise, the triangle is \emph{top-oriented}. Here, we reveal all the vertices of the descending tree $\frt_{x}^{\downarrow}$, and again all the edges emerging from them.
\end{itemize}
In both cases, our procedure is actually equivalent to reveal a certain number (maybe infinite) of triangles in $\tildebdfrT$. Glued together with the unique face of infinite degree of~$\tildebdfrT$ (that on the left of the infinite path~$\left\{ 0 \right\} \times \bbZ$), we obtain a planar map that 
we denote by~$\text{\emph{Peel}}(\bdsigma_0,\fre,\tildebdfrT)$. It is the \emph{submap of~$\tildebdfrT$ obtained after the peeling of the edge~$\fre$}. See Figure~\ref{figure:exampleupwardanddownwardrevelationfortheuniversalsubmap} for an illustration. The map which instead results from the gluing of the non-revealed faces of the SCTHP is called the \emph{complement submap} of~$\text{\emph{Peel}}(\bdsigma_0,\fre,\tildebdfrT)$ and is written~$\tildebdfrT \setminus \text{\emph{Peel}}(\bdsigma_0,\fre,\tildebdfrT)$. 
\end{definition}

We \emph{label} vertices both in~$\text{Peel}(\bdsigma_0,\fre,\tildebdfrT)$ and in the complement submap. In the former, we simply keep labels of the underlying SCTHP. In the latter, heights are unchanged but we translate in a trivial way the first coordinates so that we get a labelling in bijection with $\bbN \times \bbZ$.
Note that the submap $\text{Peel}(\sigma_0,\fre,\tildebdfrT)$ has two faces of infinite degree, which surround triangles. We define the \textit{right boundary} as the path bordering the right infinite face. 

\begin{remark}
\label{remark:otherconstructionpeelingstep}
\NoEndMark An alternative way to construct $\text{\emph{Peel}}(\bdsigma_0,\fre,\tildebdfrT)$ is as follows. If the triangle bordering~$\fre$ is bottom-oriented, we set $\bdgamma$ as the leftmost (vertical) path crossing the edge $(x,z)$, where $z$ is the right neighbour of $y$. Otherwise, we set $\bdgamma$ as the leftmost (vertical) path crossing the edge $(z,y)$, where $z$ is this time the right neighbour of $x$. In both cases, the submap $\text{\emph{Peel}}(\bdsigma_0,\fre,\tildebdfrT)$ is obtained by gluing together all the faces of $\tildebdfrT$ located on the left of~$\bdgamma$. See Figure~\ref{figure:exampleupwardanddownwardrevelationfortheuniversalsubmap} for an illustration.
\end{remark}

As announced earlier, the above peeling procedure highlights a spatial Markov property of SCTHP, summed up in two points: 
\begin{enumerate}
\item the submap $\text{Peel}(\bdsigma_0,\fre,\tildebdfrT)$ and its complement $\tildebdfrT \setminus \text{Peel}(\bdsigma_0,\fre,\tildebdfrT)$ are independent;
\item the complement submap $\tildebdfrT \setminus \text{Peel}(\bdsigma_0,\fre,\tildebdfrT)$ is distributed as a SCTHP.
\end{enumerate}
Both are simple consequence of the alternative construction of SCTHP, via the collection of ascending and descending trees emerging from vertices with same given height. The peeling operation can actually be regarded as follows. It consists of removing from the latter collection an ascending or a descending tree associated to some vertex located on the left boundary of $\tildebdfrT$ (either~$x$ or $y$, see Definition~\ref{def:peelinguniversalsubmap}). The new map that we get from the amputated collection is merely the complement submap $\tildebdfrT \setminus \text{Peel}(\bdsigma_0,\fre,\tildebdfrT)$. Proposition~\ref{prop:ascendingdescendingtreesHDHPareGW} ensures that this map is distributed as a SCTHP and is independent of $\text{Peel}(\bdsigma_0,\fre,\tildebdfrT)$.

\begin{figure}[!h]
\centering 
\subfigure[First step: an upward revelation.]{
\centering
\includegraphics[scale=0.52]{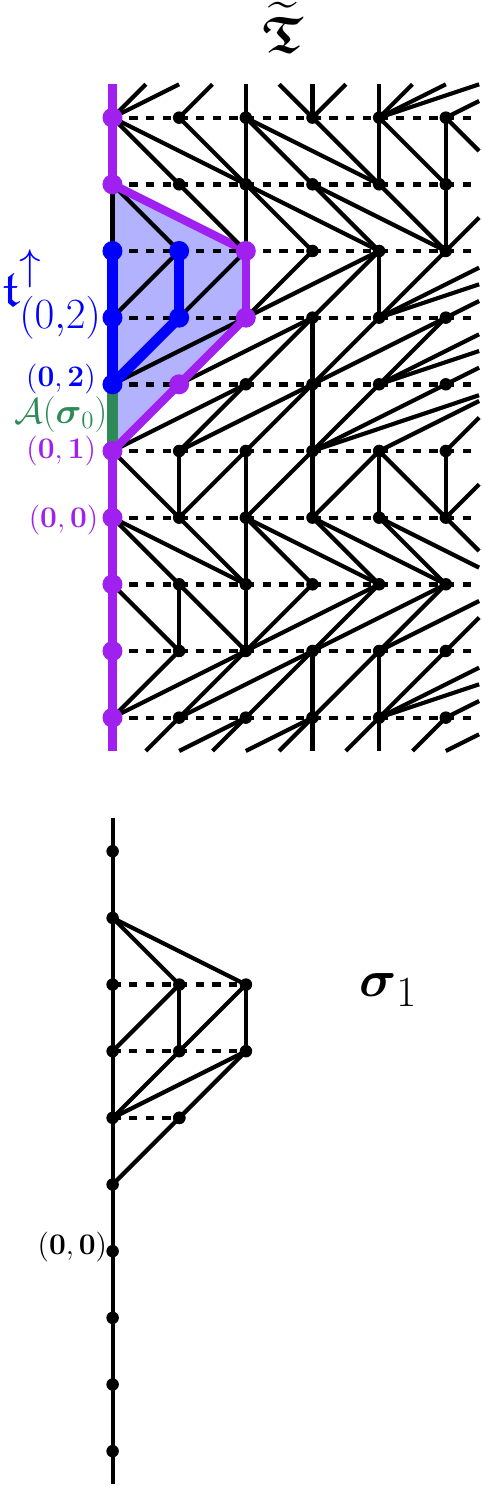}
}
\quad
\subfigure[2nd step: an upward revelation.]{
\includegraphics[scale=0.52]{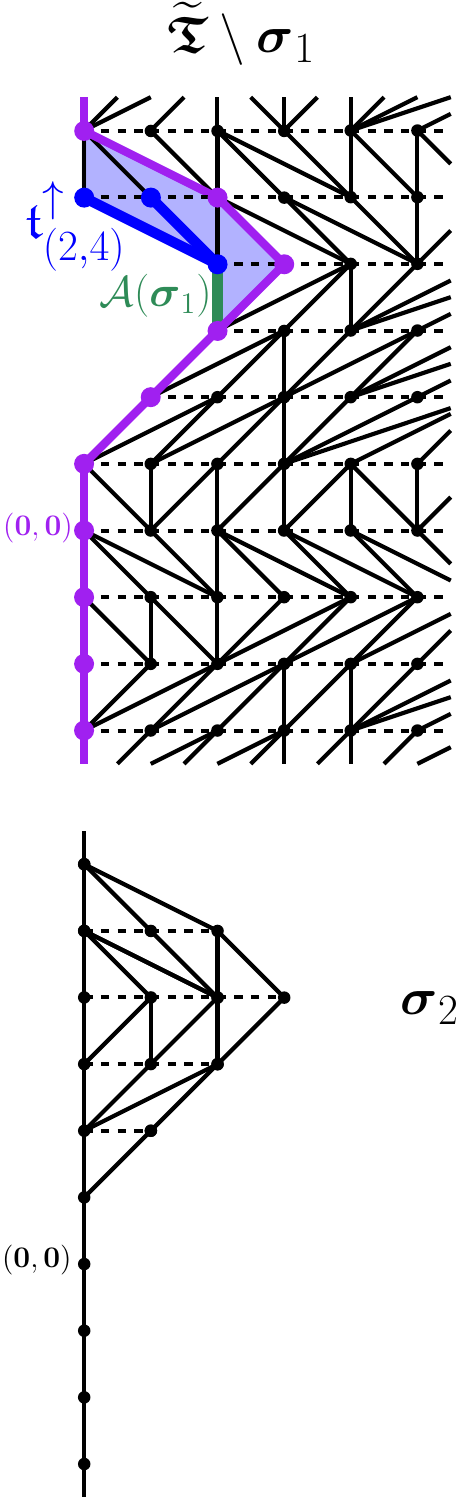}
}
\quad
\subfigure[Third step: a downward revelation.]{
\includegraphics[scale=0.52]{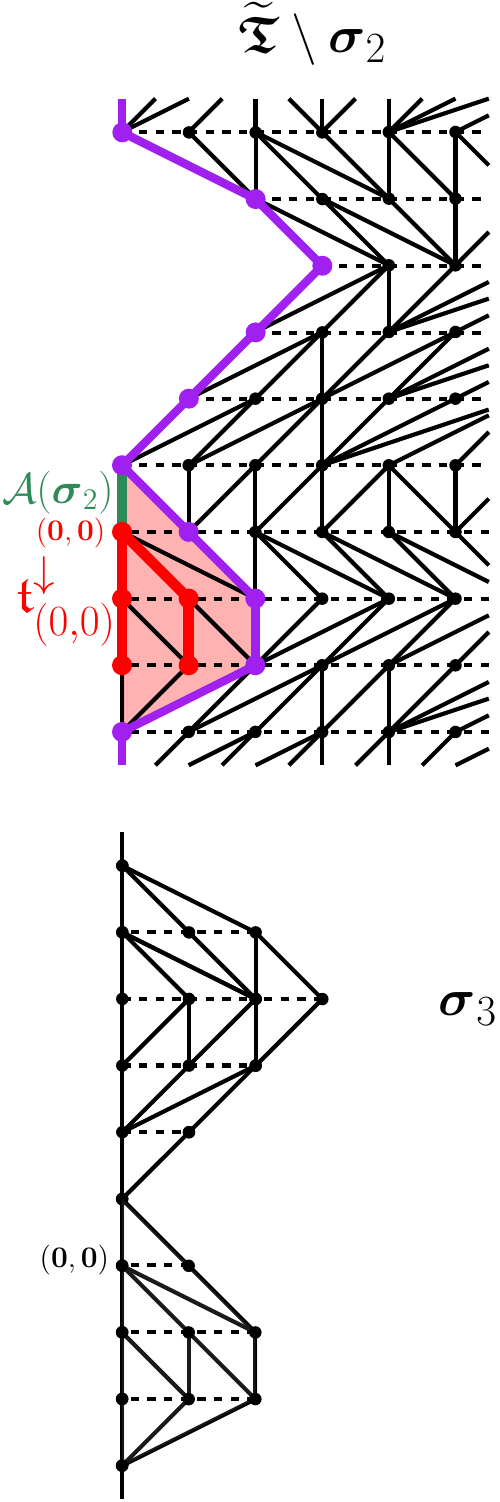}
}
\caption{An example of the first three steps of a peeling exploration for some algorithm $\calA$.}
\label{figure:illustrationdeterministiclamepeelingexploration}
\end{figure}

The second point of the spatial Markov property makes easy to iterate the peeling pattern of Definition~\ref{def:peelinguniversalsubmap}. We just need to choose an edge on the right boundary of $\text{Peel}(\bdsigma_0,\fre,\tildebdfrT)$, then peel the latter in the complement submap. It reveals new faces that we glue to our current submap, along the right boundary. With a deterministic function---the \textit{peeling algorithm}---that associates to any submap an edge on its right boundary, we can indefinitely pursue the process---or, say, the \textit{peeling exploration}. See Figure~\ref{figure:illustrationdeterministiclamepeelingexploration} for an illustration. The edge to peel can also be picked randomly, with a randomness source independent of the current submap. This is actually what we will do in the next section to study percolation on SCTHP.

\section{Percolation on SCTHP}
\label{sec:directedpercoontoymodelSCTHP}

Before dealing with SCT, a preliminary step is to grasp oriented percolation in the context of SCTHP, \textbf{with horizontal edges playing no role} and the others directed from bottom to the top. We perform on the latter a Bernoulli bond percolation process of parameter $p \in [0,1]$, as in SCT. We write $\widetilde{\bbP}_p$ for the \emph{annealed} distribution that averages on the random maps and the percolation process, and $\widetilde{\bbE}_p$ for the expectation computed under $\widetilde{\bbP}_p$.
We first affirm that a phase transition occurs in SCTHP and is similar to that in SCT---see Theorem~\ref{theorem:phasetransitionSCQ}:

\begin{theorem}\label{theorem:phasetransitiontoymodel}
\NoEndMark Let $\tildebdcalC$ be the percolation cluster of the origin in $\tildebdfrT$ and $\widetilde{\Theta}(p)$ the \emph{annealed} probability that $\tildebdcalC$ is infinite.
We have:
\begin{align}
\label{eq:probaclusterorigininfiniteSCTHP}
\widetilde{\Theta}(p) = \left\{
    \begin{array}{ll}
        0 & \mbox{if } p \leq p_c \\
	 >0 & \mbox{if } p > p_c.
    \end{array}
\right.
\end{align}
~where the quantity $p_c$ is defined in Equation~\eqref{theorem:exactexpressioncriticalpointdirectedpercolationforSCQ}. 
\end{theorem}

The method to prove the above theorem is based on the spatial Markov property raised in Section~\ref{sec:peelingexplorationSCTHPsmp} and is broadly inspired of the approach adopted in previous works like \cite{angel2003growth,angel2005scaling, angel2015percolations}.  To sum up, we will conceive a smart exploration of the cluster $\tildebdcalC$---roughly by following its leftmost interface ---from which we will exhibit some real random walk whose the asymptotic behaviour---does it diverge or not?---actually governs the percolation phenomenon. Successive moves of the walk will be interpreted as resulting of peeling actions. That will allow us to derive the exact distribution of the process and its drift. Details are outlined in Section~\ref{sec:explorationclusterRW}.

We then go deeper into the analysis of the phase transition and compute several related annealed exponents which are as expected---at least for two of them---the same as in the SCT model: 

\begin{theorem} 
\label{theorem:criticalexponentstheoremSCTHP}
\NoEndMark We have:
\begin{enumerate}[(i)]
\item $\widetilde{\bbP}_p(\lvert \tildebdcalC \rvert  \geq n)= \calO\big(c^n\big)$ for some $c=c(p) \in (0,1)$ when $p<p_c$; 
\item $\widetilde{\bbP}_{p_c}(\lvert \tildebdcalC \rvert  \geq n) \sim A \cdot n^{-1/2}$ for some $A>0$;
\item $\widetilde{\Theta}(p) \sim a \cdot (p-p_c)$ for some $a>0$, as $p \xrightarrow[]{} (p_c)^+$.
\end{enumerate}
\end{theorem}
The first estimate is commonly called "sharpness of the phase transition". Note that we do not prove in these pages any equivalent for SCT, although we believe that it is true in such context as well. All the three asymptotics stated in Theorem~\ref{theorem:criticalexponentstheoremSCTHP} are computed with the help of the exploration of~$\tildebdcalC$ above mentioned. We notably use the fact that the size of the cluster is roughly proportional to the length of the first non-negative excursion of the associated random walk, which makes the estimates on the tail distribution of $\lvert \tildebdcalC \rvert$ look more familiar. Indeed, the length of a non-negative excursion is for instance known to be exponentially-tailed when the drift of the random walk is negative (a simple large deviations result), and polynomially-tailed when the drift is null~\cite{vatutin2009local}. The third item is the easiest to get, see the end of Section~\ref{sec:explorationclusterRW}. In Section~\ref{sec:volumepercocluster}, we prove the first and second items.

Our third concern is to investigate the geometry of a large critical cluster and demonstrate that the Brownian continuum random tree emerges---like in SCT---as its scaling limit:

\begin{theorem}\label{theorem:scalinglimittoymodel}
\NoEndMark At $p=p_c$, we have:
\begin{align}
\label{theorem:scalinglimittoymodelconvergenceequality}
\big( n^{-1/2} \cdot \tildebdcalC \ \big\vert \ \vert \tildebdcalC \rvert \geq n \big) \xrightarrow[n \to +\infty]{d_{GH}} \kappa' \cdot \bdcalT_{\geq 1},
\end{align}
where $\kappa'=\kappa'(\alpha)$ is a positive number only depending on $\alpha$ and $\bdcalT_{\geq 1}$ is the CRT of mass greater than~1 \cite{aldous1991continuum, aldous1993continuum}. The convergence \eqref{theorem:scalinglimittoymodelconvergenceequality} holds in distribution for the Gromov--Hausdorff distance.
\end{theorem}
Again, the proof rests on a clever use of the cluster's exploration. While the convergence in~\eqref{theorem:scalinglimittoymodelconvergenceequality} involves metric spaces, we will rewrite it in terms of random processes instead. The result will decisively stem from the Donsker's theorem. Details are tackled in Section~\ref{sec:scalinglimitpercocluster}. Definitions of the Gromov--Hausdorff topology and of the CRT are also recalled in the latter.

\subsection{Exploration of the percolation cluster via a peeling process}\label{sec:explorationclusterRW}
A simple way to find an infinite \textbf{\emph{directed}} path---provided that it exists--- in the cluster $\tildebdcalC$ is by tracing its contours from left to right. We start by following the leftmost \emph{directed} path emerging from the origin vertex, as far as possible. At some point, we may reach a dead end---a vertex, say~$x$, with no offspring in the cluster: all the opportunities to go further have been exhausted. A clear consequence is that at our current height in the SCTHP, any infinite directed path in the cluster necessarily crosses an edge located strictly on our right. The leftmost such edge is that connecting the first right neighbour of $x$ to the rightmost offspring of $x$---say $z$ and $y$. We look now for a way to reach the latter edge $(z,y)$ from $(0,0)$. 

\begin{figure}[!h]
    \centering
\includegraphics[width=140pt]{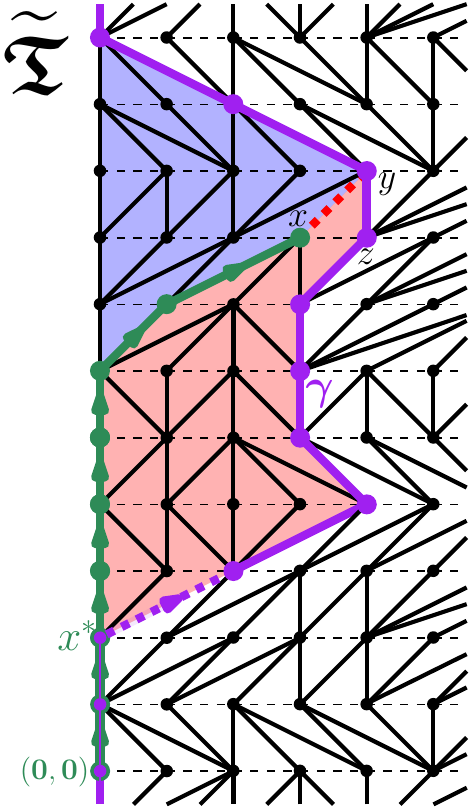}
    \caption{In green is the leftmost directed path in $\tildebdcalC$, going from the origin to $x$. Filled in light blue is a part of $\tildebdfrT$ which is inaccessible from $(0,0)$ in a directed fashion, given the first steps of the exploration. The area filled in light red may contain vertices of $\tildebdcalC$, but they are not visited by the process. In purple is the path $\gamma$, which is the leftmost crossing the edge~$(z,y)$. The walk restarts at $x^{*}$, the highest already visited vertex which is met by $\gamma$, and continues on the right of the path.}
\label{fig:exampleclusterexploration}
\end{figure}

For this purpose, it is natural to consider $\gamma$ the leftmost directed path in the map passing through the edge $(z,y)$. The friendly situation is when $\gamma$ meets somewhere downward a vertex $x^*$ (take the highest one), on the segment explored by the walk from $(0,0)$ to~$x$. In such case, we get a possible path from the origin to~$z$. So we restart the walk at $x^*$ and continue the exploration on the right of $\gamma$, according to the same principle (following the leftmost directed path up to a dead end, etc). 
We insist again on the fact that an infinite path in $\tildebdcalC$ cannot be found anywhere else but on the right of $\gamma$. By dropping to $x^*$, some vertices of the cluster are set aside by the walk, but we do know that they are in finite number, enclosed between the explored segment $[x^*,x]$ and~$\gamma$. See Figure~\ref{fig:exampleclusterexploration} for an illustration.

An unfriendly situation may also happen when $\gamma$ do not intersect vertices visited earlier. The finding is in fact more definitive: the cluster $\tildebdcalC$ do not merely contain any infinite directed path since then, there is no way to get through $\gamma$ from the origin and connect itself to an infinite directed component.

\begin{paragraph}{The link with peeling}A crucial point in our analysis is to read the above exploration of $\tildebdcalC$ through the prism of peeling.
All along the process, is in fact extended a \textit{submap} of $\tildebdfrT$ containing for sure a finite portion of the cluster, while we may still find an infinite component in its complement.
Indeed, when we are at some vertex~$x$, we successively check from left to right whether the edges emerging from $x$ are in $\tildebdcalC$ or not. This operation is performed until we perhaps find one open. Each time that we find a closed edge---say $(x,y)$---we derive a set of vertices becoming inaccessible in a directed fashion from $x$ and, a fortiori, from the origin.  Those are exactly the vertices in the ascending tree of $y$. So, in the case of such event, the submap in question can be enlarged like after an upward revelation. When we finally end up not having any possibility to go further from~$x$, the walker dropping to the vertex~$x^*$ defined in the foregoing, we enlarge this time our current submap with the finite area bounded on the left by~the explored~segment~$[x^*,x]$ and by~$\gamma$ on the right. It is not hard to see that such operation actually corresponds to a downward revelation. See indeed Remark~\ref{remark:otherconstructionpeelingstep}.
\end{paragraph}

\begin{paragraph}{The peeling exploration of $\tildebdcalC$}
To be more formal, from a SCTHP $\tildebdfrT$, we define an infinite sequence $\displaystyle{\big(\bdcalV_n;\vfre_n; \bdSigma_n \big)_{n \geq 0}}$ of vertices, directed edges and submaps, obtained by applying a \textit{randomized} peeling algorithm. Indeed, contrary to our example in Figure~\ref{figure:illustrationdeterministiclamepeelingexploration}, the algorithm $\calA$ will also be a function of random external data---independent of the SCTHP---provided by the Bernoulli percolation process.
We start with $$\bdcalV_0=(0,0), \ \vfre_0=(0,0) \rightarrow (0,1) \text{ and } \bdSigma_0=\left\{0\right\} \times \bbZ \text{ (the left boundary of $\tildebdfrT$)}.$$ 
 
The update of the sequence depends on the state of the edge $\vfre_0$ in $\tildebdcalC$:
\begin{itemize}
\item if $\vfre_0$ is \textit{open} in $\tildebdcalC$, we simply update the vertex sequence as $$\bdcalV_{1}:=\bdcalV_{0}+\vfre_0,$$
i.e.~$\bdcalV_{1}$ is the next vertex on the right boundary of $\bdSigma_0$---here $(0,1)$.
The edge $\vfre_{1}$ is the edge emerging from $\bdcalV_{1}$ on the same boundary---namely $(0,1) \to (0,2)$. The submap does not change.
\item otherwise, if $\vfre_0$ is \textit{closed} in $\tildebdcalC$, we peel the (not directed) edge $\fre_0$. The new submap is then $$\bdSigma_{1}=\text{Peel}(\bdSigma_{0},\fre_n,\tildebdfrT).$$ 
In case of an \textit{upward revelation}, the current vertex does not change
but the new directed edges $\vfre_{1}$ is now the second emerging from $\bdcalV_0$ (and is on the right boundary of $\bdSigma_{1}$). It is~$(0,0) \to~(1,1)$. See Figure~\ref{figure:downwardupwardrevelationsvertexupdateduringclusterexploration}(a) for an illustration.

In case of an \textit{downward revelation}, triangles added to the current submap are directly got from the descending tree of $\bdcalV_0$ in $\tildebdfrT$ (see Definition~\ref{def:peelinguniversalsubmap}). If the latter has height $h\geq 0$, we fix~$\bdcalV_{1}$ as the rightmost vertex of $\bdSigma_{1}$ (or $\bdSigma_{0}$) with height $\bfh(\bdcalV_0)-1-h$. Here, it is of course the vertex $(0,-1-h)$. The edge $\vfre_{1}$ is set as that emerging from $\bdcalV_{1}$ to its offspring which is on the right boundary of $\bdSigma_1$. See Figure~\ref{figure:downwardupwardrevelationsvertexupdateduringclusterexploration}(b) for an illustration.
\end{itemize}

The update is well-defined since the height of the descending tree $\frt_{(0,0)}^{\downarrow}$ is finite almost surely, given it is distributed as a subcritical Galton--Watson tree according to Proposition~\ref{prop:ascendingdescendingtreesHDHPareGW}. The algorithm is then iterated. This is possible because, by construction, the vertex $\bdcalV_n$ and the edge~$\vfre_n$ are always on the right boundary of $\bdSigma_n$. See Figure~\ref{figure:downwardupwardrevelationsvertexupdateduringclusterexploration}(c) for an illustration of a complete exploration of the cluster $\tildebdcalC$. 
\end{paragraph}

\begin{paragraph}{The height process}
Derived from the peeling exploration, the height process $\bdcalH:=\big(h(\bdcalV_n)\big)_{n \geq 0}$ plays a key role in our work. For instance, we show below that it condenses many information about the size of the cluster $\tildebdcalC$.
Given the spatial Markov property, combined with the fact that updates during the peeling exploration all obey the same rules, it is straightforward that the height process is a random walk with i.i.d.~increments. Their common distribution can be explicitly computed:
\begin{samepage}
\begin{proposition}
\label{prop:distributionincrementsheightprocess}
\NoEndMark Recall that $\bfh(\cdot)$ is the function returning the height of any rooted tree. Then we have:
\begin{align}
\label{eq:distributionincrementsheightprocess}
\widetilde{\bbP}_p(h(\bdcalV_1) = h) = \left\{
    \begin{array}{ll}
        p & \mbox{for } h=1 \\
        \alpha \cdot (1-p) & \mbox{for } h=0 \\
	 (1-\alpha) \cdot (1-p) \cdot \bbP(\bfh(\frt ^{\downarrow}) = -h-1) & \mbox{for } h<0
    \end{array}
\right.
\end{align}
where $\frt ^{\downarrow} \overset{d}{\sim} \mathbf{GW}_{1-\alpha}$. 
\end{proposition}
\end{samepage}
The proof of  \eqref{eq:distributionincrementsheightprocess} is totally transparent if we carefully read the instructions guiding updates during the peeling exploration, and if we remember that the graph structure of the SCTHP is independent of the percolation process. From \eqref{eq:distributionincrementsheightprocess} and Proposition~\ref{prop:heightsubcriticalGWtrees}, we easily deduce that the distribution of increments has a finite first moment and:
\begin{align}
\label{eq:firstmomentincrementheightprocess}
\widetilde{\bbE}_p[h(\bdcalV_1)]=p-(1-p)\cdot (1-\alpha)\cdot (1+\bbE[\bfh(\frt ^{\downarrow})])=p-(1-p)\cdot (1-\alpha)\cdot \sum_{n \geq 0} \frac{m-1}{m^{n+1}-1}.
\end{align}
Recall that $m$ is the mean of a geometric law of parameter $\alpha$, so is equal to $\alpha \cdot (1-\alpha)^{-1} >1$. 
\end{paragraph}
\begin{paragraph}{Percolation is governed by the drift of $\bdcalH$: proof of Theorem~\ref{theorem:phasetransitiontoymodel}.}
Let $T$ be the first hitting time of \textit{negative} integers $\displaystyle{\{-1,-2,\dotsc\}}$ by the height process $\bdcalH$. We have:
\begin{align}
\label{eq:sizeclusterlengthexcursion}
\lvert \tildebdcalC \rvert = +\infty \quad \text{ if and only if } \quad T=+\infty.
\end{align}
Indeed, vertices visited up to the time right before $T$ are all in the cluster. At time $T<+\infty$, we drop to negative heights after a downward revelation. It means that the right boundary of the new submap $\bdSigma_T$---or the path $\gamma$ defined at the beginning of the current section---did not meet any vertex visited earlier. And as it has been argued, it implies that $\tildebdcalC$ is finite. See Figure~\ref{figure:downwardupwardrevelationsvertexupdateduringclusterexploration}(c) for an illustration. The converse statement is trivial. 
With \eqref{eq:sizeclusterlengthexcursion}, the proof of Theorem~\ref{theorem:phasetransitiontoymodel} is reduced to the well-known fact that:
$$\widetilde{\bbP}_p(T=+\infty)>0 \text{ if and only if } \widetilde{\bbE}_p[h(\bdcalV_1)]>0.$$
Given \eqref{eq:sizeclusterlengthexcursion}, this yields the conclusion of Theorem~\ref{theorem:phasetransitiontoymodel}. \qedd \end{paragraph}

\begin{paragraph}{The off-critical percolation probability}
Let $T'$ be the first hitting time of \textit{non-positive} integers $\displaystyle{\{0,-1,-2,\dotsc\}}$ by the height process $\bdcalH$. In the case $p>p_c$, the drift of the latter is positive and a standard result on skip-free ascending walk---which are random walks with integer valued increments and~$+1$ as the only possible positive one, like is~$\bdcalH$---ensures that:
\begin{align}
\widetilde{\bbP}_p(T'=+\infty)=\widetilde{\bbE}_p[h(\bdcalV_1)]. \nonumber
\end{align}
See for instance \cite[Theorem~3]{addario2008ballot}. When $T=+\infty$, the height process may touch $0$ a finite number of times $\mathbf{n}$ before never come back again. Thanks to the strong Markov property, the distribution of~$\mathbf{n}$ is geometric. Its parameter is the probability for $\bdcalH$ to carry out a positive excursion and finish it at $0$. As a consequence:
 \begin{align}
\label{eq:offcriticalprobadecompo}
\widetilde{\Theta}(p)= \widetilde{\bbE}_p[\mathbf{n}] \cdot \widetilde{\bbE}_p[h(\bdcalV_1)]=\big(1-\widetilde{\bbP}_p(T'<+\infty; \ h(\bdcalV_{T'})=0)\big)^{-1} \cdot \widetilde{\bbE}_p[h(\bdcalV_1)].
\end{align}
Since $T'<+\infty$ almost surely at $p=p_c$, we clearly have \begin{align}
\widetilde{\bbP}_p(T'<+\infty; \ h(\bdcalV_{T'})=0) \xrightarrow[p \to p_c]{} \widetilde{\bbP}_{p_c}(h(\bdcalV_{T'})=0)<1. \nonumber
\end{align}
According to \eqref{eq:firstmomentincrementheightprocess}, we also have 
$\widetilde{\bbE}_p[h(\bdcalV_1)] \underset{p \to p_c}{\sim} c \cdot (p-p_c)$ for some $c>0$. Together with \eqref{eq:offcriticalprobadecompo}, the last two estimates imply Theorem~\ref{theorem:criticalexponentstheoremSCTHP}(iii). \qedd

\end{paragraph}

\subsection{The volume of the percolation cluster}\label{sec:volumepercocluster}

\begin{figure}[!h]
\centering
\subfigure[Updates after an upward revelation occuring at time~$t=~6$. The area filled in light blue is inaccessible from the origin vertex: we remove it and pursue the exploration on the complement submap.]{
\includegraphics[height=185pt]{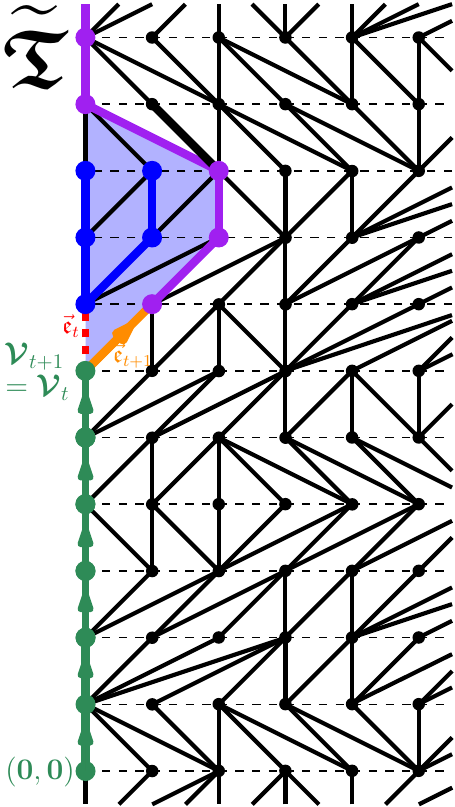}
}
\quad 
\subfigure[Updates after a downward revelation occuring at time~$n=~12$. The area filled in light red may contain vertices of the cluster (among those stroke in red), but not visited by the process.]{
\includegraphics[height=185pt]{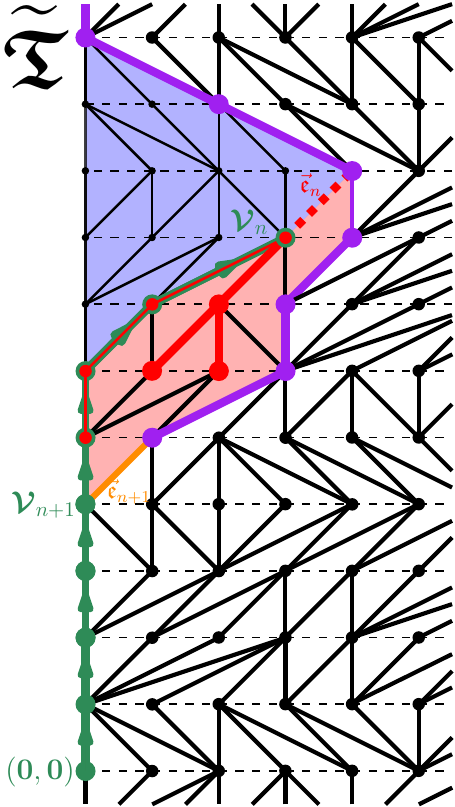}
}
\quad
\subfigure[A complete exploration of $\tildebdcalC$. Filled in light blue/red are triangles removed after an upward/downward revelation. The cluster may contain vertices in red areas above the null height. Dashed red lines are for the edges which turned out to be closed after checking.]{
\includegraphics[height=185pt]{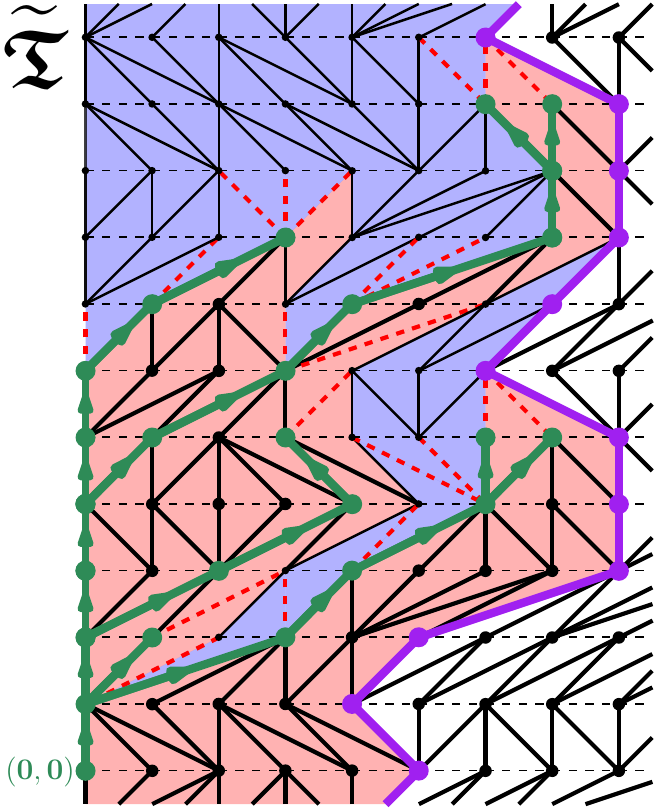}
}
\quad
\subfigure[Graph of the height process along the exploration of~$\tildebdcalC$.]{
\includegraphics[height=155pt]{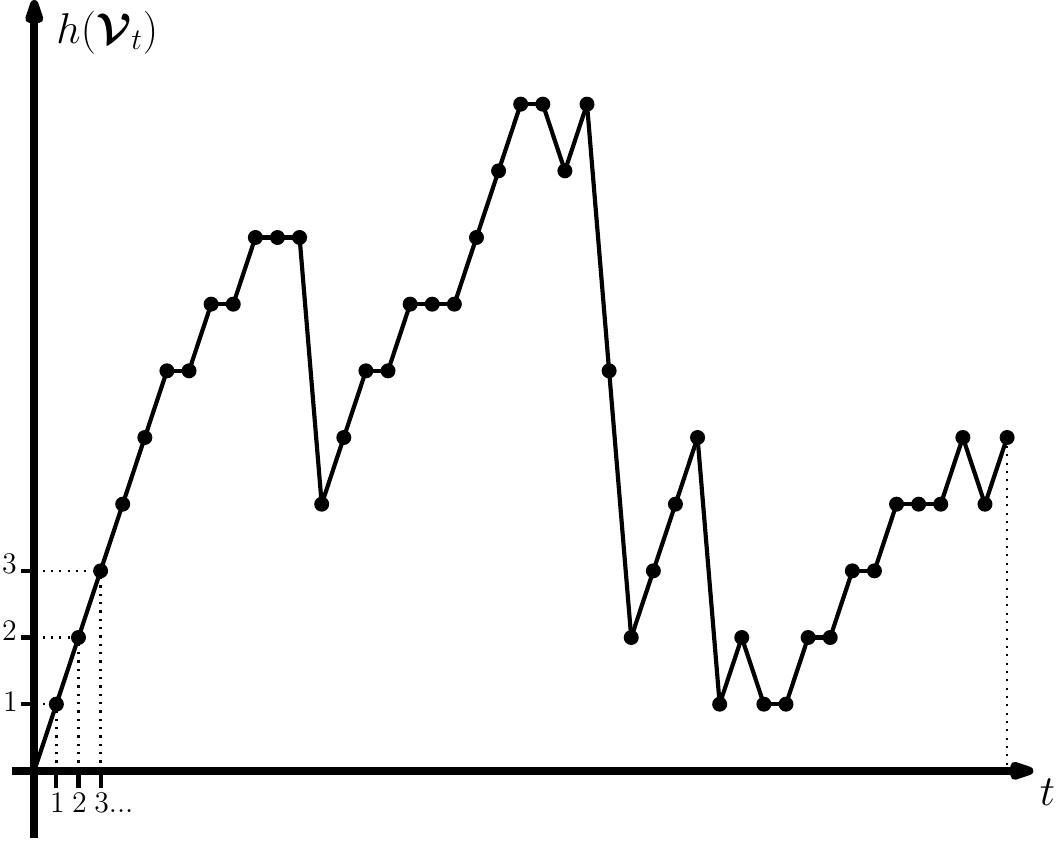}
}
\quad \quad  \quad \quad 
\subfigure[The trace of the walk is tree-shaped.]{
\includegraphics[height=155pt]{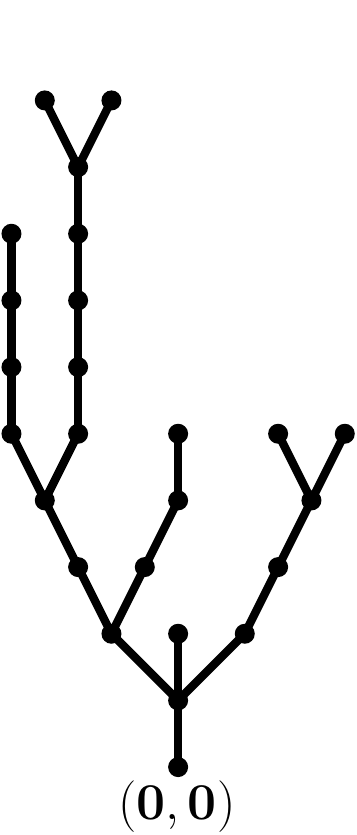}
}
\caption{Illustration of the exploration of~$\tildebdcalC$.}
\label{figure:downwardupwardrevelationsvertexupdateduringclusterexploration}
\end{figure}

Recall that we have defined earlier $T$ as the first hitting time of negative numbers by the height process $\bdcalH$. At the critical threshold $p=p_c$, the cluster $\tildebdcalC$ is finite almost surely and is fully contained in the submap $\bdSigma_T$. Actually, it is merely the union of vertices explored up to the time right before $T$---with a new vertex visited at each $+1$ jump in the associated height process--- plus some others left behind, at times of downward revelation, or negative jump for $\bdcalH$. See Figure~\ref{figure:downwardupwardrevelationsvertexupdateduringclusterexploration}(c). The latter vertices are all enclosed in areas whose volume is distributed as the total progeny of a subcritical Galton--Watson tree with offspring law~$\bdmu_{1-\alpha}$. See Figure~\ref{figure:exampleupwardanddownwardrevelationfortheuniversalsubmap}(b). 
It is well known that the latter quantity is exponentially-tailed, by considering for instance its generative function. All these facts allow us to write the following equality, holding in distribution:
\begin{align}
\label{eq:distributionequalityvolumecluster}
\lvert \tildebdcalC \rvert \overset{(d)}{=} 1+\sum_{t=1}^{T} \bdtheta_t.
\end{align}
where $(\bdtheta_t)_{t \geq 1}$ is a sequence of i.i.d.~random variables, satisfying:
$$
\bdtheta_t = \left\{
    \begin{array}{ll}
	1 & \mbox{if } h(\bdcalV_{t})-h(\bdcalV_{t-1}) = 1; \\
	0 & \mbox{if } h(\bdcalV_{t})-h(\bdcalV_{t-1}) = 0; \\
       \geq 0 & \mbox{otherwise}.
    \end{array}
\right.
$$
The random variable $\bdtheta_t$ is furthermore, in the third case (corresponding to a negative jump in the height process), stochastically dominated by the volume of a Galton--Watson tree with offspring law~$\bdmu_{1-\alpha}$, conditioned on having height $\lvert h(\bdcalV_{t})-h(\bdcalV_{t-1}) \rvert - 1$. 
In \eqref{eq:distributionequalityvolumecluster}, pay attention to the fact that $\bdtheta_{T}$ plays a special role since vertices in~$\tildebdcalC$ must have a nonnegative height. However, it does not change anything for our further analysis since the stochastic domination holds even more for~$\bdtheta_T$.

Set now $\displaystyle{\kappa:=\bbE[\bdtheta_1]<+\infty}$. Equation~\ref{eq:distributionequalityvolumecluster} suggests that the cluster $\tildebdcalC$ roughly contains $\kappa \cdot T$ vertices, as $T \to +\infty$. This is indeed true and it turns out to be a crucial point in order to prove the two first items of Theorem~\ref{theorem:criticalexponentstheoremSCTHP}. Large deviations arguments will help us to make it rigorous.

Before entering into details, we introduce a new notation which will simplify the presentation:

\begin{definition}
\label{def:notasymptoe}
\NoEndMark We write $x_n = oe_{\delta}(n)$ for some sequence~$(x_n)$ and $\delta>0$ if there exists constants~$c,C>0$ such that $x_n \leq C e^{-c \cdot n^{\delta}}$ for any $n$.
\end{definition}

\begin{remark}
\label{remark:notationepsilon}
\NoEndMark If $x_n = oe_{\delta}(n)$ and if $y_n := \sum_{k \geq n} x_k$, then $y_n = oe_{\delta}(n)$.
\end{remark}

\subsubsection{The sharpness of the phase transition}\label{sec:sharpnessphasetransitionSCTHP}
\begin{paragraph}{}
Fix $p<p_c$. We aim to prove the first item of Theorem~\ref{theorem:criticalexponentstheoremSCTHP}.  Let $n \geq 1$ and~$n' \leq n/2\kappa$. We have:
\begin{align}
\widetilde{\bbP}_p\Big( \lvert \tildebdcalC \rvert \geq n; \ T=n' \Big) \underset{\eqref{eq:distributionequalityvolumecluster}}{\leq} \widetilde{\bbP}_p\Big( 1+\sum_{t=1}^{n'} \bdtheta_t \geq n \Big) \leq \widetilde{\bbP}_p\Big( 1+\sum_{t=1}^{\lfloor n/2\kappa \rfloor} \bdtheta_t \geq n \Big). \nonumber
\end{align}
The last inequality comes from the fact that as $n'$ increases, so does $\widetilde{\bbP}_p\Big( 1+\sum_{t=1}^{n'} \bdtheta_t \geq n \Big)$. 
By a standard large deviations argument, we know that:
\begin{align}
\widetilde{\bbP}_p\Big( 1+\sum_{t=1}^{\lfloor n/2\kappa \rfloor} \bdtheta_t \geq n \Big) = oe_{1}(n). \nonumber
\end{align}
Given both inequalities above, by summing over $n' \leq n/2\kappa$, we get:
\begin{align}
\widetilde{\bbP}_p\Big( \lvert \tildebdcalC \rvert \geq n; \ T \leq n/2\kappa\Big)= oe_{1}(n). \nonumber
\end{align}
Also, it is clear that $\widetilde{\bbP}_p\Big(T > n/2\kappa\Big) \leq \widetilde{\bbP}_p\Big(h(\bdcalV_{\lfloor n/2\kappa \rfloor})\geq 0\Big)$, and the latter probability is~$oe_{1}(n)$ because the drift of the height process is negative when $p<p_c$ (so the event $\left\{ h(\bdcalV_{\lfloor n/2\kappa \rfloor})\geq 0 \right\}$ is a large deviations event). Finally, it holds that:
\begin{align}
\widetilde{\bbP}_p\Big( \lvert \tildebdcalC \rvert \geq n\Big) \leq \widetilde{\bbP}_p\Big( \lvert \tildebdcalC \rvert \geq n; \ T \leq n/2\kappa\Big) + \widetilde{\bbP}_p\Big(T > n/2\kappa\Big), \nonumber
\end{align}
and the two quantities on the right-hand side of the inequality are $oe_{1}(n)$, as it has just been argued. So is then $\widetilde{\bbP}_p\Big( \lvert \tildebdcalC \rvert \geq n\Big)$, which is the expected conclusion. \qedd
\end{paragraph}

\subsubsection{The volume of a critical cluster}\label{sec:volumecriticalclusterSCTHP}
\begin{paragraph}{}
We focus now on the second item of Theorem~\ref{theorem:criticalexponentstheoremSCTHP} with~$p=p_c$. It is well known that the distribution of the hitting time $T$ is polynomially-tailed as the underlying random walk has no drift (like here). The speed of decreasing is explicit, see for instance \cite[Theorem 11]{vatutin2009local} applied to an integer-valued random walks with square-integrable increments:
\begin{proposition}
\label{prop:taildistributionlengthnnexcursion}
There exists $A>0$ such that:
$$\widetilde{\bbP}_{p_c}(T=n) \underset{n \to +\infty}{\sim} A \cdot n^{-3/2} \quad \text{and } \quad \widetilde{\bbP}_{p_c}(T \geq n) \underset{n \to +\infty}{\sim} 2A \cdot n^{-1/2}.$$ 
\end{proposition}
Now we recall a classic moderate deviations result on exponentially-tailed random variables: 
\begin{lemma}
\label{lemma:largedeviationresultfinitesecondmoment}
\NoEndMark Let $(X_t)_{t \geq 1}$ be a sequence of i.i.d.~random variables with mean~$\overline{x}$ and an exponentially-tailed distribution, meaning that $\bbP(X_1 \geq n) = oe_{1}(n)$. Then, for any $C>0$:
\begin{align}
\bbP\Big(\Big\lvert \frac{1}{n} \sum_{t=1}^{n} X_t - \overline{x} \Big\rvert > C \cdot n^{-1/4}\Big) = oe_{1/2}(n). \nonumber
\end{align} 
\end{lemma}
This is a particular case of~\cite[Lemma 1.12]{le2005random}. With the help of the latter lemma, we are going to prove that: 
\begin{align}
\label{eq:symmetricdifferenceclusterlengthexc}
\widetilde{\bbP}_{p_c}\big( (\lvert \tildebdcalC \rvert \geq \kappa \cdot n ) \bigtriangleup ( T \geq n )\big)=\calO\Big( n^{-3/4}\Big)=o(n^{-1/2})=o\Big(\widetilde{\bbP}_{p_c}( T \geq n)\Big).
\end{align}
The symbol $\bigtriangleup$ denotes the symmetric difference of two sets. Together with Proposition~\ref{prop:taildistributionlengthnnexcursion}, the estimate \eqref{eq:symmetricdifferenceclusterlengthexc} is enough to get the conclusion since it implies that
$$\widetilde{\bbP}_{p_c}\Big( \lvert \tildebdcalC \rvert \geq \kappa \cdot n \Big) \underset{n \to +\infty}{\sim} \widetilde{\bbP}_{p_c}\Big(T \geq n\Big).$$

Fix~$C>0$. Proposition~\ref{prop:taildistributionlengthnnexcursion} ensures that 
\begin{align}
\label{eq:concentrationoflengthexcursionaroundn}
\widetilde{\bbP}_{p_c}\Big(\lvert T-n \rvert \leq C \cdot n^{3/4}\Big)=\calO\Big(n^{-3/4}\Big)=o(n^{-1/2})=o\Big(\widetilde{\bbP}_{p_c}( T \geq n)\Big).
\end{align}
Moreover, for any~$n \geq 1$ and~$n' \geq n +C \cdot n^{3/4}$:
\begin{align}
\widetilde{\bbP}_{p_c}\Big( \lvert \tildebdcalC \rvert < \kappa \cdot n; \ T=n' \Big) &\underset{\eqref{eq:distributionequalityvolumecluster}}{\leq} \widetilde{\bbP}_{p_c}\Bigg( \frac{1+\sum_{t=1}^{n'} \bdtheta_t}{n'} < \kappa \cdot n/n'\Bigg) \nonumber \\ &\leq \widetilde{\bbP}_{p_c}\Bigg( \frac{1+\sum_{t=1}^{n'} \bdtheta_t}{n'} < \kappa \cdot n/(n+C \cdot n^{3/4})\Bigg). \nonumber
\end{align}
There exists some $C'>0$ such that $n / (n+C \cdot n^{3/4}) \leq 1- \kappa^{-1} \cdot C' \cdot n^{-1/4}$ for any $n$. Since $n' \geq n$, and given the above inequality, we get that:
\begin{align}
\widetilde{\bbP}_{p_c}\Big( \lvert \tildebdcalC \rvert < \kappa \cdot n; \ T=n' \Big) \leq \widetilde{\bbP}_{p_c}\Bigg( \frac{1+\sum_{t=1}^{n'} \bdtheta_t}{n'} < \kappa -  C' \cdot (n')^{-1/4} \Bigg) \underset{Lemma~\ref{lemma:largedeviationresultfinitesecondmoment}}{=} oe_{1/2}(n'). \nonumber
\end{align}
By summing over $n'$, we deduce that:
\begin{align}
\nonumber
\widetilde{\bbP}_{p_c}\Big( \lvert \tildebdcalC \rvert < \kappa \cdot n; \ T \geq n+C \cdot n^{3/4}\Big) = oe_{1/2}(n).
\end{align}
Then, by using \eqref{eq:concentrationoflengthexcursionaroundn}:
\begin{align}
\label{eq:firstestimatesymmetricdifference}
\widetilde{\bbP}_{p_c}\Big( \lvert \tildebdcalC \rvert < \kappa \cdot n; \ T \geq n\Big)=\calO\Big(n^{-3/4}\Big).
\end{align}
Furthermore, for any $n' \leq \lfloor n-C \cdot n^{3/4} \rfloor :=N$, we have:
\begin{align}
\widetilde{\bbP}_{p_c}\Big( \lvert \tildebdcalC \rvert \geq \kappa \cdot n; \ T=n' \Big) \underset{\eqref{eq:distributionequalityvolumecluster}}{\leq} \widetilde{\bbP}_{p_c}\Big( 1+\sum_{t=1}^{n'} \bdtheta_t \geq \kappa \cdot n \Big) \leq \widetilde{\bbP}_{p_c}\Big( 1+\sum_{t=1}^{N} \bdtheta_t \geq \kappa \cdot n \Big). \nonumber
\end{align}
The last inequality comes from the fact that as $n'$ increases, so does $ \widetilde{\bbP}_p\Big( 1+\sum_{t=1}^{n'} \bdtheta_t \geq \kappa \cdot n \Big)$. Since there exists some $C'>0$ such that $n/N \geq 1+\kappa^{-1} \cdot C' \cdot n^{-1/4}$, we derive that:
\begin{align}
\widetilde{\bbP}_{p_c}\Big( \lvert \tildebdcalC \rvert \geq \kappa \cdot n; \ T=n' \Big) \leq \widetilde{\bbP}_{p_c}\Bigg(\frac{1+\sum_{t=1}^{N} \bdtheta_t}{N} \geq \kappa + C' \cdot n^{-1/4} \Bigg)\underset{Lemma~\ref{lemma:largedeviationresultfinitesecondmoment}}{=}oe_{1/2}(n). \nonumber
\end{align}
By summing over $n'$, we obtain that:
\begin{align}
\widetilde{\bbP}_{p_c}\Big( \lvert \tildebdcalC \rvert \geq \kappa \cdot n; \ T < n-C \cdot n^{3/4} \Big) = oe_{1/2}(n).\nonumber 
\end{align}
Then, by using again \eqref{eq:concentrationoflengthexcursionaroundn}:
\begin{align}
\label{eq:secondestimatesymmetricdifference}
\widetilde{\bbP}_{p_c}\Big( \lvert \tildebdcalC \rvert \geq \kappa \cdot n; \ T < n\Big)=\calO\Big(n^{-3/4}\Big).
\end{align}
The estimates \eqref{eq:firstestimatesymmetricdifference} and \eqref{eq:secondestimatesymmetricdifference} are together equivalent to \eqref{eq:symmetricdifferenceclusterlengthexc}. \qedd
\end{paragraph}

\subsection{Scaling limit of the critical percolation cluster}\label{sec:scalinglimitpercocluster}

We have shown that the height process condenses many information on the size of $\tildebdcalC$. It turns out that it captures the geometry as well. Consider the graph made of vertices and edges of the SCTHP, visited by the peeling exploration up to the time right before $T$. This graph is tree-shaped since is basically the result of the gluing of non-intersecting paths of the map. See~Figure~\ref{figure:downwardupwardrevelationsvertexupdateduringclusterexploration}. Henceforth, we call it \textit{tree skeleton}, denoted by~$\bdcalS$. We choose $(0,0)$ as root vertex. 

As it has been explained at the beginning of the previous section, the percolation cluster is the union of the tree skeleton and some other vertices left behind. 
This suggests to divide the proof of Theorem~\ref{theorem:scalinglimittoymodel} in two parts. Following a common method---see for instance~\cite{duquesne2003limit, le2005random}---based on the coding of rooted trees via discrete processes, we are going to demonstrate that the (properly rescaled) tree skeleton admits some close to \eqref{theorem:scalinglimittoymodelconvergenceequality} scaling limit, as the length $T$ of the exploration tends to infinity. Then, we will extend the convergence to the whole cluster, by observing that the "left-behind" vertices are wiped out at large scale, because are confined in small areas. Finally, we use results of Section~\ref{sec:volumepercocluster} to transform the conditioning on $T$ into a conditioning on~$\lvert \tildebdcalC \rvert$.

\begin{paragraph}{Scaling limit of the tree skeleton}

The purpose of the paragraph is to establish the following convergence, holding in distribution at $p=p_c$ for the Gromov--Hausdorff distance:
\begin{align}
\label{eq:scalinglimitSCTHPtreeskeleton}
\big( n^{-1/2} \cdot \bdcalS \ \big\vert \ T \geq n \big) \xrightarrow[n \to +\infty]{d_{GH}} \sigma \cdot \bdcalT_{\geq 1},
\end{align}
where $\sigma=\sigma(\alpha)$ is the standard deviation of $h(\bdcalV_1)$ and $\bdcalT_{\geq 1}$ is the CRT of mass greater than~1.

Let us start by clarifying notions and objects. Any finite (or \textit{discrete}) tree can be viewed as a compact metric space by endowing its set of vertices with the graph distance.  
The Gromov--Hausdorff distance is a distance on the set of metric spaces, inducing some notion of convergence. If~$(E,\bfd)$ and $(E',\bfd')$ are two \emph{compact} metric spaces, the Gromov--Hausdorff distance between them is set as
\begin{align}
\label{eq:gromovhausdorffdistance}
\bfd_{GH}\Big((E,\bfd),(E',\bfd')\Big)=\inf\Big\{\bfd_H \big(\bdphi(E),\bdphi'(E')\big) \Big\},
\end{align} 
where the infimum is taken over all metric spaces $(F,\bddelta)$ and isometric embeddings $\bdphi: E \to F$ and~$\bdphi': E' \to F$, and where $\bfd_H$ designates the Hausdorff distance between compact sets in $F$.

The limit in distribution \eqref{eq:scalinglimitSCTHPtreeskeleton} of the tree skeleton w.r.t. $\bfd_{GH}$ is claimed to be the~CRT (up to a positive factor) of mass greater than $1$. The CRT is a \emph{real tree}, that is roughly a \emph{compact} metric space containing no cycle within it. A very convenient method to represent real trees is by using nonnegative real functions with compact support in $\bbR_+$. Let~$g$ be as such, satisfying also~$g(0)=0$. We set for every $0 \leq s \leq t$:
$$\bfd_{g}(s,t):=g(s)+g(t)-2 \inf_{u \in [s,t]} g(u),$$
and the following equivalence relation $\sim_g$ on $[0,+\infty[$:
$$s \sim_g t \Longleftrightarrow \bfd_g(s,t)=0.$$
The quotient space $\calT_g := [0,+\infty[ / \sim_g$ endowed with the distance $\bfd_g$ is a compact real tree. In Figure~\ref{figure:examplerealtreecodedbypiecewiselinearfunction}, we show how is derived a real tree from a piecewise linear function. The CRT of mass greater than $1$ is the real tree obtained by taking for $g=\bfe_{\geq 1}$, where $\bfe_{\geq 1}$ is the Brownian excursion  of duration greater than $1$ (we have~$\bfe_{\geq 1}(t)=0$ for~$t>t_{exc}$ and some~$t_{exc} \geq 1$). See \cite{duquesne2005probabilistic, le2005random} for more details.

Discrete trees---like $\bdcalS$---may also be viewed as real trees if we imagine them in the plane as an union of line segments of length one (the edges), equipped with the obvious distance (the length of the shortest path between two vertices). As real tree, the tree skeleton is isomorphic to that coded by a linearly interpolated version of the height process.
Set for instance for all $t \geq 0$:
\begin{align}
\label{eq:linearlyinterpolatedheightprocessdef}
\bdcalH(t) = \left\{
    \begin{array}{ll}
       \Big(h(\bdcalV_n) + (t-n)\cdot (h(\bdcalV_{n+1})-h(\bdcalV_{n}))\Big)_{+} & \mbox{for } n \leq T-1 \mbox{ and } t \in [n,n+1]; \\ 
       0 & \mbox{elsewhere}.
    \end{array}
\right.
\end{align}
Then, up to some isometry, we have:
\begin{align}
\label{eq:representationtreeskeletonasrealtree}
(\bdcalS,\bfd)=(\calT_{\bdcalH},\bfd_{\bdcalH}),
\end{align}
where $\bfd$ stands for the shortest-path distance in $\bdcalS$ viewed as real tree.
The statement is clearly supported by Figure~\ref{figure:downwardupwardrevelationsvertexupdateduringclusterexploration} and~\ref{figure:examplerealtreecodedbypiecewiselinearfunction}. 

Once \eqref{eq:representationtreeskeletonasrealtree} is known, the scaling limit \eqref{eq:scalinglimitSCTHPtreeskeleton} merely results from a suited application of the Donsker's theorem. Indeed, a consequence of the latter is that the continuous-time (and properly rescaled) height process defined in \eqref{eq:linearlyinterpolatedheightprocessdef} converges in distribution as follows:
\begin{align}
\label{eq:scalinglimitdonskertheorem}
\bigg( \big( n^{-1/2} \cdot \bdcalH (n \cdot t ) \big)_{t \geq 0} \ \bigg\vert \ T \geq n \bigg) \xrightarrow[n \to +\infty]{(d)} \sigma \cdot \bfe_{\geq 1},
\end{align}
where $\sigma=\sigma(\alpha)$ is the standard deviation of $h(\bdcalV_1)$ and $\bfe_{\geq 1}$ a Brownian excursion of duration greater than $1$. Then it becomes clear---by using for instance \cite[Lemma 2.4]{le2005random}---that:
\begin{align}
\big( n^{-1/2} \cdot \bdcalS \ \big\vert \ T \geq n \big)=\big( \calT_{n^{-1/2} \cdot \bdcalH} \ \big\vert \ T \geq n \big) \xrightarrow[n \to +\infty]{d_{GH}} \sigma \cdot \calT_{\bfe_{\geq 1}}, \nonumber
\end{align}
where the convergence holds in distribution for the Gromov--Hausdorff distance. This is exactly the expected conclusion. \qedd

\begin{figure}[!h]
\centering 
\subfigure[Graph of a piecewise linear $g$.]{
\centering
\includegraphics[scale=0.45]{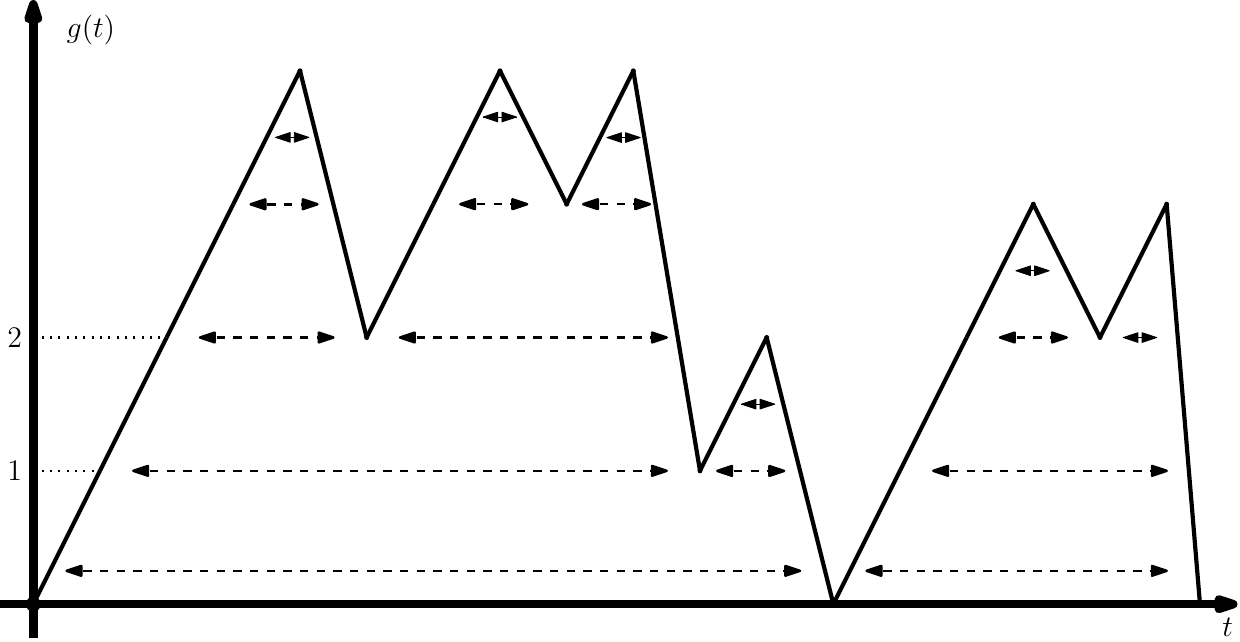}
}
\quad
\subfigure[The real tree coded by $g$.]{
\includegraphics[scale=0.45]{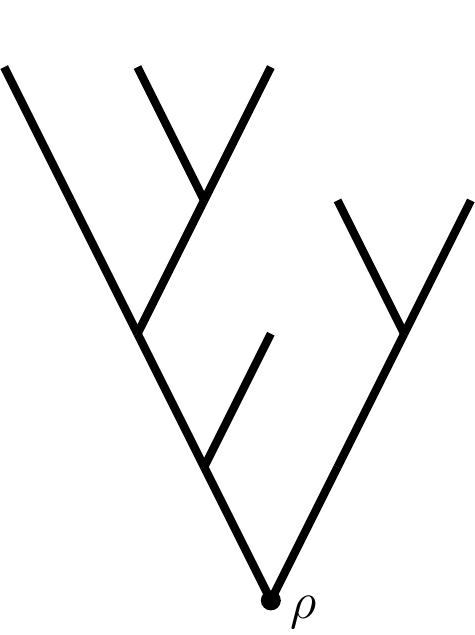}
}
\caption{Construction of a real tree via the coding by a piecewise linear function. Imagine that we apply glue all along the bottom side of the curve, then press it on both extremities. }
\label{figure:examplerealtreecodedbypiecewiselinearfunction}
\end{figure}
\end{paragraph}

\begin{paragraph}{A large cluster is roughly the tree skeleton}
Now we aim to extend \eqref{eq:scalinglimitSCTHPtreeskeleton} to the whole cluster.
It is enough to prove that
\begin{align}
\label{eq:GHdistancetreeskeletonpercoclustertendstozero}
\big( n^{-1/2} \cdot \bfd_{GH}\big(\tildebdcalC, \bdcalS \big) \ \big\vert \ T \geq n \big) \xrightarrow[n \to +\infty]{(d)} 0,
\end{align}
where the convergence holds in distribution. Since $\tildebdcalC$ and $\bdcalS$ are two subsets of the same SCTHP, we have:
$$\bfd_{GH}\Big(\tildebdcalC, \bdcalS \Big) \leq \bfd_{H}(\tildebdcalC, \bdcalS),$$
where $\bfd_H$ is the Hausdorff distance between compact sets of $\bdTheta$. See indeed \eqref{eq:gromovhausdorffdistance}. By construction, the tree skeleton is contained in the percolation cluster. Some vertices are in the latter, but not in the former. As it has already been explained, they are those left behind by the peeling exploration when occurs a downward revelation, or a drop in the height process. The distance between such vertices and the tree skeleton cannot exceed the amplitude of the jump carried out. See Figure~\ref{fig:exampleclusterexploration} or Figure~\ref{figure:downwardupwardrevelationsvertexupdateduringclusterexploration}. As a consequence:
\begin{align}
\label{eq:hausdorffdistpercoclustertreeskeletondominatesbymaxjump}
\bfd_{H}(\tildebdcalC, \bdcalS) \leq \max_{t \leq T-1} \big\{\lvert h(\bdcalV_{t+1})-h(\bdcalV_t)\rvert\big\}. 
\end{align} 
The distribution of $\lvert h(\bdcalV_1) \rvert$ is exponentially tailed. Then, for any $\epsilon >0$, any integer $n \leq n' \leq n^{1+\epsilon}$ and any real $C>0$:
\begin{align}
&\widetilde{\bbP}_{p_c}\bigg(\max_{t \leq T-1} \{\lvert h(\bdcalV_{t+1})-h(\bdcalV_t)\rvert\} > C \cdot \log{n}; \ T=n'\bigg) \nonumber \\ &\leq \widetilde{\bbP}_{p_c}\bigg(\max_{t \leq n'-1} \{\lvert h(\bdcalV_{t+1})-h(\bdcalV_t)\rvert\} > C \cdot \log{n}\bigg) \nonumber \\ &= 1 - \big(1-\widetilde{\bbP}_{p_c}(\lvert h(\bdcalV_1)\rvert > C \cdot \log{n}) \big)^{n^{1+\epsilon}}. \nonumber
\end{align} 
By using Propositions~\ref{prop:heightsubcriticalGWtrees} and~\ref{prop:distributionincrementsheightprocess}, we get the estimate
$$\widetilde{\bbP}_{p_c}\big(\lvert h(\bdcalV_1)\rvert > C \cdot \log{n}\big) = \calO\big( m^{-C \cdot \log{n}} \big) =  \calO\big(n^{- C \cdot \log{m}}\big).$$
Recall that $m = \alpha \cdot (1-\alpha)^{-1} > 1$. We can choose $C$ large enough so that $$\widetilde{\bbP}_{p_c}\big(\lvert h(\bdcalV_1)\rvert > C \cdot \log{n}\big)=o\big( n^{-5/2-2 \epsilon} \big).$$ A classic asymptotic expansion then ensures that:
\begin{align}
&1 - \big(1-\widetilde{\bbP}_{p_c}(\lvert h(\bdcalV_1)\rvert > C \cdot \log{n}) \big)^{n^{1+\epsilon}} \nonumber \\ &= 1-e^{n^{1+\epsilon}\cdot \log{(1-\widetilde{\bbP}_{p_c}(\lvert h(\bdcalV_1)\rvert > C \cdot \log{n}))}} =o\big( n^{-3/2-\epsilon} \big). \nonumber
\end{align}
From the foregoing, we deduce that:
\begin{align}
\label{eq:estimate1maxjump}
&\widetilde{\bbP}_{p_c}\bigg(\max_{t \leq T-1} \{\lvert h(\bdcalV_{t+1})-h(\bdcalV_t)\rvert\} > C \cdot \log{n}; \ n \leq T \leq n^{1+\epsilon} \bigg) \nonumber \\ &= \sum_{n \leq n' \leq n^{1+\epsilon}} \widetilde{\bbP}_{p_c}\bigg(\max_{t \leq T-1} \{\lvert h(\bdcalV_{t+1})-h(\bdcalV_t)\rvert\} > C \cdot \log{n}; \ T=n'\bigg) \nonumber \\ &\leq n^{1+\epsilon} \cdot \Big(1-\widetilde{\bbP}_{p_c}\big(\lvert h(\bdcalV_1)\rvert > C \cdot \log{n}\big) \Big)^{n^{1+\epsilon}} = o\big( n^{-1/2} \big)=o\Big( \widetilde{\bbP}_{p_c}(T \geq n) \Big).
\end{align}
Furthermore, according to Proposition~\ref{prop:taildistributionlengthnnexcursion}:
\begin{align}
\label{eq:estimate2maxjump}
\widetilde{\bbP}_{p_c}\big(T > n^{1+\epsilon}\big)=\calO\big( n^{-1/2 -\epsilon/2}\big)=o\big( n^{-1/2} \big)=o\Big( \widetilde{\bbP}_{p_c}(T \geq n) \Big).
\end{align}
Given \eqref{eq:estimate1maxjump} and \eqref{eq:estimate2maxjump}, we eventually derive that
$$\widetilde{\bbP}_{p_c}\bigg(\max_{t \leq T-1} \{\lvert h(\bdcalV_{t+1})-h(\bdcalV_t)\rvert\} > C \cdot \log{n} \ \vert \ T \geq n\bigg) \xrightarrow[n \to +\infty]{} 0,$$
for $C>0$ large enough.
With \eqref{eq:hausdorffdistpercoclustertreeskeletondominatesbymaxjump}, this yields the limit \eqref{eq:GHdistancetreeskeletonpercoclustertendstozero}.
\end{paragraph}

\begin{paragraph}{Scaling limit of the percolation cluster: proof of Theorem~\ref{theorem:scalinglimittoymodel}}
Let $F$ be a measurable and bounded function on the set of compact metric spaces, which is moreover continuous with respect to the Gromov--Hausdorff distance. We have:
\begin{align}
\label{eq:scalinglimiteqpercocluster}
\widetilde{\bbE}_{p_c} \big[ F(n^{-1/2} \cdot \tildebdcalC) \ \big\vert \ \lvert \tildebdcalC \rvert \geq n \big] &= \frac{1}{\widetilde{\bbP}_{p_c}(\lvert \tildebdcalC \rvert \geq n)} \cdot \widetilde{\bbE}_{p_c} \big[ F(n^{-1/2} \cdot \tildebdcalC) \cdot \mathbf{1}_{\lvert \tildebdcalC \rvert \geq n} \big] \nonumber \\ &= \frac{1}{\widetilde{\bbP}_{p_c}(T \geq n/\kappa)} \cdot \widetilde{\bbE}_{p_c} \big[ F(n^{-1/2} \cdot \tildebdcalC) \cdot \mathbf{1}_{T \geq n/\kappa} \big] +o(1).
\end{align}
In the last estimate, we use that $F$ is bounded and \eqref{eq:symmetricdifferenceclusterlengthexc}. The limit \eqref{eq:GHdistancetreeskeletonpercoclustertendstozero} implies that:
\begin{align}
\widetilde{\bbE}_{p_c} \big[ F(n^{-1/2} \cdot \tildebdcalC) \ \big\vert \ T \geq n/\kappa \big] \xrightarrow[n \to +\infty]{} \bbE[F(\kappa^{1/2} \cdot \sigma \cdot \bdcalT_{\geq 1})]. \nonumber
\end{align}
Together with \eqref{eq:scalinglimiteqpercocluster}, we get the expected conclusion for $\kappa'=\kappa^{1/2} \cdot \sigma$. \qedd
\end{paragraph}

\section{Percolation on SCT}
\label{sec:directedpercoonSCQ}
\begin{figure}[!h]
    \centering 
\subfigure[The leftmost (in blue) and the rightmost (in red) paths emerging from some vertex $v$ in a piece of a SCT~$\bdfrT$, represented as a cone.]{
\includegraphics[scale=0.3]{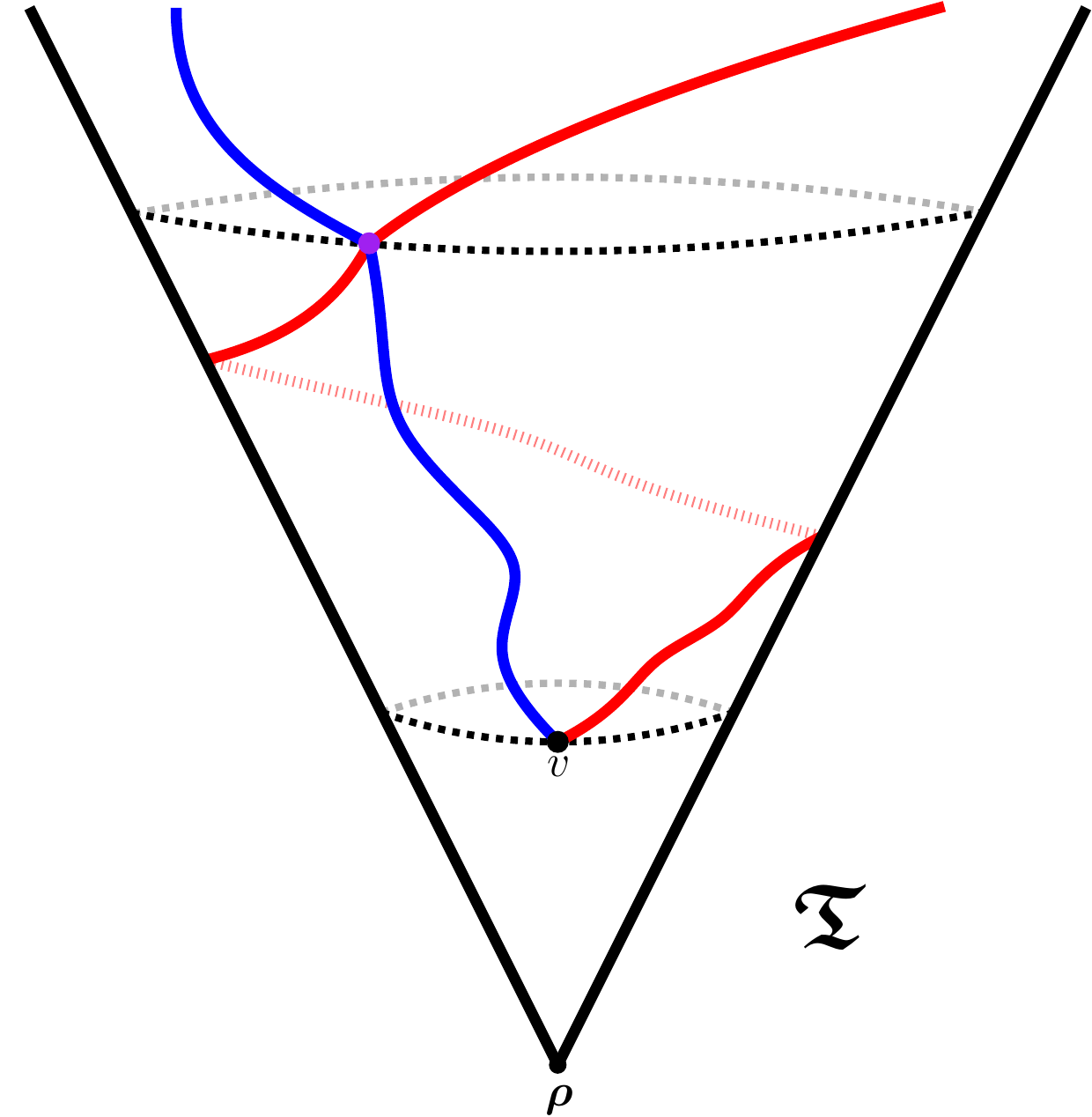}
}
\quad
\subfigure[The leftmost (in blue) and the rightmost (in red) paths emerging from $(0,0)$ in a piece of a SCTHP~$\tildebdfrT$.]{
\centering
\includegraphics[scale=0.3]{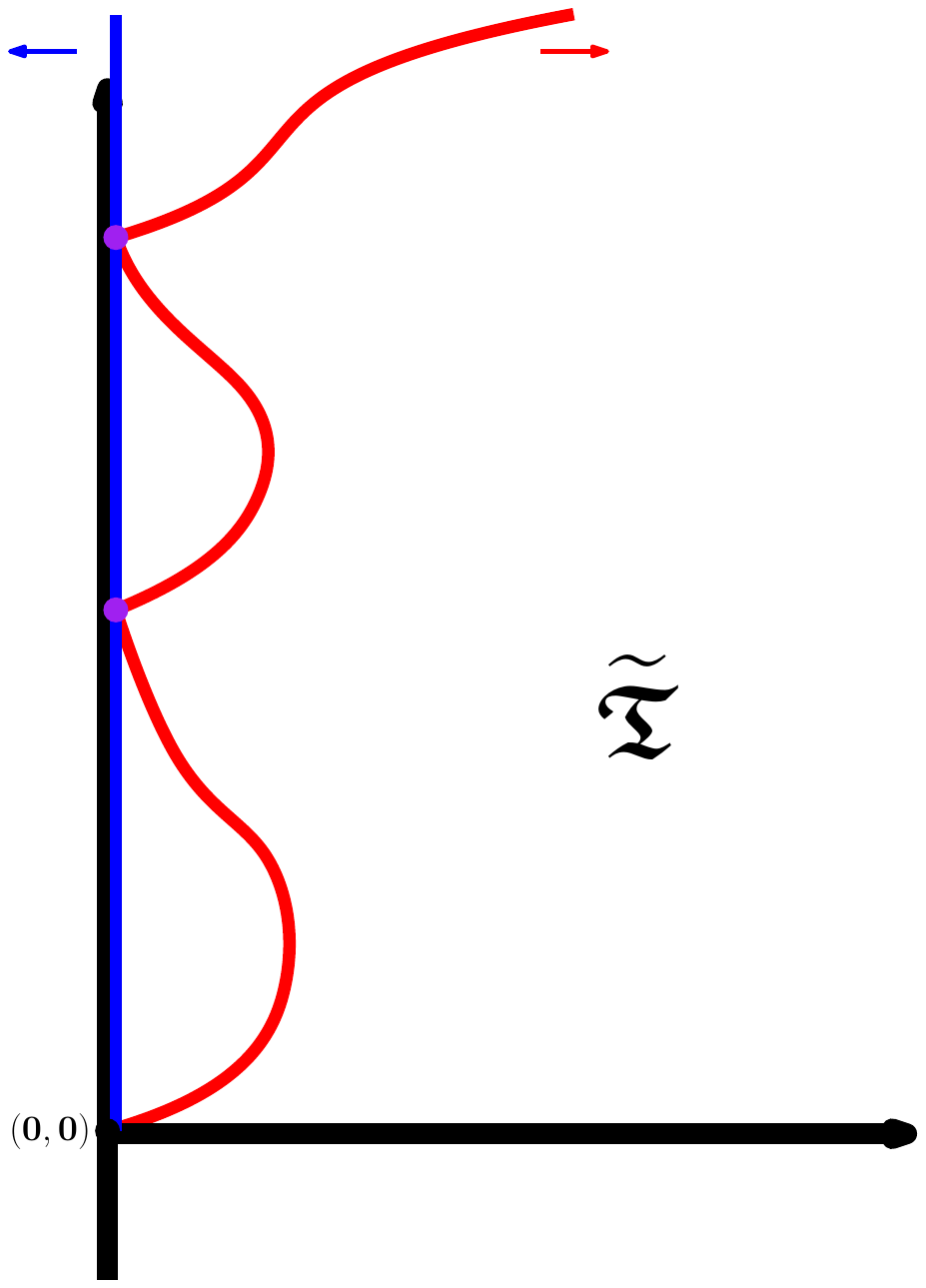}
}
    \caption{In a SCT, the leftmost path emerging from a vertex $v$ can at some point touch on its left flank the rightmost path, after "rotating" around the cone. Such situation is excluded in the SCTHP model which is infinite in the horizontal direction.}
\label{fig:exampleleftmostandrightmostpathsmerging}
\end{figure}

We aim to transfer back to the SCT model the results obtained in Section~\ref{sec:directedpercoontoymodelSCTHP} on SCTHP. Far away from the root, the landscape in a SCT around some given vertex looks a lot like in a SCTHP. There are however two significant differences. In the former, the rightmost path starting from some vertex $v$ can rotate around the underlying Galton--Watson tree, and finally merge up at some point with the leftmost path. Such a thing cannot happen in a SCTHP. See Figure~\ref{fig:exampleleftmostandrightmostpathsmerging}. The second important difference lies in the fact that we define a SCT from a supercritical Galton--Watson tree \textbf{conditioned to survive}. The conditioning kills the independance property of subtrees emerging from cousins, as it exists in the non conditioned version of the model or, to continue the comparison, between ascending trees in a SCTHP (see Proposition~\ref{prop:ascendingdescendingtreesHDHPareGW}). In Section~\ref{sec:decompositionGWhardboundariesSCQ}, we present a powerful trick to compensate such lack of independence and creates "hard boundaries" in SCT. We derive from it statements on the percolation probability function outlined in Theorem~\ref{theorem:phasetransitionSCQ} and Theorem~\ref{theorem:criticalexponentstheoremSCQ}~(ii), and also the existence in the supercritical phase of infinitely many disjoint infinite clusters. In Section~\ref{sec:largecriticalpercoclustersSCT}, we prove results on large critical clusters, namely Theorem~\ref{theorem:criticalexponentstheoremSCQ}(i) and the scaling limit of Theorem~\ref{theorem:scalinglimitSCQ}. The followed method significantly relies on a generalization of Theorem~\ref{theorem:criticalexponentstheoremSCTHP}~(ii) and Theorem~\ref{theorem:scalinglimittoymodel} to unions of critical clusters in a SCTHP, which are the subject of the closing Section~\ref{sec:intermediatesectionSCTHPtoSCT}.

\subsection{SCT and wide trees}\label{sec:decompositionGWhardboundariesSCQ}

\begin{figure}[!h]
    \centering 
\subfigure[A piece of a wide tree $\bfT$ rooted in $\rho$. In orange are the successively removed vertices and their associated subtrees.]{
\centering
\includegraphics[scale=0.3]{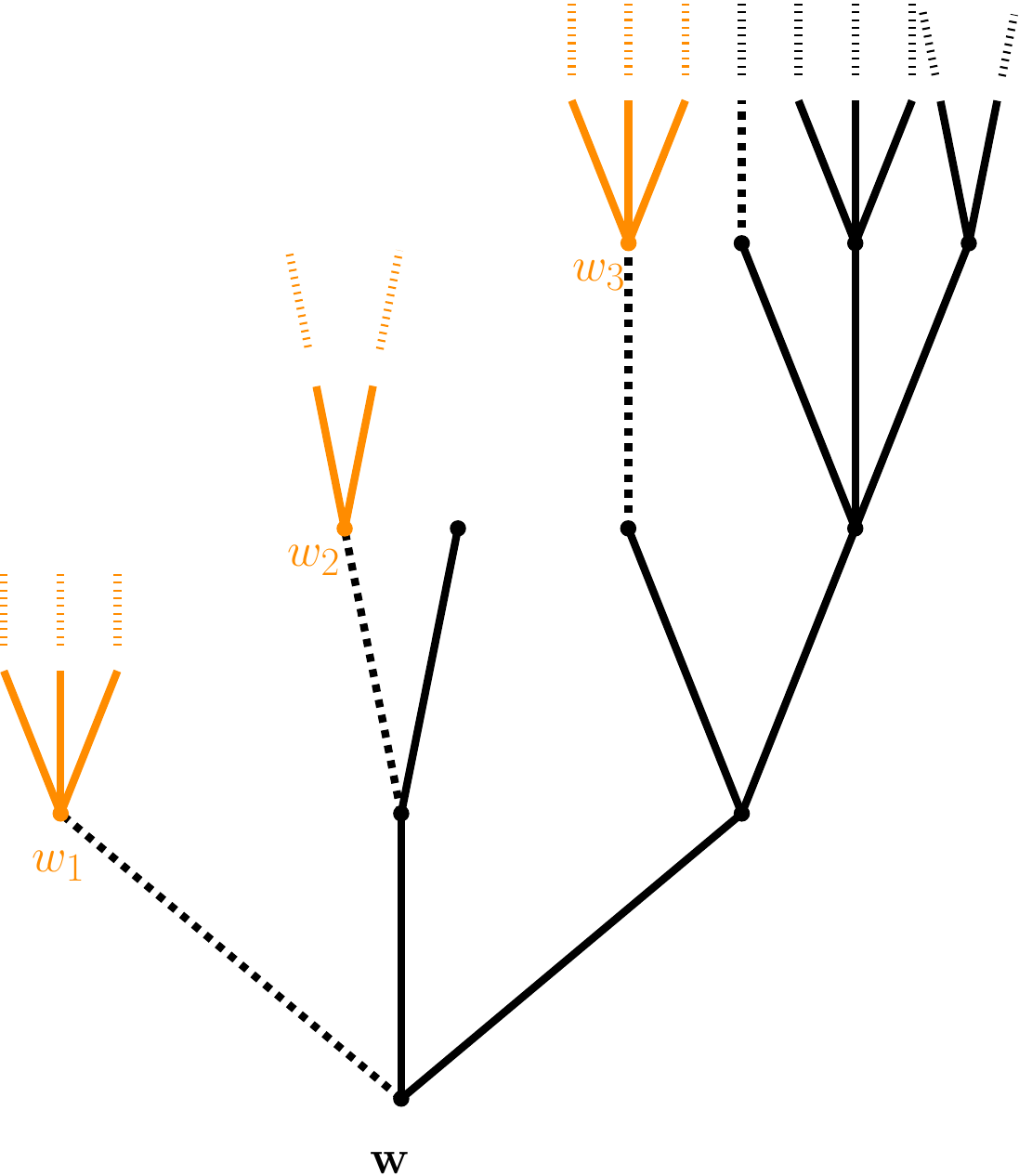}
}
\quad
\subfigure[A piece of a SCT $\bdfrT$ where a subtree $\bfT_v$ is a wide tree (filled in gray).]{
\includegraphics[scale=0.3]{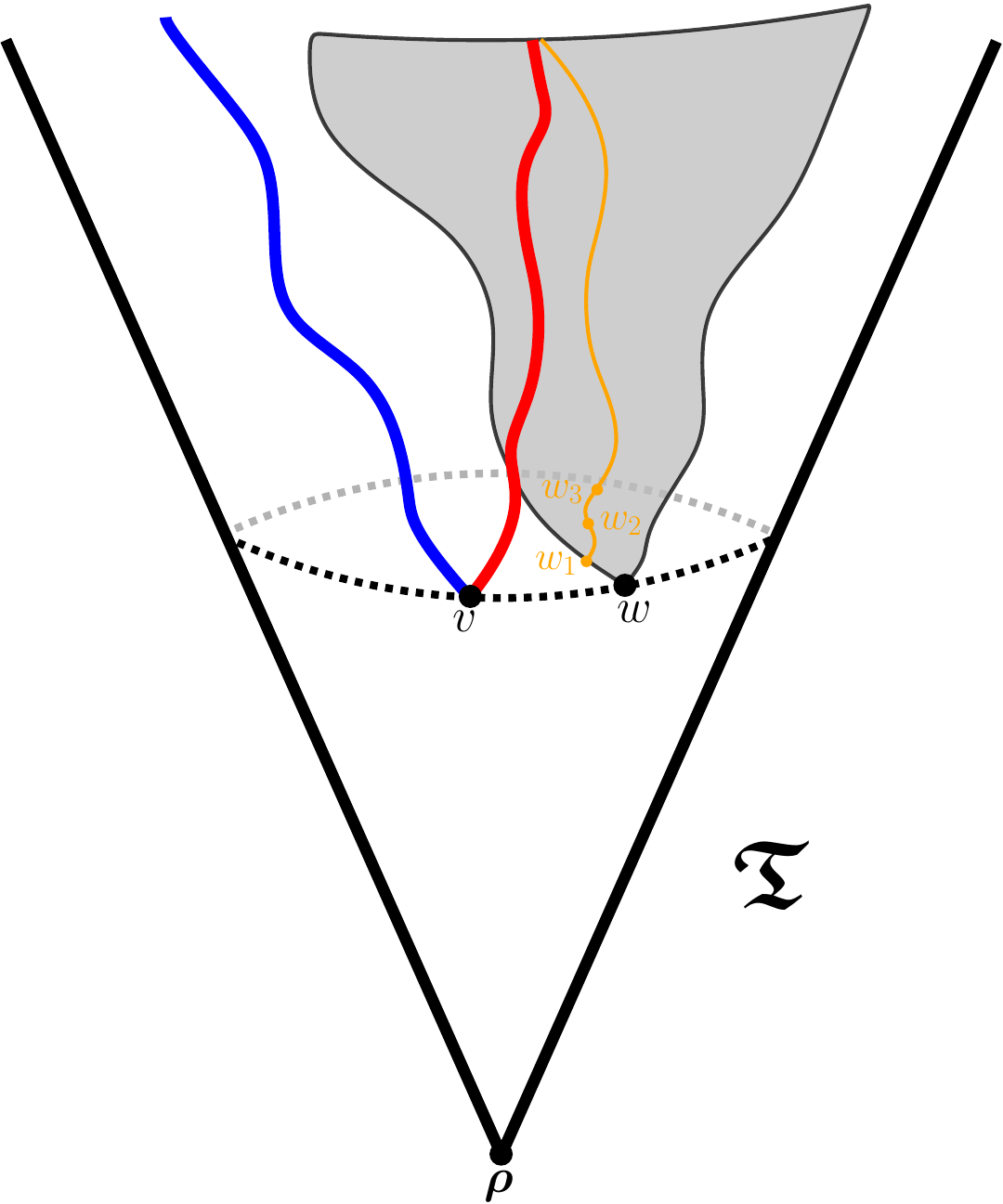}
}
    \caption{In a SCT, wide trees---like $\bfT_v$ on the right figure---create hard boundaries. A path emerging from any cousin of $v$ cannot overcome (but can reach) the orange line made of the successively removed vertices of $\bfT_v$, as we were checking whether the latter is wide or not. See for instance the rightmost path starting at~$w \neq v$, coloured in red.}
\label{fig:examplewidetreehardboundary}
\end{figure}

In a SCT, far away from the root vertex, the bad situation depicted in Figure~\ref{fig:exampleleftmostandrightmostpathsmerging}(a) does not occur in most cases. It is very likely that some middle subtree emerging from a cousin of $v$ is "wide enough" to prevent it. What is at stake is to clarify the term "wide enough". Take any (finite or infinite) plane tree $t$. Start by removing---if possible---the leftmost children of the root, as well as the subtree emerging from it. Iterate indefinitely the procedure by consistently removing, at each step, the leftmost children in $t$ of the vertices still present in the previous generation. If the subtree of $t$ obtained after such pruning, is infinite, we say that $t$ is \emph{wide}. See Figure~\ref{fig:examplewidetreehardboundary}(a) for an illustration of the pruning operation. In this example, the tree is wide since the part made of the edges in black continuous line---that we call the \emph{wide component}---is infinite.

When $t$ is a supercritical Galton--Watson tree conditioned to survive, with a geometric offspring distribution like in our model, it is wide with a positive chance:

\begin{lemma}
\label{lemma:probawidetreesuperGW}
\NoEndMark We have $\bbP(\bfT_\infty \text{ is \emph{wide}})>0$.
\end{lemma}

In fact, the result holds for any supercritical Galton--Watson tree, regardless of the offspring distribution. See \cite[Lemma 3]{budzinski2019supercritical} for a proof of the statement. The author actually deals with a more general mechanism, in which are removed at each step a certain number of children---not necessarily a single one like us, nor even a constant number.

As argued in the foregoing, our interest for wide trees results from our wish to prevent events like that represented on Figure~\ref{fig:exampleleftmostandrightmostpathsmerging}. Our observation is that wide trees fulfill their role by creating "hard boundaries" in the SCT as soon as they arise. Assume indeed that the subtree $\bfT_v$ emerging from some vertex~$v \in \bfT_\infty$ is a wide tree. In the SCT derived from~$\bfT_\infty$, any path starting from any cousin $w$ of $v$---a vertex which is at same height as $v$ in $\bfT_\infty$---can hit but cannot go through the infinite line made of the vertices bordering on its left the wide component of $\bfT_v$ (in orange on Figure~\ref{fig:examplewidetreehardboundary}).  Note that the latter line is by construction an existing path in the~SCT.

Furthermore, some stunning fact happens when is found---in a proper way, though---a wide tree in the SCT. The landscape around it is actually the same as that witnessed in a SCTHP. Indeed, if we condition for instance the underlying Galton--Watson tree $\bfT_\infty$---which is already conditioned to be infinite---so that at a given height $h \geq 1$, the subtree emanating from the leftmost vertex is wide, then the subtrees emanating from its cousins turn out to be i.i.d.~with common distribution~$\mathbf{GW}_\alpha$. This is simply because the forced condition on the former subtree is more binding than that on~$\bfT_\infty$ (to be infinite) and, in fact, bring with it. In doing so, we consequently lift all constraints on the other subtrees emanating from height $h$. Meanwhile, we remark that requiring a Galton--Watson tree to be wide does not demand anything on the subtrees emerging from the "removed vertices"---those in orange on Figure~\ref{fig:examplewidetreehardboundary}(a)---since they are by construction outside the wide component. The latter subtrees are then distributed as usual $\mathbf{GW}_\alpha$. 

To put it in a nutshell, as soon as we condition $\bfT_\infty$ on the event that at height $h$, the subtree emanating from the leftmost vertex is wide, the environment seen in the SCT by its cousins, looking upward, is as if they were in a SCTHP. In particular, their associated (directed) percolation cluster is distributed as $\tildebdcalC$, the percolation cluster of the origin in a SCTHP. 

\begin{paragraph}{Consequences on the percolation probability}
The foregoing observations are enough to prove first the lower bound of Theorem~\ref{theorem:criticalexponentstheoremSCQ}(ii).
Indeed, the event on which we condition occurs with a positive probability. Because the structure of $\bfT_\infty$ strictly above height~$h$ is independent of that below, conditionally on the population at height~$h$, it is still the case if in addition, we impose that there exists a second vertex at height $h$ and that the root of the SCT is connected to it. Given what it has just been argued, the probability that an infinite cluster emerges from the second vertex, conditionally on the previous event, is exactly~$\widetilde{\Theta}(p)$, where~$\widetilde{\Theta}$ is the percolation probability function of the SCTHP model. We deduce that 
\begin{align}
\label{eq:lowerboundprobabilitypercofunction}
\Theta(p) \geq a_1 \cdot \widetilde{\Theta}(p),
\end{align}
for some $a_1 \in (0,1)$. By using the third point of Theorem~\ref{theorem:criticalexponentstheoremSCTHP}, we conclude as expected that~$\Theta(p)$ is greater than a linear function in $p-p_c$ for $p>p_c$. It also implies that~$\Theta(p)>0$ for such~$p$, as stated in Theorem~\ref{theorem:phasetransitionSCQ}.

Obviously, the subtree emanating from the leftmost vertex is wide with a non zero probability, but not equal to one. So our method does not work as it is to reverse the inequality~\ref{eq:lowerboundprobabilitypercofunction} between~$\Theta(p)$ and~$\widetilde{\Theta}(p)$, which would be enough to establish the upper bound of Theorem~\ref{theorem:criticalexponentstheoremSCQ}(ii).  Nevertheless, there exists a smart way to choose a vertex in $\bfT_\infty$---not just the leftmost at some given height---where we end up in a similar position:

\begin{lemma}
\label{lemma:procedurechoosewidetree}
\NoEndMark By some explicit procedure, we find almost surely a vertex $\bfv^* \in \bfT_\infty$ such that:
\begin{enumerate}
\item the subtree emanating from $\bfv^*$ is \emph{wide};
\item it is disconnected from the previous generation of vertices in $\bfT_\infty$ after the percolation process is performed on the SCT;
\item the number of cousins of $\bfv^*$ has a finite mean;
\item conditionally on $\bfv^*$ and its number of cousins, the subtrees in $\bfT_\infty$ emanating from the latter are i.i.d.~with common distribution $\mathbf{GW}_\alpha$.
\end{enumerate}
\end{lemma}

Admit temporarily the lemma to see how it helps us to get the desired result. Given the second point, an infinite path emerging from the root vertex of the SCT necessarily visits some cousin~$w$ of~$\bfv^*$. The percolation cluster associated to $w$ is then infinite. Thanks to the first and fourth points, it is distributed as the cluster of the origin in a SCTHP, conditionally on $\bfv^*$ and its number of cousins. By an union bound argument, we obtain:
$$\Theta(p) \leq \widetilde{\Theta}(p) \cdot \bbE\big[ \text{number of cousins of $\bfv^*$} \big].$$
Then we use the third point of Lemma~\ref{lemma:procedurechoosewidetree}, together with Theorem~\ref{theorem:criticalexponentstheoremSCTHP}, to get the upper bound of Theorem~\ref{theorem:criticalexponentstheoremSCQ}(ii). Note that the above inequality also implies that $\Theta (p) = 0$ for~$p \leq p_c$, which completes the proof of the first part of Theorem~\ref{theorem:phasetransitionSCQ}. \qedd

\begin{prove}[Proof of Lemma~\ref{lemma:procedurechoosewidetree}]
The broad idea of the proof is that a wide subtree can always be found in~$\bfT_\infty$ by going far enough from the root vertex. We proceed as follows. We start by choosing a vertex of the tree $\bfT_\infty$, which is---as required in Lemma~\ref{lemma:procedurechoosewidetree}---disconnected from the previous generation of vertices after the percolation process is performed on the SCT $\bdfrT$. For instance, the leftmost one at minimal distance from the root vertex. Given the independance between the percolation process and the map, it is clear that such vertex always exists. We denote it by $\bfv_1$. 

The distance between~$\bfv_1$ and the root vertex is typically very short. It explains why the number of cousins of~$\bfv_1$ has a finite mean. Indeed, on the event~$\{\text{height of } \bfv_1=h\}$, a significant amount of edges cannot be \emph{closed}, more exactly as much as the number of vertices below height~$h-1$ having only one ancestor in the map, a number which is roughly equal to the total population at height~$h-1$. 
So, the number of cousins exponentially grows with the height \cite{kesten1966limit, lyons1995conceptual}, admittedly, but open a quantity of edges of same order has an overtaking super-exponential cost. 

Now we reveal whether or not the subree $\bfT_{\bfv_1}$ emanating from $\bfv_1$ is \emph{wide}. It is with positive probability, and if so, we stop the algorithm and simply set $\bfv^*=\bfv_1$. Indeed, the first three conditions of Lemma~\ref{lemma:procedurechoosewidetree} are obviously satisfied by $\bfv_1$. The fourth too, because the wideness of the subtree ensures that the underlying Galton--Watson tree $\bfT_\infty$ is indeed infinite. We had already met this argument to prove the inequality~\eqref{eq:lowerboundprobabilitypercofunction}.

When the subtree of $\bfv_1$ is not wide, its \emph{wide component}---that is the black part of the tree depicted on Figure~\ref{fig:examplewidetreehardboundary}(a)---is finite. Assume that its height is $H \geq 0$ and that of~$\bfv_1$ is~$h \geq 1$. Then we repeat the pattern that we have just described above height $h+H+1$. More precisely, we look now for a vertex~$\bfv_2$, at distance at least $h+H+1$ from the root vertex, disconnected from the previous generation of vertices after the percolation process is performed on the SCT, like was~$\bfv_1$. We choose again the leftmost one, with minimal height, and check whether or not the subtree $\bfT_{\bfv_2}$ is wide. If it is, then~$\bfv^*=\bfv_2$. Otherwise, we go further in the tree $\bfT_\infty$, and so on.

We must be certain that a vertex $\bfv^*$ satisfying the four conditions of Lemma~\ref{lemma:procedurechoosewidetree} will eventually be found after a finite number of steps. Vertices above height $h+H+1$ all belong to subtrees emerging either from cousins of~$\bfv_1$ or from the orange "removed vertices" of~$\bfT_{\bfv_1}$ (see Figure~\ref{fig:examplewidetreehardboundary}(a)). Conditionally on their number, these subtrees are i.i.d.~Galton--Watson trees with common distribution~$\mathbf{GW}_\alpha$. At least one of them has however to be infinite, in order to guarantee that the underlying tree~$\bfT_\infty$ is itself infinite. Their structure is independent of that of~$\bfT_\infty$ below height~$h$, and also of the finite wide component of~$\bfT_{\bfv_1}$---again, conditionally on their number. 
Therefore, the environment on which is performed the second step of the algorithm is (conditionally) independent of the environment revealed during the initial step. Its distribution is not that of the environment seen from the root vertex as we were looking for~$\bfv_1$, since the latter is built on a single Galton--Watson tree conditioned to be infinite, namely~$\bfT_\infty$. Nevertheless, the probability for~$\bfT_{\bfv_2}$ to be wide is lower bounded by a positive constant, uniformly in the number of subtrees emerging from height~$h+H+1$. By iterating the argument---in the unfavourable case where~$\bfT_{\bfv_2}$ is not wide---and by using the second Borel--Cantelli lemma, we conclude that the algorithm stops almost surely. \qedd
\end{prove}
\end{paragraph}

\begin{paragraph}{The number of clusters in the supercritical phase}
We close the section with the demonstration of the last part of Theorem~\ref{theorem:phasetransitionSCQ}, claiming that the percolated map is fragmented into infinitely many infinite \emph{non-directed} clusters, as $p_c <p<1$.
We insist on the fact that the statement holds for the \textit{non directed} Bernoulli percolation process, contrary to the other results. 

Lemma~\ref{lemma:procedurechoosewidetree} tells that~$\bfT_\infty$ always contains a wide subtree. There exists actually far more than that. At some given height, a positive proportion of the emerging subtrees are wide. Given the exponential growth of~$\bfT_\infty$, it shapes more and more regions, isolated from each other, where new infinite clusters can expand freely. 

To be more specific, fix~$n \geq 1$ and let~$N \geq 1$ be the number of vertices at height~$n$. 
Assume that some vertex~$v$ among the latter gives birth to (at least) four individuals. Force now the percolation process to \emph{close} the three edges connecting~$v$ to its three rightmost children. Then we impose that subtrees emanating from the rightmost and the third rightmost are \emph{wide}. After such operation, the directed cluster of the child located between them, that we denote by $w$, is (almost) completely disconnected from the rest of the map. Subtrees on the flanks being wide, its distribution is that of~$\tildebdcalC$. 
Therefore, the cluster is infinite with probability~$\widetilde{\Theta}(p)$, independently of the other geometric requirements on~$v$ and its progeny. The probability of the whole scenario is~$q':=\alpha^4 \cdot (1-p)^3 \cdot \bbP(\bfT \text{ is wide})^2 \cdot \widetilde{\Theta}(p)>0.$

By construction, the \emph{non directed} clusters crossing, at height~$n+1$, vertices in the same configuration as~$w$, are disjoint.  For any $1 \leq K \leq N$, find at least $K$ such vertices occurs with probability
$$1-\sum_{k=0}^{K-1} \binom{N}{k} (q')^k \cdot (1-q')^{N-k},$$
which tends to one as $N \to +\infty$. 

The conclusion comes from the fact that~$N$ diverges with~$n$ in~$\bfT_\infty$. Getting more than~$K$ vertices like~$w$ at height~$n+1$ becomes certain as~$n \to +\infty$.
It means that the SCT contains at least~$K$ disjoint infinite non-directed clusters with probability one. Since the latter is true for any~$K \geq 1$, there are in fact infinitely many.
\end{paragraph}

\subsection{Large critical percolation clusters in SCT}\label{sec:largecriticalpercoclustersSCT}

In this section, we demonstrate the first item of Theorem~\ref{theorem:criticalexponentstheoremSCQ}, which provides an asymptotic on the tail distribution of the size of a critical cluster, then the scaling limit of Theorem~\ref{theorem:scalinglimitSCQ}. We mainly aim to show that clusters in a SCT are close in distribution to those observed in a SCTHP. The strategy employed strongly relies on the quick emergence of wide trees in~$\bfT_\infty$. We will make reference at several stages to a result---see Theorem~\ref{theorem:largecriticalmultipleclusters}---on large unions of critical clusters in a SCTHP, which generalizes in a sense Theorem~\ref{theorem:criticalexponentstheoremSCTHP}~(ii) and Theorem~\ref{theorem:scalinglimittoymodel}. We choose to postpone to Section~\ref{sec:intermediatesectionSCTHPtoSCT} the rigorous assertion, as well as the many technicalities of its proof, since it would unnecessarily burden the reasoning set out below. 

\begin{paragraph}{The volume of a critical cluster in a SCT}
We start by noting that we easily obtain some equivalent of~\eqref{eq:lowerboundprobabilitypercofunction} involving $\bbP_{p_c}(\lvert \bdcalC \rvert  \geq n)$, that is:
\begin{align}
\label{eq:lowerboundtaildistributionpercoclusterSCT}
\bbP_{p_c}(\lvert \bdcalC \rvert  \geq n) \geq a_1 \cdot \widetilde{\bbP}_{p_c}(\lvert \tildebdcalC \rvert  \geq n),
\end{align}
where $a_1 \in (0,1)$. Like for~\eqref{eq:lowerboundprobabilitypercofunction}, the idea is to condition on a smartly chosen event having a non-zero chance to occur, under which the mass of $\bdcalC$ is mostly supported by some sub-cluster distributed as~$\tildebdcalC$. This can be achieved for instance, by conditioning on the joint events "the root vertex is only connected to its second leftmost children" and "the subtree emanating from the leftmost children of the root vertex of $\bfT_\infty$ is \emph{wide}". Indeed, the former event ensures that the size of the percolation cluster emanating from the second leftmost children is at least $n-1$, while the former event implies that its distribution is that of $\tildebdcalC$, as argued in the previous section. Together with Theorem~\ref{theorem:criticalexponentstheoremSCTHP}(ii), the inequality~\eqref{eq:lowerboundtaildistributionpercoclusterSCT} implies that~$\bbP_{p_c}(\lvert \bdcalC \rvert  \geq n)=\Omega\big(n^{-1/2}\big)$.

The approach to get an upper bound on $\bbP_{p_c}(\lvert \bdcalC \rvert  \geq n)$---and thus complete the proof of Theorem~\ref{theorem:criticalexponentstheoremSCQ}~(i)---is to say that wide trees quickly appear in $\bfT_\infty$, that is not so far from the root vertex, so that most of the mass of $\bdcalC$ is concentrated in a region of the SCT akin to a SCTHP. We will clarify in a specific statement what it precisely means. Let us first introduce several notations. 

For any $h \geq 1$, we write from now on~$\bfB_h$ to designate the ball of radius $h$ of $\bfT_\infty$ centered at the root vertex, or by definition the set of vertices of $\bfT_\infty$ whose height is at most $h$. We denote by~$\partial \bfB_r$ its boundary, that is the set of vertices whose height is exactly $h$. Now we set for any~$N \geq 1$:
\begin{align}
\bfH_N := 1+\inf\{h \geq 1 \text{ s.t. } \lvert \partial \bfB_h \rvert > N \}, \nonumber
\end{align}
which is the height (up to a constant) where the population of vertices first exceeds~$N$. The random variable~$\bfH_N$ is almost surely finite because $\partial \bfB_h \to +\infty$ as $h \to +\infty$. Finally, for any~$h, N \geq 1$, we define the event $W_{h}^{N}$ that:
\begin{enumerate}
\item the inequality~$\lvert \partial \bfB_h \rvert > 2N$ holds;
\item there exists two vertices~$v$ and~$w$, among the $N$ leftmost vertices of $\partial \bfB_h$ for the former, among the~$N$ rightmost for the latter, which are disconnected from~$\partial \bfB_{h-1}$ (in the SCT subject to the percolation process) and such that the subtrees emanating from them are \emph{wide}.
\end{enumerate}
Otherwise said, under the event~$W_{h}^{N}$, the population at height~$h$ exceeds a certain given level, and two wide subtrees emerge, one from a vertex in the left half of~$\partial \bfB_h$, the other from the right half.  

We come to the heralded statement:

\begin{lemma}
\label{lemma:findtwowidetrees}
\NoEndMark
For any $\beta > 0$ large enough, we have:
\begin{align}
\label{eq:firstestimatesizeballboundary}
\bbP \bigg( m/2 \leq \frac{\lvert \partial \bfB_{\bfH_{\beta \log{n}}} \rvert}{\beta \log{n}} \leq 4m^2 \bigg) = 1-\underset{n \to \infty}{o}\big(n^{-1/2}\big), 
\end{align}
and for any $\delta>0$:
\begin{align}
\label{eq:secondestimatevolumeballsubpol}
\bbP\Big(\lvert \bfB_{\bfH_{\beta \log{n}}} \rvert \leq  n^{\delta}\Big) = 1-o(n^{-1/2}).
\end{align}
Moreover, there exists some $C>0$ such that:
\begin{align}
\label{eq:emergencewidetrees}
\bbP_{p_c} \Big(W_{\bfH_{\beta \log{n}}}^{C \log{n}} \Big) = 1-\underset{n \to \infty}{o}(n^{-1/2}).
\end{align}
\end{lemma}
We postpone the proof of the lemma to the end of the current section. The asymptotics \eqref{eq:secondestimatevolumeballsubpol} and~\eqref{eq:emergencewidetrees} together say that as (at least) two wide trees emerge in $\bfT_\infty$ at the (random) height~$\bfH_{\beta \cdot \log{n}}$ with probability~$1-~o(n^{-1/2})$, the volume of $\bfT_\infty$ standing below remains subpolynomial in $n$. In such a situation, most of the mass of a large cluster $\bdcalC$---that is whose the size is larger than~$n$---would necessarily be concentrated in a region of the SCT looking like a SCTHP. 

Given our preliminary result affirming that~$\bbP_{p_c}(\lvert \bdcalC \rvert  \geq~n)=~\Omega\big(n^{-1/2}\big)$, the main consequence of Lemma~\ref{lemma:findtwowidetrees} is that we can work under the event
$$E=E_{\beta,n}=E_{\beta,C,\delta,n}:=\Big( m/2 \leq \frac{\lvert \partial \bfB_{\bfH_{\beta \log{n}}} \rvert}{\beta \log{n}} \leq 4m^2 \Big) \cap \Big(\lvert \bfB_{\bfH_{\beta \log{n}}} \rvert \leq  n^{\delta}\Big) \cap W_{\bfH_{\beta \log{n}}}^{C \log{n}}$$
to study a large cluster, since it implies that:
\begin{align}
\label{eq:firstasymptoticprooftaildistrcluster}
\bbP_{p_c}\big(\lvert \bdcalC \rvert  \geq n\big) \underset{n \to +\infty}{\sim} \bbP_{p_c}\big(\lvert \bdcalC \rvert  \geq n \ \big\vert \ E_{\beta,n} \big).
\end{align} 

On the event $E$, the tree $\bfT_\infty$ contains at most~$n^{\delta}$ vertices up to height $\bfH_{\beta \log{n}}$. If the cluster~$\bdcalC$ contains more than $n$ vertices, there would be at least $n-n^{\delta}$ among them, located above or at the height~$\bfH_{\beta \log{n}}$.
On~$E$ again, there exists two distinct vertices of $\partial \bfB_{\bfH_{\beta \log{n}}}$ from which emanate wide subtrees of $\bfT_\infty$. They are also disconnected from the previous generation in the percolated~SCT. Then, they cannot belong to~$\bdcalC$. Write~$\bfw_\ell$ and~$\bfw_r$ for the leftmost and the rightmost such vertices, respectively. Conditionally on them, subtrees arising from their cousins are mutually independent.
The distribution of those being on the left of $\bfT_{\bfw_\ell}$, resp.~on the right of~$\bfT_{\bfw_r}$, depends on whether or not their root vertex is connected to~$\partial \bfB_{\bfH_{\beta \log{n}}-1}$. For those which are not, they are distributed as~$\mathbf{GW}_\alpha$ but conditioned to be \emph{non wide}. For the others, the distribution is simply~$\mathbf{GW}_\alpha$. In the middle area, between~$\bfT_{\bfw_\ell}$ and~$\bfT_{\bfw_r}$, subtrees are all distributed as~$\mathbf{GW}_\alpha$.

The primary consequence of the last point is that the environment seen by vertices flanked by~$\bfw_\ell$ on their left and by~$\bfw_r$ on their right, is like in a~SCTHP (when they look upward). It is slightly more complex elsewhere given the \emph{non wide} conditioning on some subtrees. That brings us to split the part of~$\bdcalC$ exceeding height~$\bfH_{\beta \log{n}}$ into two pieces, denoted by $\bdcalC_{mid}$ and~$\bdcalC_{out}$. The former gathers sub-clusters emanating from vertices of the \emph{middle area} (between~$\bfw_\ell$ and~$\bfw_r$), while the latter does with those coming from the left of~$\bfT_{\bfw_\ell}$ or from the right of~$\bfT_{\bfw_r}$. 
The graphs~$\bdcalC_{mid}$ and~$\bdcalC_{out}$ are conditionally independent, because the subtrees on which they stand are themselves. 

The (random) number of vertices of~$\bdcalC$ intersecting the middle area is a measurable function of the graph structure below height $\bfH_{\beta \log{n}}$ of the percolated SCT. Given its value $N \geq 1$ and the occurrence of the event~$E$, the graph~$\bdcalC_{mid}$ is distributed as~$\tildebdcalC_{N}:=\bigcup_{i=1}^{N} \tildebdcalC_{v_i}$, where we write~$\tildebdcalC_{v_i}$ to designate the directed percolation cluster associated to the vertex~$v_i$ in a SCTHP, and where~$(v_i)_i$ is a sequence of~$N$ distinct vertices, all located on the~$x$-axis~$\bbN \times \{0\}$. In Section~\ref{sec:intermediatesectionSCTHPtoSCT}, we demonstrate---see Theorem~\ref{theorem:largecriticalmultipleclusters}---that~$\widetilde{\bbP}_{p_c}(\lvert \tildebdcalC_{N} \rvert  \geq n)$ is at most of order~$\calO\big(\log{(n)} \cdot n^{-1/2}\big)$ if the number~$N$ of clusters does not grow with~$n$ faster than logarithmically. Since~$N \leq  \lvert \partial \bfB_{\bfH_{\beta \log{n}}}\rvert \leq 4m^2 \beta \log{n}$ on the event~$E$, 
we get by using the latter result that:
\begin{align}
\label{eq:asymptclustermidout}
\bbP_{p_c}\big(\lvert \bdcalC_{mid}\rvert \geq n \ \big\vert \ E_{\beta,n} \big) = \calO\big(\log{(n)} \cdot n^{-1/2}\big).
\end{align}
The same asymptotic comparison holds even more for $\bdcalC_{out}$. Indeed, the potential \emph{non wide} conditioning on some subtrees outside the middle area merely pushes down (in distribution) their volume, and consequently, the number of vertices accessible in the SCT from their root vertex. 

Remember now that the total mass $\lvert \bdcalC_{mid} \rvert + \lvert \bdcalC_{out} \rvert$ is greater than~$n-n^{\delta}$ on the event $\big\{\lvert \bdcalC \rvert \geq n\big\}$. One of the two parts captures more than a half of it. From \eqref{eq:asymptclustermidout} and their conditional independence, we derive that it receives in fact, as $n \to +\infty$, an overwhelming proportion of the total mass. We mean that conditionally on both events~$E$ and $\big\{\lvert \bdcalC \rvert \geq n\big\}$, we have either 
$$\lvert \bdcalC_{mid} \rvert \geq n-2 n^{\delta} \quad \text{and} \quad \lvert \bdcalC_{out} \rvert \leq n^{\delta},$$ 
or the converse (by switching roles of $\bdcalC_{mid}$ and $\bdcalC_{out}$), with probability tending to one. Hence, we get from \eqref{eq:firstasymptoticprooftaildistrcluster}:
\begin{align}
\label{eq:secondasymptoticprooftaildistrcluster}
&\bbP_{p_c}\big(\lvert \bdcalC \rvert  \geq n \ \big\vert \ E_{\beta,n} \big) \nonumber \\ &\underset{n \to +\infty}{\sim} \bbP_{p_c}\big(\lvert \bdcalC_{mid} \rvert  \geq n-2n^\delta; \ \lvert \bdcalC \rvert  \geq n  \ \big\vert \ E_{\beta,n} \big)+\bbP_{p_c}\big(\lvert \bdcalC_{out} \rvert  \geq n-2n^\delta; \ \lvert \bdcalC \rvert  \geq n \ \big\vert \ E_{\beta,n} \big).
\end{align}
However, chances that $\bdcalC_{out}$ overtakes $\bdcalC_{mid}$ are slim. Indeed, as the event $E$ happens, the number of vertices on the left of $\bfw_\ell$ or on the right of $\bfw_r$ is at most $2C \log{n} -2$, while there are more than~$m \beta /2 \log{n} - 2C \log{n}$ in the middle area. By reading carefully through the proof of Theorem~\ref{theorem:largecriticalmultipleclusters} detailed in Section~\ref{sec:intermediatesectionSCTHPtoSCT}, it turns out that in a SCTHP, the mass of the union~$\tildebdcalC_{N}$ of critical clusters, conditioned to be large, is actually monopolized by a single sub-cluster. It is not at all shared equally. The same holds within $\displaystyle{\bdcalC_{mid} \cup \bdcalC_{out}}$. 
Since the distribution of a SCT is rotationally invariant, the probability that the sub-cluster arises from the middle area then increases up to one as $\beta \to +\infty$. It means that:
\begin{align}
\label{eq:chancesoutovertakesmid}
\lim_{\beta \to +\infty }\limsup_{n \to +\infty} \frac{\bbP_{p_c}\big(\lvert \bdcalC_{out} \rvert  \geq n-2n^\delta; \ \lvert \bdcalC \rvert  \geq n  \ \big\vert \ E_{\beta,n} \big)}{\bbP_{p_c}\big(\lvert \bdcalC_{mid} \rvert  \geq n-2n^\delta; \ \lvert \bdcalC \rvert  \geq n  \ \big\vert \ E_{\beta,n} \big)} =0.
\end{align}
Once again in Section~\ref{sec:intermediatesectionSCTHPtoSCT}---see Remark~\ref{remark:probalargemulticlustersaroundn}, we argue why the "$2n^\delta$" factor can be removed in the event involving $\bdcalC_{mid}$. Finally, given that~$\big\{\lvert \bdcalC_{mid} \rvert  \geq n\big\} \subset \big\{\lvert \bdcalC \rvert  \geq n\big\}$, we derive from the asymptotic~\eqref{eq:secondasymptoticprooftaildistrcluster} and the limit~\eqref{eq:chancesoutovertakesmid} the following inequalities:
\begin{align}
\label{eq:thirdasymptoticprooftaildistrcluster}
1 \leq \liminf_{n \to +\infty} \frac{\bbP_{p_c}\big(\lvert \bdcalC \rvert  \geq n \ \big\vert \ E_{\beta,n} \big)}{\bbP_{p_c}\big(\lvert \bdcalC_{mid} \rvert  \geq n \ \big\vert \ E_{\beta,n} \big)} \leq \limsup_{n \to +\infty}\frac{\bbP_{p_c}\big(\lvert \bdcalC \rvert  \geq n \ \big\vert \ E_{\beta,n} \big)}{\bbP_{p_c}\big(\lvert \bdcalC_{mid} \rvert  \geq n \ \big\vert \ E_{\beta,n} \big)} \leq 1+M_\beta,
\end{align}
with $M_\beta \xrightarrow[\beta \to +\infty]{} 0$.
Combined with \eqref{eq:asymptclustermidout}, this completes the proof of Theorem~\ref{theorem:criticalexponentstheoremSCQ}(i). \qedd
\end{paragraph}

\begin{paragraph}{Scaling limit of a critical cluster in a SCT}
We roughly adopt the same strategy to prove the scaling limit of Theorem~\ref{theorem:scalinglimitSCQ}.
Take $F$ a measurable and bounded function on the set of compact metric spaces, which is moreover uniformly continuous with respect to the Gromov--Hausdorff distance. On the event
$$E_{\beta,n} \cap \big\{\lvert \bdcalC \rvert \geq n \big\} \cap \big\{\lvert \bdcalC_{mid} \rvert \geq n-2n^{\delta} \big\} \cap \big\{\lvert \bdcalC_{out} \rvert \leq n^\delta \big\},$$
we have for $\delta <1/2$: 
$$\bfd_{GH}\Big(n^{-1/2} \cdot \bdcalC, n^{-1/2} \cdot \bdcalC_{mid} \Big) \underset{\bdcalC_{mid} \ \subset \  \bdcalC}{\leq}   n^{-1/2} \cdot 2 n^\delta \xrightarrow[n \to +\infty]{} 0.$$ 
Since $F$ is uniformly continuous, we deduce that:
\begin{align}
&\bbE_{p_c}\Big[F(n^{-1/2} \cdot \bdcalC) \cdot  \mathbf{1}_{\lvert \bdcalC_{mid} \rvert \geq n-2n^\delta} \mathbf{1}_{\lvert \bdcalC_{out} \rvert \leq n^\delta} \mathbf{1}_{\lvert \bdcalC\rvert \geq n}\ \Big\vert \ E_{\beta,n} \Big] \nonumber \\ &= \bbE_{p_c}\Big[F(n^{-1/2} \cdot \bdcalC_{mid}) \cdot \mathbf{1}_{\lvert \bdcalC_{mid} \rvert \geq n-2n^\delta} \mathbf{1}_{\lvert \bdcalC_{out} \rvert \leq n^\delta} \mathbf{1}_{\lvert \bdcalC\rvert \geq n}\ \Big\vert \ E_{\beta,n} \Big] + o\Big(\bbP_{p_c}\big(\lvert \bdcalC \rvert  \geq n\big)\Big). \nonumber
\end{align}
The same remains true by switching roles of $\bdcalC_{mid}$ and $\bdcalC_{out}$.
From the boundedness of $F$ and the arguments developed up to \eqref{eq:secondasymptoticprooftaildistrcluster}, we get that:
\begin{align}
\label{eq:firstasymptoticscalinglimitlargecluster}
&\bbE_{p_c}\Big[F(n^{-1/2} \cdot \bdcalC) \cdot \mathbf{1}_{\lvert \bdcalC\rvert \geq n}\Big] =\bbE_{p_c}\Big[F(n^{-1/2} \cdot \bdcalC_{mid}) \cdot \mathbf{1}_{\lvert \bdcalC_{mid} \rvert \geq n-2n^\delta} \cdot  \mathbf{1}_{\lvert \bdcalC\rvert \geq n}\ \Big\vert \ E_{\beta,n} \Big] \nonumber \\ &+ \bbE_{p_c}\Big[F(n^{-1/2} \cdot \bdcalC_{out}) \cdot \mathbf{1}_{\lvert \bdcalC_{out} \rvert \geq n-2n^\delta} \cdot  \mathbf{1}_{\lvert \bdcalC\rvert \geq n}\ \Big\vert \ E_{\beta,n} \Big] + o\Big(\bbP_{p_c}\big(\lvert \bdcalC \rvert  \geq n\big)\Big).
\end{align}
The boundedness of $F$ once again, the limit~\eqref{eq:chancesoutovertakesmid} and Remark~\ref{remark:probalargemulticlustersaroundn}---that we had already invoked to obtain~\ref{eq:thirdasymptoticprooftaildistrcluster}---together ensure that:
\begin{align}
\label{eq:upperboundlargeclusteroutsidearea}
\limsup_{n \to +\infty}\bigg\vert \frac{\bbE_{p_c}\big[F(n^{-1/2} \cdot \bdcalC_{out}) \cdot \mathbf{1}_{\lvert \bdcalC_{out} \rvert \geq n-2n^\delta} \cdot  \mathbf{1}_{\lvert \bdcalC\rvert \geq n}\ \big\vert \ E_{\beta,n} \big]}{\bbP_{p_c}\big(\lvert \bdcalC_{mid} \rvert  \geq n \ \big\vert \ E_{\beta,n}\big)} \bigg\vert  \leq M_\beta \cdot \vert\vert F \vert\vert_\infty,
\end{align}
where $M_\beta$ tends to $0$ as $\beta \to +\infty$.

As explained in the paragraph before~\eqref{eq:asymptclustermidout}, conditionally on the event~$E_{n, \beta}$, the graph $\bdcalC_{mid}$ is distributed as
$\tildebdcalC_{N}$. Recall that~$N$ corresponds to the number of vertices of $\bdcalC$ located in the middle area at height $\bfH_{\beta \log{n}}$, and is lower than~$4m^2 \beta \log{n}$ on the event~$E$. In Section~\ref{sec:intermediatesectionSCTHPtoSCT}---see Theorem~\ref{theorem:largecriticalmultipleclusters}, we prove that under such conditions, the following scaling limit holds:
\begin{align}
\nonumber
 \bbE_{p_c}\Big[F(n^{-1/2} \cdot \bdcalC_{mid}) \ \Big\vert \ \lvert \bdcalC_{mid} \rvert \geq n; \ E_{\beta,n} \Big] \xrightarrow[n \to +\infty]{} \bbE\Big[F(\kappa' \cdot \bdcalT_{\geq 1})\Big],
\end{align}
where $\kappa'=\kappa'(\alpha)$ is the positive number defined in Theorem~\ref{theorem:scalinglimittoymodel} and $\bdcalT_{\geq 1}$ is the CRT of mass greater than~1. Then:
\begin{align}
\label{eq:scalinglimitlargeclustermiddlearea}
&\bbE_{p_c}\Big[F(n^{-1/2} \cdot \bdcalC_{mid}) \cdot \mathbf{1}_{\lvert \bdcalC_{mid} \rvert \geq n-2n^\delta} \cdot  \mathbf{1}_{\lvert \bdcalC\rvert \geq n}\ \Big\vert \ E_{\beta,n} \Big] \cdot \bbP_{p_c}\big(\lvert \bdcalC_{mid} \rvert  \geq n \ \big\vert \ E_{\beta,n} \big)^{-1} \nonumber \\ &\underset{\text{Remark~\eqref{remark:probalargemulticlustersaroundn}}}{=}  \bbE_{p_c}\Big[F(n^{-1/2} \cdot \bdcalC_{mid}) \ \Big\vert \ \lvert \bdcalC_{mid} \rvert \geq n; \ E_{\beta,n} \Big] + o(1) \xrightarrow[n \to +\infty]{} \bbE\Big[F(\kappa' \cdot \bdcalT_{\geq 1})\Big].
\end{align}
The inequalities in \eqref{eq:thirdasymptoticprooftaildistrcluster}, the asymptotic \eqref{eq:firstasymptoticscalinglimitlargecluster}, the upper bound \eqref{eq:upperboundlargeclusteroutsidearea} and the limit \eqref{eq:scalinglimitlargeclustermiddlearea} finally together imply that:
\begin{align}
\limsup_{n \to +\infty} \bbE_{p_c}\Big[F(n^{-1/2} \cdot \bdcalC) \ \Big\vert \ \lvert \bdcalC \rvert \geq n \Big] \leq \bbE\Big[F(\kappa' \cdot \bdcalT_{\geq 1})\Big] + M_\beta \cdot \vert\vert F \vert\vert_\infty \nonumber
\end{align}
and 
\begin{align}
\liminf_{n \to +\infty} \bbE_{p_c}\Big[F(n^{-1/2} \cdot \bdcalC) \ \Big\vert \ \lvert \bdcalC \rvert \geq n \Big] \geq (1+M_\beta)^{-1} \cdot \bbE\Big[F(\kappa' \cdot \bdcalT_{\geq 1})\Big] - M_\beta \cdot (1+M_\beta)^{-1} \cdot \vert\vert F \vert\vert_\infty \nonumber
\end{align}
We conclude by taking $\beta \xrightarrow[]{} +\infty$ to get $M_{\beta} \xrightarrow[]{} 0$. \qedd 
\end{paragraph}

\begin{prove}[Proof of Lemma~\ref{lemma:findtwowidetrees}] 
It remains to prove the three asymptotics stated in Lemma~\ref{lemma:findtwowidetrees} to complete the work. We start with ~\eqref{eq:firstestimatesizeballboundary}. For sake of simplicity, we remove the conditioning on $\bfT_\infty$ (imposing to be infinite). We will write~$\bbP_\infty$ instead of~$\bbP$ to indicate when the conditioning is back. 
By definition of~$\bfH_{\beta \log{n}}$, we have~$\lvert \partial \bfB_{\bfH_{\beta \log{n}}-1} \rvert \geq~\lceil \beta \log{n} \rceil$. So:
\begin{align}
\bbP\Big(\lvert \partial \bfB_{\bfH_{\beta \log{n}}} \rvert < \frac{m}{2} \beta \log{n} \ \Big\vert \ \bfH_{\beta \log{n}}<+\infty \Big) \leq \bbP\Big(\sum_{i=1}^{ \lceil \beta \log{n} \rceil} \bfX_i < \frac{m}{2} \beta \log{n}\Big), \nonumber
\end{align}
where the $(\bfX_i)_i$ are i.i.d.~random variables with geometric distribution $\bdmu_\alpha$.
According to the large deviations theory, the right-hand side of the inequality is~$o(n^{-1/2})$ for~$\beta$ large enough, because~$\bdmu_\alpha$ is an exponentially-tailed distribution. We deduce that:
\begin{align}
\label{eq:firstestimatelemmawidetreesfirstpart}
\bbP_\infty\bigg(\frac{\lvert \partial \bfB_{\bfH_{\beta \log{n}}} \rvert}{\beta \log{n}} < \frac{m}{2}\bigg) = \underset{n \to \infty}{o}\big(n^{-1/2}\big).
\end{align}
Similarly, by using this times that $\partial \bfB_{\bfH_{\beta \log{n}}-2} \leq \lfloor \beta \log{n} \rfloor$, we get for $\beta$ large enough:
\begin{align}
\bbP\Big(\lvert \partial \bfB_{\bfH_{\beta \log{n}-1}} \rvert > 2m \beta \log{n} \ \Big\vert \ \bfH_{\beta \log{n}}<+\infty \Big) \leq \bbP\Big(\sum_{i=1}^{ \lfloor \beta \log{n} \rfloor} \bfX_i > 2m \beta \log{n}\Big) = \underset{n \to \infty}{o}\big(n^{-1/2}\big), \nonumber
\end{align}
and:
\begin{align}
\label{eq:firstestimatelemmawidetreesintermediate}
&\bbP\Big(\lvert \partial \bfB_{\bfH_{\beta \log{n}}} \rvert > 4m^2 \beta \log{n} \ \Big\vert \ \bfH_{\beta \log{n}}<+\infty; \ \lvert \partial \bfB_{\bfH_{\beta \log{n}-1}} \rvert \leq 2m \beta \log{n}\Big) \nonumber \\ &\leq \bbP\Big(\sum_{i=1}^{\lfloor 2m \beta \log{n}\rfloor} \bfX_i > 4m^2 \beta \log{n}\Big)= \underset{n \to \infty}{o}\big(n^{-1/2}\big).
\end{align}
Since~$\bbP_\infty(\bfH_{\beta \log{n}}<+\infty)=1$, we derive from the two last asymptotics that:
\begin{align}
\label{eq:firstestimatelemmawidetreesfsecpart}
\bbP_\infty\bigg(\frac{\lvert \partial \bfB_{\bfH_{\beta \log{n}}} \rvert}{\beta \log{n}} > 4m^2\bigg) = \underset{n \to \infty}{o}\big(n^{-1/2}\big).
\end{align}
By combining with \eqref{eq:firstestimatelemmawidetreesfirstpart}, we obtain the expected asymptotic~\eqref{eq:firstestimatesizeballboundary}, 

We continue with the proof of \eqref{eq:secondestimatevolumeballsubpol}. Since $\lvert \partial \bfB_{h} \rvert \leq \beta \log{n}$ for every $h < \bfH_{\beta \log{n}}-1$, the volume of $\bfB_{\bfH_{\beta \log{n}}}$ is bounded by $$\bfH_{\beta \log{n}} \cdot \beta \log{n}+ \partial \bfB_{\bfH_{\beta \log{n}}}.$$
Hence, as we condition on the event $\big\{\lvert \partial \bfB_{\bfH_{\beta \log{n}}} \rvert \leq 4m^2 \beta \log{n} \big\}$, if $\lvert \bfB_{\bfH_{\beta \log{n}}} \rvert > n^{\delta}$ for some~$\delta >0$, it holds that
$\bfH_{\beta \log{n}} > n^{\delta /2}$ for $n$ large enough. In particular, this means that $\partial \bfB_{n^{\delta/2}-1} \leq \beta \log{n}$. Now we state that if a function~$\Phi : \bbN^* \mapsto ]0,+\infty[$ satisfies $\Phi(n)=o(m^n)$, then we have:
\begin{align}
\bbP\big(1 \leq \lvert \partial \bfB_{n} \rvert \leq \Phi(n)\big)=\calO\big(\Phi(n) \cdot m^{-n}\big). \nonumber
\end{align}
The asymptotic can be proved via an explicit computation of the generative function of~$\lvert \partial \bfB_{n} \rvert$, by mimicking the strategy employed in the proof of Proposition~\ref{prop:heightsubcriticalGWtrees}. Just set $u_0=s$ for $s \in [0,1]$ instead of $u_0=0$. Otherwise, see \cite{fleischmann2007lower} where the authors deal with such asymptotic in a much wider context. In our situation, the above statement implies that:
\begin{align}
&\bbP\Big(\lvert \bfB_{\bfH_{\beta \log{n}}} \rvert > n^{\delta} \ \Big\vert \ \lvert \partial \bfB_{\bfH_{\beta \log{n}}} \rvert \leq 4m^2 \beta \log{n}; \ \bfH_{\beta \log{n}}<+\infty \Big) \nonumber \\ 
&\underset{\text{for $n$ large enough}}{\leq} \bbP\big(\partial \bfB_{n^{\delta/2}-1} \leq \beta \log{n}\big) = \calO\big(\log{n} \cdot m^{-n^\delta}\big)=o(n^{-1/2}), \nonumber
\end{align}
for any $\delta>0$. Given \eqref{eq:firstestimatelemmawidetreesfsecpart}, this completes the proof of \eqref{eq:secondestimatevolumeballsubpol}.

We finally turn our attention to the asymptotic \eqref{eq:emergencewidetrees}. Again, we temporarily remove the conditioning on $\bfT_\infty$. 
Assume that $\partial \bfB_{h}> \beta \log{n}$ for some height $h \geq 1$ and large~$n$. Choose a vertex~$v$ in $\bfB_{h}$ and assume that it gives birth to (at least) two children. Then, we force the edge connecting~$v$ to its rightmost child to be \emph{closed}, after the percolation process has been performed. We also impose that the subtree emanating from the latter is \emph{wide}. Since the closed edge was the only possible path in the~SCTfrom~$\partial \bfB_h$ to~$v$, we have just isolated a vertex of $\partial \bfB_{h+1}$ from the previous generation, with a wide subtree arising from it. By Lemma~\ref{lemma:probawidetreesuperGW}, it has of course a positive probability to occur, namely $$q:=~\alpha^2 (1-p) \cdot~\bbP(\bfT \text{ is \emph{wide}}),$$ 
where $\bfT \overset{(d)}{\sim} \mathbf{GW}_{\alpha}$. The probability that such event happens for at least one of the $K \leq \beta \log{n}$ leftmost vertices of $\partial \bfB_h$ is~$1-(1-~q)^K$. Consequently, it does with probability~$1-o(n^{-1/2})$ when~$K=~c \log{n}$ for sufficiently large $0 <c< \beta$. The leftmost vertex for which it occurs has, by definition, at most $c \log{n}$ vertices on its left in~$\partial \bfB_h$---including itself. The number of vertices on the left of its rightmost child in $\partial \bfB_{h+1}$ is thus stochastically dominated by the sum of $\lfloor c \log{n} \rfloor$ i.i.d.~random variables with common distribution~$\bdmu_\alpha$. By a large deviation argument
---see \eqref{eq:firstestimatelemmawidetreesintermediate}, the sum does not exceed $2m \cdot~c \log{n}$ with probability~$1-o(n^{-1/2})$. It ensures that the rightmost child is among the~$2m \cdot~c \log{n}$ leftmost vertices of~$\partial \bfB_{h+1}$. A symmetric argument allows us to find another such vertex among the~$2m \cdot~c \log{n}$ rightmost vertices. This puts an end to the demonstration given that~$\bbP_\infty(\bfH_{\beta \log{n}}<+\infty)=1$. \qedd
\end{prove}

\subsection{Unions of critical percolation clusters in a SCTHP}\label{sec:intermediatesectionSCTHPtoSCT}

The next pages are dedicated to describe more fully a result on unions of critical percolation clusters in SCTHP to which we have repeatedly referred in the previous section. It can be read as a generalization of the second item of Theorem~\ref{theorem:criticalexponentstheoremSCTHP} and Theorem~\ref{theorem:scalinglimittoymodel}.  

Set~$\displaystyle{( v_k =~(i_k,0))_{k \geq 1}}$ an infinite sequence of vertices all located on the~$x$-axis~$\bbN \times \{0\}$, ordered from left to right---that is $i_k < i_{k+1}$ for every $k \geq 1$. We denote by $\tildebdcalC_{v_k}$ the directed cluster associated to the vertex~$v_k$, and write for any $N \geq 1$:
$$\tildebdcalC_{N}:=\bigcup_{1 \leq k \leq N} \tildebdcalC_{v_k},$$
which is the union of the~$N$ directed percolation clusters emanating from the~$N$ leftmost vertices of the sequence. 

We aim to look into the tail distribution of~$\tildebdcalC_{N}$ and its geometry at large scale in the critical regime~$p=p_c$, as the number~$N$ of clusters goes to infinity, but at a rate remaining slow compared to the volume of $\tildebdcalC_{N}$, though. Our statement is the following:

\begin{theorem}
\label{theorem:largecriticalmultipleclusters}
\NoEndMark
Let~$(N_n)_{n \geq 1}$ be a sequence of integers such that $N_n =~\calO\big( n^{\epsilon}\big)$ for some~$\epsilon \in~(0,1/20)$. We have that:
\begin{align}
\label{theorem:taildistributionmulitplecriticalpercoclusters}
\widetilde{\bbP}_{p_c}(\lvert \tildebdcalC_{N_n} \rvert  \geq n) = \calO\big(N_n \cdot n^{-\frac{1}{2}}\big)
\end{align}
and also:
\begin{align}
\label{theorem:scalinglimittoymodelconvergenceequalitymultipleclusters}
\big( n^{-1/2} \cdot \tildebdcalC_{N_n} \ \big\vert \ \vert \tildebdcalC_{N_n} \rvert \geq n \big) \xrightarrow[n \to +\infty]{d_{GH}} \kappa' \cdot \bdcalT_{\geq 1},
\end{align}
where $\kappa'=\kappa'(\alpha)$ is the positive number defined in Theorem~\ref{theorem:scalinglimittoymodel} and $\bdcalT_{\geq 1}$ is the CRT of mass greater than~1 \cite{aldous1991continuum, aldous1993continuum}. The convergence \eqref{theorem:scalinglimittoymodelconvergenceequalitymultipleclusters} holds in distribution for the Gromov--Hausdorff distance.
\end{theorem}

Neither the scaling limit \eqref{theorem:scalinglimittoymodelconvergenceequalitymultipleclusters} nor the asymptotic comparison~\eqref{theorem:taildistributionmulitplecriticalpercoclusters} depend on the specific choice of vertices~$v_k$ behind the definition of~$\tildebdcalC_{N}$, in the sense that for the latter, we obtain an universal upper bound on~$N_n^{-1} n^{1/2} \cdot \widetilde{\bbP}_{p_c}(\lvert \tildebdcalC_{N_n} \rvert  \geq n)$. Also, both assertions remain true by taking random numbers~$N_n=\bfN_n$, provided that they are generated independently of the map, the percolation process, and are still almost surely a~$\calO\big( n^{\epsilon}\big)$ as~$n \to +\infty$. Note finally that Theorem~\ref{theorem:criticalexponentstheoremSCTHP}(ii) trivially implies that $\widetilde{\bbP}_{p_c}(\lvert \tildebdcalC_{N_n} \rvert  \geq n)=\Omega\big(n^{-1/2}\big)$, since we have $$\widetilde{\bbP}_{p_c}(\lvert \tildebdcalC_{N_n} \rvert  \geq n) \geq \widetilde{\bbP}_{p_c}(\lvert \tildebdcalC_{v_1} \rvert  \geq n) = \widetilde{\bbP}_{p_c}(\lvert \tildebdcalC \rvert  \geq n).$$

The convergence~\eqref{theorem:scalinglimittoymodelconvergenceequalitymultipleclusters} may seem a bit confused, since~$\tildebdcalC_{N_n}$ will sometimes contain several disjoint components, and so not be compact. For two vertices $v, w \in \tildebdcalC_{N_n}$ belonging to disjoint components, say for instance associated to~$v_k$ and $v_{k'}$, we set $$\bfd(v,w):=\bfd(v,v_k)+\bfd(w,v_{k'}).$$
By doing so, we turn $(\tildebdcalC_{N_n},\bfd )$ into a compact metric space.

\begin{figure}[!h]
    \centering
\includegraphics[scale=0.5]{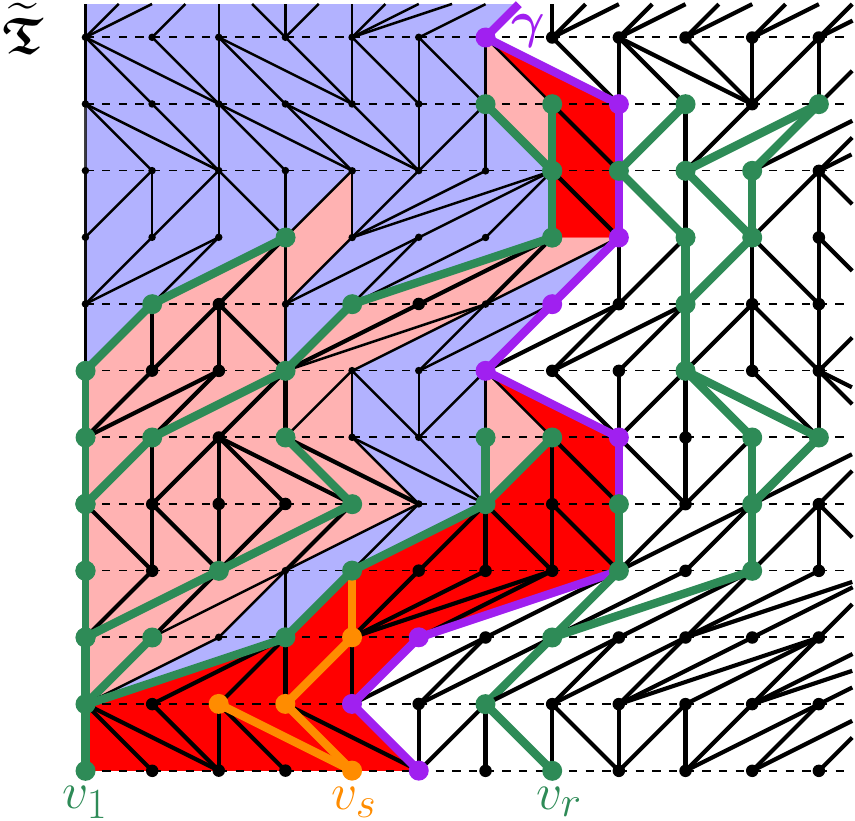}
    \caption{Illustration of the generalized peeling exploration. We start by exploring the tree skeleton of $v_1$ (in green). The revealed submap $\bdSigma_T$ is made of the blue and red areas, corresponding to the successively performed upward and downward revelations. From the vertex~$v_s$ on the left of the purple path~$\bdgamma$ (the right boundary of $\bdSigma_T$), we can only reach either~$\tildebdcalC_{v_1}$ or the bottom red area enclosed between the tree skeleton and~$\bdgamma$. Some vertices in the latter category are drawn in orange. The exploration then refreshes at~$v_r$. We repeat the above procedure to reveal the unknown part of~$\tildebdcalC_{v_r}$ located on the right of~$\bdgamma$. Some vertices in the interior of~$\bdSigma_T$ might be accessible from $v_r$, but they are necessarily in red areas bordering~$\bdgamma$.}
\label{fig:examplerefreshvertexgeneralizedpeelingexplo}
\end{figure}

\begin{paragraph}{An iterative peeling exploration}
The proof of the above theorem rests on a generalization of the peeling exploration designed in Section~\ref{sec:explorationclusterRW}. The principle is to discover one by one the clusters~$\tildebdcalC_{v_k}$ for $k \leq N_n$ as follows.

The first vertex $v_1$ may not be the origin vertex~$(0,0)$. In such case, we simply remove vertices inaccessible from $v_1$ (in a directed fashion) in the underlying SCTHP. It merely consists of erasing the ascending trees associated to vertices on the~$x$-axis which are on the left of~$v_1$. The new map obtained at the end of the operation is distributed as a SCTHP with~$v_1$ as new origin vertex---see indeed Section~\ref{sec:basicpropertiesofthetoymodelSCTHP}. Once the pruning achieved (if needed), we perform the usual peeling exploration of $\tildebdcalC_{v_1}$ described in Section~\ref{sec:explorationclusterRW}. 

Since we are in the critical regime~$p=p_c$, the exploration of~$\tildebdcalC_{v_1}$ ends at some time $T<+\infty$. The cluster~$\tildebdcalC_{v_1}$ is finite, considered \emph{explored}, as well as the vertices~$v_s$ of our sequence located in the interior of $\bdSigma_T$. Indeed, paths emerging from the latter cannot overcome the right boundary of the submap, so their clusters are entirely contained in $\bdSigma_T$. We thus know everything about them. See Figure~\ref{fig:examplerefreshvertexgeneralizedpeelingexplo} for an illustration. 

We do not stop the process however. We restart at the leftmost vertex of our sequence, say~$v_r$, belonging to the complement submap of~$\bdSigma_T$. Within the latter, we faithfully repeat the procedure described above, and carry on as long as there still exists \emph{unexplored} clusters. 
\end{paragraph}

\begin{paragraph}{A Markovian process}
Vertices of our sequence are of two kinds. Either the peeling exploration refreshes at them---like~$v_r$ in Figure~\ref{fig:examplerefreshvertexgeneralizedpeelingexplo}, or they are skipped---like $v_s$. We name the former \textit{stopover vertices}. 
The set of stopover vertices is not deterministic, only~$v_1$ is always one of them. Let~$\bdcalR$ be their number. For any~$1 \leq r \leq \bdcalR$, write~$\bdcalC^*_r$ for the part of the cluster of the $r$-th stopover vertex, included in the current complement submap---that is the complement submap at the~$r$-th refreshment time.
A straightforward consequence of the spatial Markov property satisfied by SCTHP--- see Section~\ref{sec:peelingexplorationSCTHPsmp}---is that for any $1 \leq r \leq N_n$:
$$\bdcalC^*_r \overset{(d)}{\sim} \tildebdcalC \quad \text{conditionally on the event $\bdcalR \geq r$},$$
and $\bdcalC^*_r$ is independent of the part of $\tildebdcalC_{N_n}$ explored before the $r$-th refreshment time. In particular, the distribution of $\bdcalC^*_r$ does not depend on the exact vertex where this refreshment occurs.
\end{paragraph}

\begin{paragraph}{How the mass is spread across $\tildebdcalC_{N_n}$}
We claim that under conditions set in Theorem~\ref{theorem:largecriticalmultipleclusters}, only one cluster, emerging from a stopover vertex, has a macroscopic size far outstripping the others:

\begin{proposition}
\label{prop:spreadmasslargeunioncriticalpercoclusters}
\NoEndMark
There exists a positive valued function $\varphi$ satisfying $\varphi(n)=o\big( n^{1/2}\big)$ such that :
\begin{align}
\widetilde{\bbP}_{p_c}\bigg(\bigcup_{r=1}^{N_n} A_r \cap B_r^c \bigg)=o\big(n^{-1/2-\epsilon}\big) \quad  \text{ and } \quad \widetilde{\bbP}_{p_c}\big(\lvert \tildebdcalC_{N_n} \rvert  \geq n\big) \underset{n \to \infty}{\sim} \widetilde{\bbP}_{p_c}\bigg(\bigcup_{r=1}^{N_n} A_r \cap B_r \bigg), \nonumber
\end{align}
for some~$\epsilon \in (0,\frac{1}{20})$ and where~$A_r$ and~$B_r$ are events defined for any~$1 \leq r \leq N_n$ as:
$$A_r:=(\bdcalR \geq r) \cap \big(\lvert \bdcalC_r^* \rvert \geq n\big) \quad \text{and} \quad B_r:=(\bdcalR \geq r) \cap \big(\lvert \tildebdcalC_{N_n} \setminus \bdcalC_r^* \rvert < \varphi(n)\big).$$
\end{proposition}

The first asymptotic comparison exactly means that the event to see two clusters getting a significant share of the total mass of~$\tildebdcalC_{N_n}$ has a negligible probability compared to that of~$(\lvert \tildebdcalC_{N_n} \rvert\geq~n)$, which is a~$\Omega\big(n^{-1/2}\big)$ as we have already argued. The second asymptotic is merely a consequence of the first. It says that a unique cluster then emerges at large scale, as heralded. The proof of the above proposition is postponed to the end of the current section.
Proposition~\ref{prop:spreadmasslargeunioncriticalpercoclusters} is all we need to prove both items of Theorem~\ref{theorem:largecriticalmultipleclusters}.  We deal with them in separate paragraphs.
\end{paragraph}

\begin{paragraph}{The asymptotic of $\widetilde{\bbP}_{p_c}(\lvert \tildebdcalC_{N_n} \rvert  \geq n)$}
As a direct corollary of Proposition~\ref{prop:spreadmasslargeunioncriticalpercoclusters}:
\begin{align}
\label{eq:secondversionasymptoticspreadmass}
\widetilde{\bbP}_{p_c}\big(\lvert \tildebdcalC_{N_n} \rvert  \geq n\big) \underset{n \to \infty}{\sim} \widetilde{\bbP}_{p_c}\bigg(\bigcup_{r=1}^{N_n} A_r \bigg).
\end{align}
Since~$\bigcup_{r=1}^{N_n} A_r$ can be written as the disjoint union of the events~$\displaystyle{A_r \cap \big(\forall r'<r: \ \lvert \bdcalC_{r'}^*\rvert < n\big)}$, we derive from the spatial Markov property that:
\begin{align}
\label{eq:spatialmarkovpropertylargeunioncriticalpercoclusters}
\widetilde{\bbP}_{p_c}\bigg(\bigcup_{r=1}^{N_n} A_r \bigg)&=\sum_{r=1}^{N_n} \widetilde{\bbP}_{p_c}\bigg(A_r \cap \big(\forall r'<r: \ \lvert \bdcalC_{r'}^*\rvert < n-\varphi(n)\big)\bigg) \nonumber \\ &= \widetilde{\bbP}_{p_c}\big(\lvert \tildebdcalC \rvert \geq n\big) \cdot \sum_{r=1}^{N_n} \widetilde{\bbP}_{p_c}\bigg((\bdcalR \geq r) \cap \big(\forall r'<r: \ \lvert \bdcalC_{r'}^*\rvert < n\big)\bigg).
\end{align}
From this equality, we deduce the upper bound $\widetilde{\bbP}_{p_c}\bigg(\bigcup_{r=1}^{N_n} A_r \bigg) \leq N_n \cdot \widetilde{\bbP}_{p_c}\big(\lvert \tildebdcalC \rvert \geq n\big).$
Theorem~\ref{theorem:criticalexponentstheoremSCTHP}(ii) added to~\eqref{eq:secondversionasymptoticspreadmass} yield the conclusion. \qedd
\end{paragraph}

\begin{remark}
\label{remark:probalargemulticlustersaroundn}
\NoEndMark
The asymptotic comparison~\eqref{eq:secondversionasymptoticspreadmass} and the equality~\eqref{eq:spatialmarkovpropertylargeunioncriticalpercoclusters} together ensure that for any function~$\psi$ such that~$\psi(n)=o(n)$:
\begin{align}
1 \leq \frac{\widetilde{\bbP}_{p_c}\Big(\lvert \tildebdcalC_{N_n} \rvert  \geq n\Big)}{\widetilde{\bbP}_{p_c}\Big(\lvert \tildebdcalC_{N_n} \rvert  \geq n+\psi(n)\Big)} \leq \frac{\widetilde{\bbP}_{p_c}\Big(\lvert \tildebdcalC \rvert  \geq n\Big)}{\widetilde{\bbP}_{p_c}\Big(\lvert \tildebdcalC \rvert  \geq n+\psi(n)\Big)} \underset{\text{Theorem~\ref{theorem:criticalexponentstheoremSCTHP}(ii)}}{=} 1+o(1). \nonumber
\end{align}
Thus, we have~$\widetilde{\bbP}_{p_c}\Big(n \leq \lvert \tildebdcalC_{N_n} \rvert < n+\psi(n)\Big) = o\Big( \widetilde{\bbP}_{p_c}\Big(\lvert \tildebdcalC_{N_n} \rvert \geq n \Big) \Big)$ as $n \to +\infty$.
\end{remark}

\begin{paragraph}{The scaling limit of $\tildebdcalC_{N_n}$}
Let $F$ be a measurable and bounded function on the set of compact metric spaces, which is moreover uniformly continuous with respect to the Gromov--Hausdorff distance. We get from Proposition~\ref{prop:spreadmasslargeunioncriticalpercoclusters} that:
\begin{align}
\label{eq:firstestimatescalinglimitunionclusters}
\widetilde{\bbE}_{p_c} \big[ F(n^{-\frac{1}{2}} \cdot \tildebdcalC_{N_n}) \cdot \mathbf{1}_{\lvert \tildebdcalC_{N_n} \rvert \geq n} \big] = \widetilde{\bbE}_{p_c} \big[ F(n^{-\frac{1}{2}} \cdot \tildebdcalC_{N_n}) \cdot \mathbf{1}_{\cup_{r=1}^{N_n} A_r \cap B_r} \big]+o\Big( \widetilde{\bbP}_{p_c}\big(\lvert \tildebdcalC_{N_n} \rvert \geq n \big) \Big),
\end{align}
by using the fact that $F$ is bounded. Since~$\varphi(n)=o(n^{1/2})$, the events~$A_r \cap B_r$ are disjoint for large~$n$. We then obtain:
\begin{align}
\widetilde{\bbE}_{p_c} \big[ F(n^{-\frac{1}{2}} \cdot \tildebdcalC_{N_n}) \cdot \mathbf{1}_{\cup_{r=1}^{N_n} A_r \cap B_r} \big] = \sum_{r=1}^{N_n} \widetilde{\bbE}_{p_c} \big[ F(n^{-\frac{1}{2}} \cdot \tildebdcalC_{N_n}) \cdot \mathbf{1}_{A_r \cap B_r} \big]. \nonumber
\end{align}
Given the definition of the shortest-path distance $\bfd$, the Gromov--Hausdorff distance between~$\tildebdcalC_{N_n}$ and~$\bdcalC_r^*$, which is a subset of the former if it exists, can be upper bounded as follows:
$$\bfd_{GH}\Big(\tildebdcalC_{N_n},\bdcalC_r^*\Big) \leq \bfd_{H}\Big(\tildebdcalC_{N_n},\bdcalC_r^*\Big) \leq \lvert \tildebdcalC_{N_n} \setminus \bdcalC_r^* \rvert,$$
where $\bfd_H$ is the Hausdorff distance between compact sets of $\tildebdcalC_{N_n}$. 
It implies that on the event~$A_r \cap~B_r$:
\begin{align}
\bfd_{GH}\Big(n^{-\frac{1}{2}} \cdot \tildebdcalC_{N_n}, n^{-\frac{1}{2}} \cdot \bdcalC_r^* \Big) \leq n^{-\frac{1}{2}} \cdot \varphi(n), \nonumber
\end{align}
for any $1 \leq~r \leq~N_n$.
If we denote by $\bdomega(\cdot)$ the modulus of continuity of $F$, we deduce that:
\begin{align}
\bigg\lvert \widetilde{\bbE}_{p_c} \big[ F(n^{-\frac{1}{2}} \cdot \tildebdcalC_{N_n}) \cdot \mathbf{1}_{A_r \cap B_r} \big] - \widetilde{\bbE}_{p_c} \big[ F(n^{-\frac{1}{2}} \cdot \bdcalC_{r}^*) \cdot \mathbf{1}_{A_r \cap B_r} \big] \bigg\rvert \leq \bdomega\bigg(\frac{\varphi(n)}{n^{1/2}}\bigg) \cdot \widetilde{\bbP}_{p_c}(A_r \cap B_r), \nonumber
\end{align}
again for any $1 \leq r \leq N_n$.
By using now Proposition~\ref{prop:spreadmasslargeunioncriticalpercoclusters} and the boundedness of $F$, we have: $$\widetilde{\bbE}_{p_c} \big[ F(n^{-\frac{1}{2}} \cdot \bdcalC_{r}^*) \cdot \mathbf{1}_{A_r} \big]-\widetilde{\bbE}_{p_c} \big[ F(n^{-\frac{1}{2}} \cdot \bdcalC_{r}^*) \cdot \mathbf{1}_{A_r \cap B_r} \big]=o\big(n^{-1/2-\epsilon}\big),$$
uniformly with respect to $r$. Therefore:
\begin{align}
\label{eq:secondestimatescalinglimitunionclusters}
\widetilde{\bbE}_{p_c} \big[ F(n^{-\frac{1}{2}} \cdot \tildebdcalC_{N_n}) \cdot \mathbf{1}_{\cup_{r=1}^{N_n} A_r \cap B_r} \big]= \sum_{r=1}^{N_n} \widetilde{\bbE}_{p_c} \big[ F(n^{-\frac{1}{2}} \cdot \bdcalC_{r}^*) \cdot \mathbf{1}_{A_r} \big] + o\Big( \widetilde{\bbP}_{p_c}\big(\lvert \tildebdcalC_{N_n} \rvert \geq n \big) \Big).
\end{align}
Remark indeed that~$\lim\limits_{n \to +\infty} \bdomega(\frac{\varphi(n)}{n^{1/2}})=0$ because~$F$ is uniformly continuous and~$\varphi(n)=o\big(n^{1/2}\big)$. We state now that:
\begin{align}
\label{eq:claimthereisnottwomacroclusters}
\sum_{r=1}^{N_n} \widetilde{\bbP}_{p_c}\bigg(\big(\exists r' < r  :  \lvert \bdcalC_{r'}^*\rvert \geq n \big) \cap A_r \bigg) = \calO\big(n^{-1+2\epsilon}\big),
\end{align}
which is a mere consequence of the spatial Markov property, combined with an union bound argument and Theorem~\ref{theorem:criticalexponentstheoremSCTHP}(ii). 
The function $F$ is bounded and $\epsilon<\frac{1}{20}$, so we can change~\eqref{eq:secondestimatescalinglimitunionclusters} into:
\begin{align}
\widetilde{\bbE}_{p_c} \big[ F(n^{-\frac{1}{2}} \cdot \tildebdcalC_{N_n}) \cdot \mathbf{1}_{\cup_{r=1}^{N_n} A_r \cap B_r} \big] = \sum_{r=1}^{N_n} \widetilde{\bbE}_{p_c} \bigg[ F(n^{-\frac{1}{2}} \cdot \bdcalC_{r}^*) \cdot \mathbf{1}_{A_r \cap (\forall r'<r: \ \lvert \bdcalC_{r'}^*\rvert < n)} \bigg]+o\Big( \widetilde{\bbP}_{p_c}\big(\lvert \tildebdcalC_{N_n} \rvert \geq n \big) \Big). \nonumber
\end{align}
We apply the spatial Markov property like for \eqref{eq:spatialmarkovpropertylargeunioncriticalpercoclusters} to derive that:
\begin{align}
&\sum_{r=1}^{N_n} \widetilde{\bbE}_{p_c} \bigg[ F(n^{-1/2} \cdot \bdcalC_{r}^*) \cdot \mathbf{1}_{A_r \cap (\forall r'<r: \ \lvert \bdcalC_{r'}^*\rvert < n)} \bigg] \nonumber \\ &= \widetilde{\bbE}_{p_c} \bigg[ F(n^{-1/2} \cdot \tildebdcalC) \cdot \mathbf{1}_{\lvert \tildebdcalC\rvert \geq n} \bigg] \cdot \sum_{r=1}^{N_n} \widetilde{\bbP}_{p_c}\bigg((\bdcalR \geq r) \cap \big(\forall r'<r: \ \lvert \bdcalC_{r'}^*\rvert < n\big)\bigg). \nonumber
\end{align}
The asymptotic~\eqref{eq:firstestimatescalinglimitunionclusters} added to~\eqref{eq:secondversionasymptoticspreadmass} and~\eqref{eq:spatialmarkovpropertylargeunioncriticalpercoclusters}, lead to:
\begin{align}
\widetilde{\bbE}_{p_c} \big[ F(n^{-1/2} \cdot \tildebdcalC_{N_n}) \ \big\vert \ \lvert \tildebdcalC_{N_n} \rvert \geq n \big] =&\ \widetilde{\bbE}_{p_c} \big[ F(n^{-1/2} \cdot \tildebdcalC) \ \big\vert \ \lvert \tildebdcalC \rvert \geq n \big]+o\big(1\big). \nonumber
\end{align}
The scaling limit in Theorem~\ref{theorem:largecriticalmultipleclusters} finally comes from Theorem~\ref{theorem:scalinglimittoymodel}.
\end{paragraph}

\begin{prove}[Proof of~Proposition~\ref{prop:spreadmasslargeunioncriticalpercoclusters}]
We complete the work by proving the core proposition.  

Fix~$\epsilon \in (0,1/20)$ so that~$N_n =~\calO\big( n^{\epsilon}\big)$. We begin with a preliminary statement on the maximal time spent between two refreshment events during the peeling exploration. We write~$T_r$ for the length of the $r$-th exploration. By the spatial Markov property, conditionally on~$(\bdcalR \geq r)$, the random variable~$T_r$ is distributed as the hitting time~$T$ of Proposition~\ref{prop:taildistributionlengthnnexcursion}. Our claim is:
\begin{align}
\label{eq:claimlengthmultiexcursion}
\widetilde{\bbP}_{p_c}\bigg(\max_{1 \leq r \leq \bdcalR} T_r > n^{1+7\epsilon}\bigg) = o(n^{-1/2-2\epsilon}).
\end{align}
It stems directly from an union bound argument.
With the help of~\eqref{eq:claimlengthmultiexcursion}, we gradually rule out the negligible contributions to the mass of $\tildebdcalC_{N_n}$. 

We deal first with that of vertices which are not stopover. We mean the quantity of vertices of~$\tildebdcalC_{N_n}$ only accessible from the latter. Observe on Figure~\ref{fig:examplerefreshvertexgeneralizedpeelingexplo} that they are all located in the interior of bottom "red areas", those revealed right before a refreshment time. Denote these sets by~$\bdfrA_r$ for~$1 \leq r \leq \bdcalR$. Conditionally on the event~$(\bdcalR \geq r)$, by the spatial Markov property, the volume of~$\bdfrA_r$ is stochastically dominated by that of a subcritical Galton--Watson tree~$\frt_r \overset{(d)}{\sim} \frt \overset{(d)}{\sim} \mathbf{GW}_{1-\alpha}$, whose distribution is exponentially-tailed. See the analysis of the peeling exploration algorithm in Section~\ref{sec:explorationclusterRW}. Fix now~$1 \leq r \leq N_n$ and~$t \leq n^{1+7\epsilon}$. We have that:
\begin{align}
&\widetilde{\bbP}_{p_c}\Big(\lvert \bdfrA_r \rvert > n^\delta; \ T_r = t \ \Big\vert \ \bdcalR \geq r \Big) \nonumber \\ &\leq \widetilde{\bbP}_{p_c}\Big(\lvert \frt_r \rvert > n^\delta; \ T_r = t \ \Big\vert \ \bdcalR \geq r \Big) \leq \widetilde{\bbP}_{p_c}\Big(\lvert \frt \rvert > n^\delta \ \Big\vert \ \bdcalR \geq r \Big) = oe_{\delta}(n), \nonumber
\end{align}
for any $\delta>0$, where~$oe$ is the notation for asymptotics introduced in Definition~\ref{def:notasymptoe}. We sum over~$t \leq n^{1+7\epsilon}$ and use \eqref{eq:claimlengthmultiexcursion} to deduce that:
\begin{align}
\widetilde{\bbP}_{p_c}\Big(\lvert \bdfrA_r \rvert > n^{\delta} \ \Big\vert \ \bdcalR \geq r\Big)=o(n^{-1/2-2\epsilon}). \nonumber
\end{align}
This holds again for any $\delta>0$. Finally, by an union bound argument, we get that:
\begin{align}
\label{eq:deviationsnoncoreverticescontrib}
\forall \delta >0, \quad \widetilde{\bbP}_{p_c}\bigg(\sum_{r=1}^{\bdcalR} \lvert \bdfrA_r \rvert > N_n \cdot n^\delta \bigg)=o(n^{-1/2-\epsilon}).
\end{align} 

A roughly similar demonstration works for vertices "left-behind" by the peeling exploration, those in the sets~$\displaystyle{\big(\tildebdcalC_{\bfv_r^*} \setminus \bdcalC_r^* \big)_r}$, where~$\bfv_r^*$ is the $r$-th stopover vertex met. This time, they might be found in several distinct "red areas" that have been revealed throughout the exploration of~$\bdcalC_{r-1}^*$.  On Figure~\ref{fig:examplerefreshvertexgeneralizedpeelingexplo}, we remark that the latter necessarily border the left boundary of the complement submap, as it is after the~$r$-th refreshment time. Their number is consequently bounded by the amount of vertices of~$\bdcalC_r^*$, visited during the $r$-th exploration and located on the left boundary. Such quantity is geometrically distributed. Indeed, if at some point of the exploration of~$\bdcalC_r^*$, an upward revelation occurs and an infinite (blue) part of the SCTHP is unveiled---because the underlying ascending tree is infinite, then we are sure that we will never meet again the left boundary in the future. This event has of course a non-zero chance to occur since ascending trees are supercritical. Thus, for any~$\delta>0$, the number of "red areas" containing vertices of~$\tildebdcalC_{\bfv_r^*} \setminus \bdcalC_r^*$ is lower than~$n^{\delta/2}$ with probability~$\displaystyle{1-oe_{\delta/2}(n)}$ as $n \to +\infty$. Furthermore, the total number of "red areas" generated throughout the $(r-1)$-th exploration is bounded by~$T_{r-1}$, so by~$n^{1+7 \epsilon}$ with probability~$\displaystyle{1-o(n^{-1/2 - 2 \epsilon})}$ according to~\eqref{eq:claimlengthmultiexcursion}. Large deviations theory ensures that the probability for the total volum of~$n^{\delta/2}$ i.i.d.~subcritical Galton--Watson trees to be greater than~$n^{\delta}$ is~$oe_{\delta}(n)$. Since there are at most~$n^{(1+7\epsilon)n^{\delta/2}}=e^{(1+7\epsilon)\cdot \log{(n)}\cdot n^{\delta/2}}$ ways to choose~$n^{\delta/2}$ "red areas" among~$n^{1+7\epsilon}$, we derive from an union bound argument that:
\begin{align}
\forall \delta >0, \quad \widetilde{\bbP}_{p_c}\Big(\lvert \tildebdcalC_{\bfv_r^*} \setminus \bdcalC_r^* \rvert > n^{\delta}; \ T_r \leq n^{1+7\epsilon}\ \Big\vert \ \bdcalR \geq r \Big) = oe_{\delta}(n). \nonumber
\end{align}
And as for~\eqref{eq:deviationsnoncoreverticescontrib}, it implies that:
\begin{align}
\label{eq:deviationscoreverticesleftbehindcontrib}
\forall \delta >0, \quad \widetilde{\bbP}_{p_c}\bigg(\sum_{r=1}^{\bdcalR} \lvert \tildebdcalC_{\bfv_r^*} \setminus \bdcalC_r^* \rvert > N_n \cdot n^\delta \bigg)=o(n^{-1/2-\epsilon}).
\end{align}

Henceforth $\delta=9 \epsilon$. With probability~$1-o(n^{-\frac{1}{2}-\epsilon})$, the mass of~$\tildebdcalC_{N_n} \setminus \bigcup_{r=1}^{\bdcalR} \bdcalC_r^*$ 
does not exceed~$2 N_n n^{9 \epsilon}$, according to~\eqref{eq:deviationsnoncoreverticescontrib} and \eqref{eq:deviationscoreverticesleftbehindcontrib}. Meanwhile, the spatial Markov property, an union bound argument and Theorem~\ref{theorem:criticalexponentstheoremSCTHP}(ii) give us that:
\begin{align}
&\widetilde{\bbP}_{p_c}\Big(\exists 1 \leq r \neq r' \leq \bdcalR \text{ s.t. } \lvert \bdcalC_r^* \rvert \geq n-4 N_n \cdot n^{9 \epsilon} \text{ and } \lvert \bdcalC_{r'}^* \rvert \geq n^{9 \epsilon} \Big) \nonumber \\ &= \calO\big( N_n^2 \cdot n^{-1/2 - 9\epsilon/2} \big) = \calO\big( n^{-1/2 -5 \epsilon /2}\big) = o(n^{-1/2 - \epsilon}). \nonumber
\end{align}
Hence:
\begin{align}
\label{eq:probatohavetwolargecoreclusters}
&\widetilde{\bbP}_{p_c}\bigg(\exists 1 \leq r \leq \bdcalR \text{ s.t. } \lvert \bdcalC_r^* \rvert \geq n-4 N_n \cdot n^{9 \epsilon} \text{ and } \big\lvert \bigcup_{r' \neq r} \bdcalC_{r'}^* \big\rvert \geq N_n \cdot n^{9 \epsilon} \bigg)  = o(n^{-1/2 - \epsilon}). 
\end{align}
Remark that~$N_n \cdot n^{9 \epsilon}=\calO(n^{10 \epsilon})=o(n^{1/2})$. The first asymptotic of Proposition~\ref{prop:spreadmasslargeunioncriticalpercoclusters} is now at hand, by just setting~$\varphi(n)=~4 N_n \cdot~n^{9 \epsilon}$, then using~\eqref{eq:deviationsnoncoreverticescontrib}, \eqref{eq:deviationscoreverticesleftbehindcontrib} and \eqref{eq:probatohavetwolargecoreclusters}. 

We turn our attention to the second asymptotic. As soon as~$\lvert \tildebdcalC_{N_n} \rvert \geq n$, and at the same time~$\big\lvert \tildebdcalC_{N_n} \setminus \bigcup_{r=1}^{\bdcalR} \bdcalC_r^* \big\rvert \leq~2 N_n \cdot n^{9 \epsilon}$, we have~$\big\lvert \bigcup_{r=1}^{\bdcalR} \bdcalC_r^* \big\rvert \geq n-2 \cdot N_n \cdot n^{9 \epsilon}$. Therefore:
$$\lvert \bdcalC_r^* \rvert \geq n/N_n - 2 \cdot n^{9 \epsilon} \geq C \cdot n^{1-\epsilon},$$ 
for some~$1 \leq r \leq \bdcalR$ and some~$C>0$. As~\eqref{eq:probatohavetwolargecoreclusters}, we can prove that:
\begin{align}
&\widetilde{\bbP}_{p_c}\bigg(\exists 1 \leq r \leq \bdcalR \text{ s.t. } \lvert \bdcalC_r^* \rvert \geq C \cdot n^{1-\epsilon} \text{ and } \big\lvert \bigcup_{r' \neq r} \bdcalC_{r'}^* \big\rvert \geq N_n \cdot n^{9 \epsilon} \bigg)  = o(n^{-1/2 - \epsilon}). \nonumber
\end{align}
Because $\widetilde{\bbP}_{p_c}\big(\lvert \tildebdcalC_{N_n} \rvert \geq n\big)=\Omega\big(n^{-1/2}\big)$, it provides the following asymptotic:
\begin{align}
\label{eq:firstequivalentproofproponelargecluster}
\widetilde{\bbP}_{p_c}\bigg(\lvert \tildebdcalC_{N_n} \rvert  \geq n; \  \bigcup_{r=1}^{N_n} A'_r \cap B_r \bigg) \underset{n \to +\infty}{\sim} \widetilde{\bbP}_{p_c}\big(\lvert \tildebdcalC_{N_n} \rvert  \geq n \big),
\end{align}
where~$A'_r$ is for any~$r$ the event~$(\bdcalR \geq r) \cap \big(\lvert \bdcalC_r^* \rvert \geq n-\varphi(n)\big)$. 

Given that~$\varphi(n)=o( n )$, the events~$(A_r \cap B_r )_r$ and~$(A'_r \cap B_r )_r$ become disjoint at large $n$. 
Moreover, both~$\widetilde{\bbP}_{p_c}\Big(\bigcup_{r=1}^{N_n} A'_r \Big)$ and~$\widetilde{\bbP}_{p_c}\Big(\bigcup_{r=1}^{N_n} A_r \Big)$ are $\Omega\big(n^{-1/2}\big)$. Added to the first asymptotic of Proposition~\ref{prop:spreadmasslargeunioncriticalpercoclusters}, we derive that:
 \begin{align}
 \label{eq:equivalentproofproponelargecluster}
\widetilde{\bbP}_{p_c}\bigg(\bigcup_{r=1}^{N_n} A'_r \cap B_r \bigg) \underset{n \to \infty}{\sim} \widetilde{\bbP}_{p_c}\bigg(\bigcup_{r=1}^{N_n} A'_r \bigg) \quad \text{and} \quad \widetilde{\bbP}_{p_c}\bigg(\bigcup_{r=1}^{N_n} A_r \cap B_r \bigg) \underset{n \to \infty}{\sim} \widetilde{\bbP}_{p_c}\bigg(\bigcup_{r=1}^{N_n} A_r \bigg).
 \end{align}
Thanks to~\eqref{eq:spatialmarkovpropertylargeunioncriticalpercoclusters} and Theorem~\ref{theorem:criticalexponentstheoremSCTHP}(ii), we know furthermore that~$\widetilde{\bbP}_{p_c}\big(\bigcup_{r=1}^{N_n} A'_r \big) \underset{n \to \infty}{\sim} \widetilde{\bbP}_{p_c}\big(\bigcup_{r=1}^{N_n} A_r \big).$
By \eqref{eq:equivalentproofproponelargecluster}, we deduce that~$\widetilde{\bbP}_{p_c}\big(\bigcup_{r=1}^{N_n} A'_r \cap B_r\big) \underset{n \to \infty}{\sim} \widetilde{\bbP}_{p_c}\big(\bigcup_{r=1}^{N_n} A_r \cap B_r \big).$
 Combined with the inclusion~$\bigcup_{r=1}^{N_n} A_r \cap B_r \subset \big\{ \lvert \tildebdcalC_{N_n} \rvert \geq n\big\}$ and~\eqref{eq:firstequivalentproofproponelargecluster}, it completes the proof of of the proposition. 
\end{prove}

\section{Open problems and future research}\label{sec:openproblems}

There are several natural ways to improve or even extend the work done: 

\begin{paragraph}{Relaxing the assumption on the offspring distribution.}
The main statements of this paper hold for SCT built from a Galton--Watson tree with a geometric offspring distribution. However, we have good reasons to think that all remain true in a wide class of probability measures, presumably those having a finite second moment. By reading carefully through the Section~\ref{sec:directedpercoonSCQ}, where we survey percolation on SCT by relying on a prior analysis of SCTHP, we realize that no argument really depend on the reproduction law, except the large deviations estimates in Lemma~\ref{lemma:findtwowidetrees} on the growth of a Galton--Watson process. The main hurdle to a generalization of our results rather lies in the toy model itself, given the loss of the spatial Markov property when the offspring distribution is not geometric anymore.
\end{paragraph}

\begin{paragraph}{Involving the horizontal edges.}
As soon as the cycles connecting vertices at the same distance from the root vertex are involved in the percolation process---see Figure~\ref{fig:defcausaltriangulations}, they obviously ease the creation of infinite clusters in the map. We expect everything to hold in such context, apart from the exact value of~$p_c$, certainly smaller. We also conjecture the existence of a second non trivial threshold~$p_c < p_u <1$, from which emerges a new percolation regime where only one infinite cluster prevails in the map. This is supported by~\cite[Proposition 1]{benjamini1996percolation} that deals with a graph model geometrically close to SCT. The primary challenge is to redesign the peeling exploration of SCTHP, since the current version does not delimit areas in which may expand infinite paths made of both upward directed edges and horizontal ones. 
\end{paragraph}

\begin{paragraph}{Large critical clusters in other hyperbolic models.}
As said in the introduction, infinitely many coexisting infinite components is guessed to be a specific feature of percolation in hyperbolic environments, already witnessed in deterministic contexts~\cite{benjamini1996percolation,hutchcroft2019percolation}, as well as in random ones~\cite{curien2019peeling, ray2014geometry}. The Brownian continuum random tree as scaling limit of large critical clusters is also believed to be a universal phenomenon occuring in every such model. But very few has been rigorously demonstrated thus far~\cite{chen2017long}. Is our method adaptable?
\end{paragraph}

\bibliographystyle{plain}
\bibliography{bibliopercoSCT}

\end{document}